\documentclass{amsart}
\usepackage{calc}
\usepackage[nomessages]{fp}
\usepackage{tikz,amsmath,amssymb}
\usepackage{mathtools}
\usepackage{tkz-euclide} 
\usetikzlibrary{calc}
\usetikzlibrary{arrows} 
\usetikzlibrary{decorations.markings}
\usepackage[english]{babel}
\usepackage[hidelinks]{hyperref}
\usepackage{amsmath,amsfonts,amssymb}
\usepackage{amsthm}
\usepackage[all]{xy}
\usepackage{ stmaryrd }
\usepackage{ mathrsfs }
\theoremstyle{plain} 

\usepackage{epsfig}
\usepackage{lmodern}
\usepackage{eurosym}
\usepackage{graphicx}
\usepackage{enumitem}
\usepackage{amsbsy}
\usepackage{rotating}
\usepackage[noabbrev,capitalise]{cleveref}

\DeclareSymbolFont{myletters}{OML}{ztmcm}{m}{it}
\DeclareMathSymbol{\uplambda}{\mathord}{myletters}{"15}

\definecolor{midnightblue}{rgb}{0.1, 0.1, 0.44}
\definecolor{plum}{rgb}{0.56, 0.27, 0.52}
\definecolor{Plum}{rgb}{0.56, 0.27, 0.52}
\definecolor{patriarch}{rgb}{0.5, 0.0, 0.5}
\definecolor{darkgreen}{rgb}{0.0, 0.2, 0.13}
\definecolor{darkcerulean}{rgb}{0.03, 0.27, 0.49}
\definecolor{jade}{rgb}{0.0, 0.66, 0.42}
\definecolor{lava}{rgb}{0.81, 0.06, 0.13}

\hypersetup{
	colorlinks=true,
	linkcolor=darkcerulean,
	citecolor=jade,
	filecolor=magenta,      
	urlcolor=cyan,
	pdfpagemode=FullScreen,
}

\newcommand{\rep}{\operatorname{rep}}
\newcommand{\Cat}{\operatorname{Cat}}
\newcommand{\Int}{\operatorname{Int}}
\newcommand{\GenJF}{\operatorname{GenJF}}
\newcommand{\GenRep}{\operatorname{GenRep}}
\newcommand{\JF}{\operatorname{JF}}

\newcommand{\add}{\operatorname{add}}
\newcommand{\Ker}{\operatorname{Ker}}

\newcommand{\Hom}{\operatorname{Hom}}
\newcommand{\Id}{\operatorname{Id}}

\newcommand{\diag}{\operatorname{diag}}
\newcommand{\cdiag}{\pmb{\operatorname{diag}}}
\newcommand{\op}{\operatorname{op}}
\newcommand{\ind}{\operatorname{ind}}

\newcommand{\GK}{\operatorname{\mathsf{GK}}}
\newcommand{\NEnd}{\operatorname{\mathsf{N}End}}
\newcommand{\wt}{\operatorname{\mathsf{wt}}}
\newcommand{\End}{\operatorname{End}}

\newcommand{\vdim}{\operatorname{\pmb{\dim}}}

\newcommand{\Adds}{\operatorname{AddS}}

\newcommand{\Supp}{\operatorname{Supp}}
\newcommand{\rev}{\operatorname{rev}}

\renewcommand{\L}{\mathbf L}
\newcommand{\R}{\mathbf R}
\newcommand{\A}{\mathbf A}
\newcommand{\B}{\mathbf B}
\newcommand{\E}{\mathbf E}

\newcommand{\RSK}{\operatorname{\mathbf{RSK}}}
\newcommand{\GRSK}{\operatorname{\pmb{\mathcal{RSK}}}}

\newcommand{\RPP}{\operatorname{RPP}}

\newcommand{\Fer}{\operatorname{Fer}}
\newcommand{\Hk}{\operatorname{Hk}}
\newcommand{\SSYT}{\operatorname{SSYT}}
\newcommand{\INC}{\operatorname{INC}}
\newcommand{\SYT}{\operatorname{SYT}}
\newcommand{\AR}{\operatorname{AR}}
\newcommand{\Mult}{\operatorname{Mult}}
\newcommand{\oplamb}{\operatorname{\pmb{\uplambda}}}
\newcommand{\opc}{\operatorname{\pmb{\mathsf{c}}}}
\newcommand{\opf}{\operatorname{\pmb{\mathsf{f}}}}
\newcommand{\opE}{\operatorname{\pmb{\mathsf{E}}}}
\newcommand{\crep}{\operatorname{\pmb{rep}}}

\newcommand{\st}{\operatorname{\pmb{\mathcal{ST}}}}
\newcommand{\ldiag}[1]{\pmb{\langle} #1 \pmb{\rangle}}
\newcommand{\LR}[1]{\pmb{[} #1 \pmb{]}}

\author[B.~Dequêne]{Benjamin Dequêne}
\address[B.~Dequêne]{UFR des Sciences, Laboratoire Amiénois de Mathématiques Fondamentales et Appliquées (LAMFA), Université de Picardie Jules Vernes (UPJV) }
\email{benjamin.dequene@u-picardie.fr}

\title[Extended RSK via quiver representations]{An extended generalization of RSK correspondence via $A$ type quiver representations}

\usepackage{thmtools}
\declaretheorem[numberwithin=section,name=Theorem,
refname={Theorem,Theorems},
Refname={Theorem,Theorems}]{theorem}
\declaretheorem[numberlike=theorem,name=Lemma,
refname={Lemma,Lemmas},
Refname={Lemma,Lemmas}]{lemma}
\declaretheorem[numberlike=theorem,name=Proposition,
refname={Proposition,Propositions},
Refname={Proposition,Propositions}]{proposition}

\declaretheorem[numberlike=theorem,name=Corollary,
refname={Corollary,Corollaries},
Refname={Corollary,Corollaries}]{corollary}

\declaretheorem[style=definition,numberlike=theorem,name=Definition,
refname={Definition,Definitions},
Refname={Definition,Definitions}]{definition}

\declaretheorem[style=definition,numberlike=theorem,name=Example,
refname={Example,Examples},
Refname={Example,Examples}]{example}

\declaretheorem[style=remark,numberlike=theorem,name=Remark,
refname={Remark,Remarks},
Refname={Remark,Remarks}]{remark}

\usepackage{todonotes}

\newcommand{\new}[1]{\textit{\textbf{\color{patriarch}{#1}}}}
\newcommand{\llrr}[1]{\llbracket #1 \rrbracket}

\begin{document}
	
	\maketitle
	
	\begin{abstract} 
		Let $\lambda=(\lambda_1 \geqslant \ldots \geqslant \lambda_k > 0)$. For any $c$ Coxeter element of $\mathfrak{S}_{\lambda_1+k-1}$, we construct a bijection from fillings of $\lambda$ to reverse plane partitions. We recover two previous generalizations of the Robinson--Schensted--Knuth correspondence for particular choices of Coxeter element depending on $\lambda$: one based on the work of, among others, Burge, Hillman, Grassl, Knuth, and uniformly presented by Gansner; the other developed by Garver, Partrias, and Thomas, and independently by Dauvergne, called Scrambled RSK.
		
		Our results in this paper develop the combinatorial consequences of our previous work of type $A_{\lambda_1+k}$ quivers.   
	\end{abstract}
	
	\tableofcontents

	\section{Introduction}
	\label{sec:Intro}

This article has a short version in the proceedings of the 36th edition of the Formal Power Series and Algebraic Combinatorics (FPSAC) Conference \cite{DFPSAC24}.

\subsection{RSK and its generalizations}
\label{ss:introRSK}

Let $n \in \mathbb{N}^*$. The \new{Robinson--Schensted correspondence} is a famous one-to-one correspondence from elements of the symmetric group $\mathfrak{S}_n$ to pairs of standard Young tableaux of the same shape and of size $n$. It first arose from the representation theory of the symmetric group, thanks to the work of Robinson \cite{R38}, before receiving a combinatorial realization via Schensted row insertions \cite{S61}. This correspondence was studied for numerous combinatorial consequences, as a combinatorial proof of a representation-theoretic identity involving the dimension of the irreducible representations of $\mathfrak{S}_n$ (see \cref{rem:RSversion}), Viennot's geometric construction \cite{V77}, plactic monoids \cite{LS81,S97}, or the Erdős--Szekeres theorem \cite{ES87}. We refer the reader to \cite{F96,S13} for more details.

The \new{Robinson--Schensted--Knuth (RSK) correspondence}, denoted by $\RSK$, is a generalization of the Robinson--Schensted correspondence, introduced by Knuth \cite{K70}, and presented as a bijection from nonnegative integer matrices to pairs of semi-standard Young tableaux of the same shape. We recover the Robinson-Schensted correspondence by restricting $\RSK$ to permutation matrices. The RSK correspondence extends many properties of the previous correspondence; for instance, its symmetry under transposition of the matrix results in interchanging the tableaux. As a remarkable consequence, we cite the Cauchy identity for symmetric functions (see \cite{St99,F96} for more details), which generalizes the representation-theoretic identity mentioned above. It also has many interpretations in different settings, using deformations and generalizations of this correspondence. We refer the reader to \cite{P01,Kr06,AF22,GRB23}.

In this paper, we focus on two of those generalizations. Gansner introduced the first \cite{Ga81Ma,Ga81Hi}, based on observations from various works of Burge \cite{B72}, Hillman--Grassl \cite{HG76} and Knuth \cite{K70}. Given a fixed nonzero integer partition $\lambda$, via Greene--Kleitman invariants \cite{GK76}, he defines a map,  denoted by $\GRSK_\lambda$, which realizes a bijection from arbitrary fillings of $\lambda$ to reverse plane partitions of $\lambda$. Garver, Patrias, and Thomas give the second one \cite{GPT19}, in terms of quiver representation theory. Independently, Dauvergne \cite{Dauv20}, in a combinatorial setting, introduced it as \emph{``Scrambled RSK"}. Here, we focus on the quiver representation theory point of view. This variant can be introduced as a family of one-to-one correspondences $(\RSK_{m,c})_{m,c}$, parametrized by orientations of an $A_n$ type quiver (seen here as a Coxeter element $c \in \mathfrak{S}_{n+1})$ --- see \cref{ss:typeACox}), and $m \in \{1,\ldots,n\}$, from $m \times (n-m+1)$ integer matrices to reverse plane partitions of $(n-m+1)^m$ (seen as $n$-tuples of integer partitions satisfying storability conditions -- see \cref{sec:stor}).

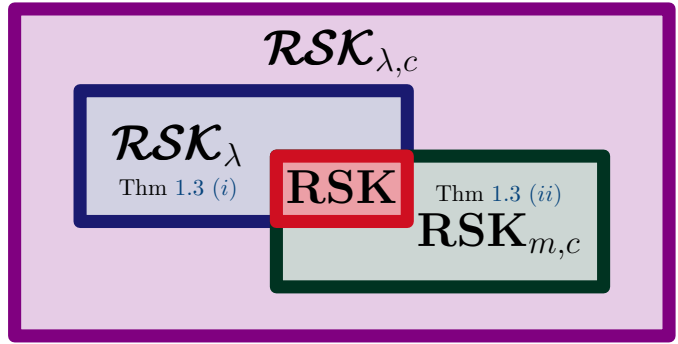
\begin{figure}[!ht]
	\centering
	\scalebox{0.87}{
		\begin{tikzpicture}
			\tkzDefPoint(0,0){a}
			\tkzDefPoint(0,5){b}
			\tkzDefPoint(10,5){c}
			\tkzDefPoint(10,0){d}
			\tkzDrawPolygon[line width = 2mm, color = patriarch, fill=patriarch!20](a,b,c,d);
			
			\tkzLabelPoint[below](5,4.8){{\Huge $\GRSK_{\lambda,c}$}};
			
			\begin{scope}[yshift=-0.25cm]

				\tkzDefPoint(1,2){a}
				\tkzDefPoint(1,4){b}
				\tkzDefPoint(6,4){c}
				\tkzDefPoint(6,2){d}
				\tkzDrawPolygon[line width = 2mm, color = midnightblue,fill=midnightblue!20](a,b,c,d);
				
				\tkzLabelPoint[below](2.5,3.6){{\Huge $\GRSK_\lambda$}};
				
				\tkzLabelPoint[below](2.5,2.8){{Thm \ref{2ndmainresult} \ref{1main}}};
				
				\tkzDefPoint(4,1){a}
				\tkzDefPoint(4,3){b}
				\tkzDefPoint(9,3){c}
				\tkzDefPoint(9,1){d}
				\tkzDrawPolygon[line width = 2mm, color = darkgreen, fill=darkgreen!20](a,b,c,d);
				
				\tkzLabelPoint[above,black](7.4,1.3){{\Huge $\RSK_{m,c}$}};
				
				\tkzLabelPoint[above,black](7.4,2.1){{Thm \ref{2ndmainresult} \ref{2main}}};
				
				\tkzDefPoint(4,2){a}
				\tkzDefPoint(4,3){b}
				\tkzDefPoint(6,3){c}
				\tkzDefPoint(6,2){d}
				\tkzDrawPolygon[line width = 2mm, color = lava, fill=lava!40](a,b,c,d);
				
				\tkzLabelPoint[black](5,2.9){{\Huge $\RSK$}};
			\end{scope}
	\end{tikzpicture}}
	\caption{\label{fig:mainidea} Illustration showing the purpose of the extended generalization $\GRSK_{\lambda,c}$.}
\end{figure}

Our main goal is to exhibit a construction of an extended generalization of $\GRSK_\lambda$, for any nonzero integer partitions $\lambda$, based on a combinatorial extraction of results from \cite{Deq23}, involving any Coxeter element $c \in \mathfrak{S}_{n+1}$, where $n$ is the hook-length of the box $(1,1)$ in $\lambda$, using the combinatorics of the $A_n$ type quivers. We denote those maps by $\GRSK_{\lambda,c}$. In \cref{fig:mainidea}, we pictured how those maps $\GRSK_{\lambda,c}$ can be seen as an extended generalization that contains the previously mentioned correspondence. We state the precise results in \cref{ss:mainres}.

\subsection{Quiver representation theory}
\label{ss:quiverrep}

In this section, we recall the setting of \cite{GPT19} and state the main result of \cite{Deq23}, which motivates our work. 

Fix $\mathbb{K}$ an algebraically closed field, and $n \geqslant 1$. Consider an $A_n$ type quiver $Q$: this is a directed graph whose underlying graph is a line with $n$ vertices. We label the vertices from $1$ to $n$, from left to right. A \emph{representation} $E$ is an assignment of a vector space $E_q$ at each vertex $q$ of $Q$, and an assignment of a linear transformation $E_\alpha$ to each arrow $\alpha$ of $Q$. We say that $E$ is \emph{finite-dimensional} whenever, for all $q \in Q_0$, $E_q$ is finite-dimensional. Given two representations $E$ and $F$, a \emph{morphism} $\phi:E \longrightarrow F$ is a collection of linear maps $(\phi_q : E_q \longrightarrow F_q)_q$ assigned to each vertex of $Q$, such that it satisfies some commutativity properties (see \cref{sec:JRstory}). Denote by $\rep_\mathbb{K}(Q)$ the category of (finite-dimensional) representations of $Q$ over $\mathbb{K}$. One can view this category as a set of representations of $Q$ equipped with morphisms between them. 

A representation is said to be \new{indecomposable} whenever it is not isomorphic to a direct sum of two nonzero representations. Write $\ind_\mathbb{K}(Q)$ for the set of isomorphism classes of indecomposable representations of $Q$. We can encode the data of $\ind_\mathbb{K}(Q)$, and the morphisms between indecomposable representations, by the \new{Auslander--Reiten quiver}, denoted by $\AR_\mathbb{K}(Q)$. It is a directed graph whose vertices are elements of $\ind_\mathbb{K}(Q)$, and arrows are irreducible morphisms between them. We recall that any $E \in \rep_\mathbb{K}(Q)$ is characterized, up to isomorphisms, by the multiplicities of its indecomposable summands. Write $\Mult(E): \ind_\mathbb{K} \longrightarrow \mathbb{N}$ for the map which associates indecomposable representations to their multiplicities. It could be seen as a filling of $\AR_\mathbb{K}(Q)$.

An endomorphism $N: E \longrightarrow E$ is said to be \emph{nilpotent} whenever, for every vertex $q$ of $Q$, $N_q$ is nilpotent. Write $\NEnd(E)$ for the set of nilpotent endomorphisms of a given representation $E$.

We define an invariant on isomorphism classes of $\rep_{\mathbb{K}}(Q)$, called the \new{generic Jordan form data} as follows. Given $E \in \rep_\mathbb{K}(Q)$, we study the set of its nilpotent endomorphisms, denoted by $\NEnd(E)$, by determining their Jordan form: it is displayed as an $n$-tuple of integer partitions. Garver, Patrias, and Thomas \cite{GPT19} proved that a (Zariski) dense open set $\Omega \subset \NEnd(E)$ exists in which all the nilpotent endomorphisms have the same Jordan form. This common Jordan form data is called the generic Jordan form data of $X$, denoted by $\GenJF(E)$ --- see \cref{thm:defGenJF} for the precise statement. Note that this invariant can be computed combinatorially (see \cref{ss:JRandCJR}).

Note that if $n > 1$, $\GenJF$ is not a complete invariant. However, we can still be interested in determining the full subcategories of $\rep_\mathbb{K}(Q)$ (closed under direct sums and summands) in which $\GenJF$ becomes complete. We call these subcategories \new{Jordan recoverable}.

Determining all the Jordan recoverable subcategories of $\rep_\mathbb{K}(Q)$ is still a difficult task. A conjecture is stated in \cite{Deq23} and is recalled in \cref{sec:Further}. Another question raised is how to recover the representation, up to isomorphisms, from its generic Jordan form data. Garver, Patrias, and Thomas described an algebraic way to do so, and they called \new{canonically Jordan recoverable} any subcategory in which their algebraic procedure succeeds. Note that any canonically Jordan recoverable subcategory is Jordan recoverable, but the converse is false, which explains the refined notion.

They prove that, for any vertex $m$ in $Q$, the subcategories additively generated by indecomposable representations $X$ such that $X_m \neq 0$, denoted by $\mathscr{C}_{Q,m}$, are canonically Jordan recoverable. Moreover, they show that $\GenJF$ can be seen as a generalization of the RSK correspondence, as they recover $\RSK$ if $Q$ is oriented such that $m$ is the only sink (respectively, only source) of $Q$. They also showed that they recover the Hillman--Grassl correspondence when $Q$ is linearly oriented (i.e., $ Q$ has a unique source and only one sink). We refer the reader to \cite[Section 6]{GPT19} for more details.

Recall that a filling of $\AR_\mathbb{K}(Q)$ corresponds to a representation $E \in \rep_\mathbb{K}(Q)$ up to isomorphism. Now, see $\GenJF$ as a map from fillings of $\AR_\mathbb{K}(Q)$ (which define, up to isomorphism,  representations of $Q$) to $n$-tuples of integer partitions. This map becomes a bijection if we restrict its domain to fillings $f$ which vanish on indecomposable representations $X$ such that $X_m = 0$, and its codomain to $n$-tuples of integer partitions that satisfy some storability conditions (see \cref{sec:stor} and \cite{Deq23}). In this case, $\GenJF$ coincides with the Dauvergne's Scrambled RSK $\RSK_{m,c}$ mentioned earlier, where $c \in \mathfrak{S}_{n+1}$ is the Coxeter element corresponding to $Q$ (see \cref{ss:mainres}).

The main result of \cite{Deq23} generalizes one of the results of \cite{GPT19} by describing all the canonically Jordan recoverable subcategories of $\rep_\mathbb{K}(Q)$.

Recall that, for any $A_n$ type quiver $Q$, the isomorphism classes of indecomposable representations are in bijection with intervals $\llrr{i,j} = \{i,i+1,\ldots,j\}$ in $\{1,\ldots, n\}$. Moreover, their indecomposable representations characterize subcategories closed under direct sums and summands. Thus, for any subcategory $\mathscr{C}$ of $\rep_\mathbb{K}(Q)$, we write $\Int(\mathscr{C})$ for the set of intervals corresponding to the indecomposable representations in $\mathscr{C}$. 

Two intervals $\llrr{i,j}$ and $\llrr{k,\ell}$ are \new{adjacent} whenever either $j +1 = k$ or $\ell+1 = i$. An interval set $\mathscr{J}$ is said to be \new{adjacency-avoiding} if no pair of intervals in $\mathscr{J}$ are adjacent.

\begin{theorem}[\cite{Deq23}] \label{thm:Deq23}
	Let $n \geqslant 1$ and $Q$ be an $A_n$ type quiver. A subcategory $\mathscr{C}$ is canonically Jordan recoverable if and only if $\Int(\mathscr{C})$ is adjacency-avoiding. 
\end{theorem}

Note that this result shows that canonical Jordan recoverability does not depend on the orientation of $Q$.

A \new{bipartition} of any given set $\A$ is a pair $(\L,\R)$ of subsets of $\A$ such that $\L \cap \R = \varnothing$ and $\L \cup \R = \A$. 

We highlight another remarkable fact from \cite{Deq23}. 
As any subset of an adjacency-avoiding interval set is adjacency-avoiding, we can focus on maximal ones. Those maximal adjacency-avoiding interval sets are parametrized by bipartitions $(\L,\R)$ of $\{2, \ldots, n\}$. Precisely, for any maximal adjacency-avoiding interval set $\mathscr{J}$, there exists a unique bipartition $(\L,\R)$ of $\{2,\ldots,n\}$ such that:
\[\mathscr{J} = \{\llrr{\ell,r-1} \mid \ell \in \L \cup \{1\} \text{ and } r \in \R \cup \{n+1\}\}\]  

\subsection{Main results}
\label{ss:mainres}

We proceed to a combinatorial extraction of the results in \cite{Deq23}. We summarize this extraction in \cref{tab:linksRSKandCJR}.
\begin{table}[!ht]
	\centering
	\scalebox{0.9}{\begin{tabular}{|c|c|}
			\hline \textbf{Combinatorial tools} & \textbf{Representation-theoretic tools} \\
			
			\hline Coxeter element of $\mathfrak{S}_{n+1}$ & Orientation of an $A_{n}$ type quiver $Q$ \\
			
			\hline Transposition in $\mathfrak{S}_{n+1}$ & Indecomposable representation in $\rep_\mathbb{K}(Q)$ \\
			
			\hline AR quiver of $c$ & AR quiver of $\rep_\mathbb{K}(Q)$ \\
			
			\hline {\begin{tabular}{c}
					Integer partition $\lambda$ with hook-length \\ of $(1,1)$ in $\lambda$ equals to $n$
			\end{tabular}} & maximal CJR subcategory $\mathscr{C}$ of $\rep_\mathbb{K}(Q)$ \\
			
			\hline Filling of $\lambda$ &  $\Mult(E)$ for some $E \in \mathscr{C}$ \\
			
			\hline Reverse plane partition of $\lambda$ & $\GenJF(E)$ for some $E \in \mathscr{C}$. \\
			\hline 
	\end{tabular}}
	\vspace*{0.5cm}
	\caption{\label{tab:linksRSKandCJR} Identifications between tools from quiver representation theory and combinatorics.}
\end{table}
The bijective link between integer partitions with $h_\lambda(1,1) = n$ and maximal canonically Jordan recoverable (CJR) subcategories of $\rep_\mathbb{K}(Q)$ is given by the parametrization with bipartitions of $\{2,\ldots,n\}$. Note also that Reading's bijection \cite[Lemma 1.7]{R07} allows us to define a Coxeter element $\opc(\lambda) \in \mathfrak{S}_{n+1}$ from such a $\lambda$.

Given such a $\lambda$, we build a one-to-one correspondence from generic Jordan form data of a representation in the category coming from $\lambda$ to reverse plane partitions of shape $\lambda$, thanks to the notion of $\opc(\lambda)$-storability for $n$-tuples of partitions (see \cref{sec:stor}).

See \cref{sec:Extension} to get the combinatorial construction of $\GRSK_{\lambda,c}$. Our main result is the following.

\begin{theorem} \label{1stmainthm}
	Let $n \geqslant 1$, $\lambda$ be an integer partition such that $h_\lambda(1,1) = n$, and $c \in \mathfrak{S}_{n+1}$ be a Coxeter element. The map $\GRSK_{\lambda,c}$ realizes a one-to-one correspondence from fillings of shape $\lambda$ to reverse plane partitions of shape $\lambda$.  
\end{theorem}

We also show some secondary results that justify the name of ``extended generalization'' of RSK.

\begin{theorem}\label{2ndmainresult}
	Let $n \geqslant 1$, $\lambda$ be an integer partition such that $h_\lambda(1,1) = n$, and $c \in \mathfrak{S}_{n+1}$ be a Coxeter element.
	\begin{enumerate}[label = $(\roman*)$,itemsep=1mm]
		\item \label{1main} If $c = \opc(\lambda)^{\pm 1}$, then $\GRSK_{\lambda,c} = \GRSK_\lambda$.
		\item  \label{2main} If $\lambda=(n-m+1)^m$ for some $m \in \{1,\ldots,n\}$, then $\GRSK_{\lambda,c} = \RSK_{m,c}$.
		\item \label{3main} If $c = (1,\ldots,n+1)$, then $\GRSK_{\lambda,c}$ coincides with the Hillman--Grassl correspondence.
	\end{enumerate}
\end{theorem}

We refer the reader to \cref{sec:Extension} for the proofs of the main theorems, and more details about the local transformation mentioned above.

\subsection{Outline of the paper}
First, in \cref{sec:RSK}, while introducing some essential notations and vocabulary, we recall the classical story of the Robinson--Schensted--Knuth correspondence and discuss a first generalization described by Gansner \cite{Ga81Hi} allowing us to obtain, for any integer partition $\lambda$, a one-to-one correspondence from fillings to reverse plane partitions of shape $\lambda$.

In \cref{sec:Coxeter}, we go through relevant reminders on Coxeter elements of $\mathfrak{S}_{n+1}$. In particular, we highlight their bijective links with bipartitions of $\{2,\ldots,n\}$, integer partitions such that the hook-length of the box $(1,1)$ is $n$, and representations of $A_n$ type quivers.

In \cref{sec:stor}, given a Coxeter element $c \in \mathfrak{S}_{n+1}$, we focus on some behavior for $n$-tuples of integer partitions, called $c$-storability. We enumerate basic combinatorial results, and we introduce diagonal transformations which, for instance, can change $c$-storability into $c'$-storability for some other Coxeter element $c' \in \mathfrak{S}_{n+1}$.

Then, in \cref{sec:JRstory}, we started to gather the pieces of the puzzle with the representation-theoretic notion of Jordan recoverability. After giving some motivations, among others, we recall operations preserving a refinement of this notion, called canonical Jordan recoverability, and a characterization of the maximal canonically Jordan recoverable subcategories for $A_n$ type quivers. This also relates to the $c$-storability properties via the $n$-tuple of integer partitions obtained by computing the generic Jordan form of the nilpotent endomorphism of any given representation in such a category.

In \cref{sec:Extension}, we prove the main result by extracting the results from the previous section into a setting allowing us to see the correspondence highlighted in \cref{tab:linksRSKandCJR}.

Finally, in \cref{sec:Further}, we gave further directions to extend this work: especially a realization of $\GRSK_{\lambda,c}$ via local transformations motivated by \cite{H14,GPT19}, and the elaboration of variations of $\GRSK_{\lambda,c}$ for other Weyl group, or for string algebras.

	\section{The Robinson--Schensted--Knuth correspondence}
	\label{sec:RSK}
	\subsection{Notations and vocabulary}
\label{s:notvoc}
This section introduces the basic objects we need throughout this paper.

\subsubsection*{Quivers and directed graphs} A \new{quiver} is a quadruplet $Q=(Q_0,Q_1,s,t)$ where $Q_0$ is a set called the \emph{vertex set}, $Q_1$ is another set called the \emph{arrow set}, and $s,t : Q_1 \longrightarrow Q_0$ are functions called \emph{source} and \emph{target functions}. Given a quiver $Q$, we denote by $Q^{\op}$ its \new{opposite quiver} defined from $Q$ by reversing all its arrows.

Let $Q=(Q_0,Q_1,s,t)$ and $\Xi = (\Xi_0,\Xi_1,\sigma,\tau)$ be two quivers. A \new{morphism of quivers} $\Psi$ is a pair of maps $(\Psi_0 : Q_0 \longrightarrow \Xi_0, \Psi_1:Q_1\longrightarrow \Xi_1)$ such that, for all $\alpha \in Q_1$, $\sigma(\Xi_1(\alpha)) = \Xi_0 (s(\alpha))$ and $\tau(\Xi_1(\alpha)) = \Xi_0 (t(\alpha))$. Such a morphism $\Psi$ is an \new{isomorphism} whenever $\Psi_0$ and $\Psi_1$ are bijective. We say that $Q$ and $\Xi$ are isomorphic in such a case.

A quiver $Q$ is said to be \new{finite} whenever $Q_0$ and $Q_1$ are finite. We say that $Q$ \new{has no multi-arrows} whenever \[\#\{\alpha \in Q_1 \mid s(\alpha)= q_1 \text{ and } t(\alpha) = q_2\} \leqslant 1\] for all pairs $(q_1,q_2) \in (Q_0)^2$. In the combinatorial setting, we call \new{directed graph} any finite quiver without multi-arrows. Since the arrows in any directed graph are uniquely determined by their source and target, we denote directed graphs by pairs $G=(G_0,G_1)$, where we view the arrow set $G_1$ as a subset of $(G_0)^2$. 

Let $G = (G_0,G_1)$ be a directed graph. A \new{path} $\gamma$ in $G$ as a finite sequence of vertices $(v_0, \ldots, v_k)$ such that $(v_i,v_{i+1}) \in G_1$. A \emph{lazy path at $v \in G_0$} is the path $(v)$. In the following, we denote by $\Pi(G)$ the set of paths in $G$. For any $\gamma = (v_0, \ldots, v_k) \in \Pi(G)$, we denote by $s(\gamma) = v_0$ its source and by $t(\gamma) = v_k$ its target. We also write $\Supp(\gamma) = \{v_0, \ldots, v_k\}$ for the \new{support of $\gamma$}. For $\ell \geqslant 1$, we extend the notion of support to $\ell$-tuples of paths $\pmb{\gamma} = (\gamma_1, \ldots, \gamma_\ell) \in \Pi(G)^\ell$ as $\Supp(\pmb{\gamma}) = \bigcup_{i=1}^\ell \Supp(\gamma_i)$.

A directed graph $G$ is said to be \new{connected} whenever for any pair $(v,w) \in G_0^2$, there exist $\ell \in \mathbb{N}^*$ and $(\gamma_1,\ldots,\gamma_\ell) \in \Pi(G)^\ell$ such that:
\begin{enumerate}[label=$\bullet$,itemsep=1mm]
	\item $v = s(\gamma_1)$, 
	\item for any $i \in \{1,\ldots,\ell-1\}$ odd, $t(\gamma_i) =t(\gamma_{i+1})$, and if $\ell$ is odd, then $t(\gamma_\ell) = w$;
	\item  for any $i \in \{1,\ldots,\ell-1\}$ even, $s(\gamma_i) = s(\gamma_{i+1})$, and if $\ell$ is even, then $s(\gamma_\ell) = w$.
\end{enumerate} 

We say that $G$ is \new{acyclic} if the only paths $\gamma$ in $G$ such that $s(\gamma) = t(\gamma)$ are the lazy ones. Call \new{antichain} of $G$ any subset of vertices $\{w_1,\ldots,w_r\} \subset G_0
$  such that there is no $\gamma \in \Pi(G)$ with $s(\gamma) = w_i $ and $t(\gamma) = w_j$ for all $1 \leqslant i, j \leqslant r$ with $i \neq j$. 

\subsubsection*{Integer partitions}
An \new{integer partition} is a finite weakly decreasing sequence $\lambda = (\lambda_1, \lambda_2, \ldots, \lambda_p)$ of positive integers. Define its \emph{size} as $|\lambda| = \lambda_1 + \ldots + \lambda_k$ and its \emph{length} by $\ell(\lambda) = p$. If needed, we can extend the definition of an integer partition to an infinite weakly decreasing sequence of nonnegative integers with finitely many nonzero entries. 

Given $a,b \in \mathbb{N}^*$, we denote by $a^b$ the integer partition $\lambda$ such that $\ell(\lambda) = b$, and  $\lambda_i = a$ for $1 \leqslant i \leqslant b$.

We endow $(\mathbb{N}^*)^2$ with the \emph{cartesian product order} $\unlhd$ defined by \[(i,j) \unlhd (i',j') \Longleftrightarrow \begin{cases}
	i \leqslant i' \\
	j \leqslant j'
\end{cases}.\]  A \new{Ferrers diagram} is a finite ideal of $((\mathbb{N}^*)^2, \unlhd)$. Recall that we have a one-to-one correspondence between Ferrers diagrams and integer partitions. For a given integer partition $\lambda$, we define the Ferrers diagram of shape $\lambda$ to be \[\Fer(\lambda) = \{(i,j) \in (\mathbb{N}^*)^2 \mid j \leqslant \lambda_i\}.\] Call \emph{box of $\lambda$} any element of $\Fer(\lambda)$. We use matrix coordinates for the boxes of any partition, meaning that we use English conventions to draw Ferrers diagrams.

Given a box $b \in \Fer(\lambda)$, we write $h_\lambda(b)$ for the \new{hook-length of $b$ in $\lambda$}, which is defined, if $b = (i,j)$, as
\[h_\lambda(b) = \# \{(u,v) \in \Fer(\lambda) \mid u \geqslant i, v \geqslant j, \text{ and } u = i \text{ or } v=j \}.\]
Explicitly, one can show that $h_\lambda(b) = \lambda_i - i + \lambda'_j - j +1$ where $\lambda'$ is the \emph{conjugate} of $\lambda$; meaning $\lambda'$ is the unique integer partition such that \[\Fer(\lambda') = \{(j,i) \mid (i,j) \in \Fer(\lambda)\}.\] In particular, we have $h_\lambda(1,1) = \lambda_1 + \ell(\lambda) -1$. In the following, for any $n \in \mathbb{N}$, we write $\Hk_n$ for the set of integer partitions such that $h_\lambda(1,1) = n$.

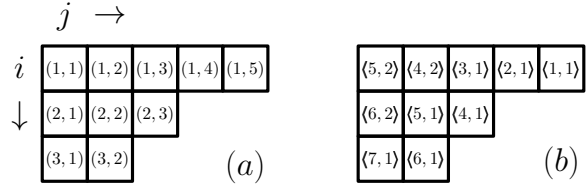
\begin{figure}[!ht]
	\centering
	\scalebox{0.6}{\begin{tikzpicture}
			\begin{scope}[xshift=0cm]
				\tkzLabelPoint[below](-0.5,0.9){{\Huge $i$}};
				\tkzLabelPoint[below](-0.5,-0.1){{\Huge $\downarrow$}};
				\tkzLabelPoint[below](0.5,2.1){{\Huge $j$}};
				\tkzLabelPoint[below](1.5,1.9){{\Huge $\rightarrow$}};
				
				\tkzDefPoint(0,0){a}
				\tkzDefPoint(0,1){b}
				\tkzDefPoint(1,1){c}
				\tkzDefPoint(1,0){d}
				\tkzDrawPolygon[line width = 0.7mm, color = black](a,b,c,d);
				
				\tkzDefPoint(1,0){a}
				\tkzDefPoint(1,1){b}
				\tkzDefPoint(2,1){c}
				\tkzDefPoint(2,0){d}
				\tkzDrawPolygon[line width = 0.7mm, color = black](a,b,c,d);
				
				\tkzDefPoint(2,0){a}
				\tkzDefPoint(2,1){b}
				\tkzDefPoint(3,1){c}
				\tkzDefPoint(3,0){d}
				\tkzDrawPolygon[line width = 0.7mm, color = black](a,b,c,d);
				
				\tkzDefPoint(3,0){a}
				\tkzDefPoint(3,1){b}
				\tkzDefPoint(4,1){c}
				\tkzDefPoint(4,0){d}
				\tkzDrawPolygon[line width = 0.7mm, color = black](a,b,c,d);
				
				\tkzDefPoint(4,1){a}
				\tkzDefPoint(4,0){b}
				\tkzDefPoint(5,0){c}
				\tkzDefPoint(5,1){d}
				\tkzDrawPolygon[line width = 0.7mm, color = black](a,b,c,d);
				
				\tkzDefPoint(0,0){a}
				\tkzDefPoint(0,-1){b}
				\tkzDefPoint(1,-1){c}
				\tkzDefPoint(1,0){d}
				\tkzDrawPolygon[line width = 0.7mm, color = black](a,b,c,d);
				
				\tkzDefPoint(1,0){a}
				\tkzDefPoint(1,-1){b}
				\tkzDefPoint(2,-1){c}
				\tkzDefPoint(2,0){d}
				\tkzDrawPolygon[line width = 0.7mm, color = black](a,b,c,d);
				
				\tkzDefPoint(2,0){a}
				\tkzDefPoint(2,-1){b}
				\tkzDefPoint(3,-1){c}
				\tkzDefPoint(3,0){d}
				\tkzDrawPolygon[line width = 0.7mm, color = black](a,b,c,d);
				
				\tkzDefPoint(0,-2){a}
				\tkzDefPoint(0,-1){b}
				\tkzDefPoint(1,-1){c}
				\tkzDefPoint(1,-2){d}
				\tkzDrawPolygon[line width = 0.7mm, color = black](a,b,c,d);
				
				\tkzDefPoint(1,-2){a}
				\tkzDefPoint(1,-1){b}
				\tkzDefPoint(2,-1){c}
				\tkzDefPoint(2,-2){d}
				\tkzDrawPolygon[line width = 0.7mm, color = black](a,b,c,d);
				
				\tkzLabelPoint(0.5,0.85){{$(1,1)$}};
				\tkzLabelPoint(1.5,0.85){{$(1,2)$}};
				\tkzLabelPoint(2.5,0.85){{$(1,3)$}};
				\tkzLabelPoint(3.5,0.85){{$(1,4)$}};
				\tkzLabelPoint(4.5,0.85){{$(1,5)$}};
				\tkzLabelPoint(0.5,-0.15){{$(2,1)$}};
				\tkzLabelPoint(1.5,-0.15){{$(2,2)$}};
				\tkzLabelPoint(2.5,-0.15){{$(2,3)$}};
				\tkzLabelPoint(0.5,-1.15){{$(3,1)$}};
				\tkzLabelPoint(1.5,-1.15){{$(3,2)$}};
				\tkzLabelPoint(4.5,-1.1){{\Huge $(a)$}};
			\end{scope}
			
			\begin{scope}[xshift=7cm, yshift=0cm]
				\tkzDefPoint(0,0){a}
				\tkzDefPoint(0,1){b}
				\tkzDefPoint(1,1){c}
				\tkzDefPoint(1,0){d}
				\tkzDrawPolygon[line width = 0.7mm, color = black](a,b,c,d);
				
				\tkzDefPoint(1,0){a}
				\tkzDefPoint(1,1){b}
				\tkzDefPoint(2,1){c}
				\tkzDefPoint(2,0){d}
				\tkzDrawPolygon[line width = 0.7mm, color = black](a,b,c,d);
				
				\tkzDefPoint(2,0){a}
				\tkzDefPoint(2,1){b}
				\tkzDefPoint(3,1){c}
				\tkzDefPoint(3,0){d}
				\tkzDrawPolygon[line width = 0.7mm, color = black](a,b,c,d);
				
				\tkzDefPoint(3,0){a}
				\tkzDefPoint(3,1){b}
				\tkzDefPoint(4,1){c}
				\tkzDefPoint(4,0){d}
				\tkzDrawPolygon[line width = 0.7mm, color = black](a,b,c,d);
				
				\tkzDefPoint(4,1){a}
				\tkzDefPoint(4,0){b}
				\tkzDefPoint(5,0){c}
				\tkzDefPoint(5,1){d}
				\tkzDrawPolygon[line width = 0.7mm, color = black](a,b,c,d);
				
				\tkzDefPoint(0,0){a}
				\tkzDefPoint(0,-1){b}
				\tkzDefPoint(1,-1){c}
				\tkzDefPoint(1,0){d}
				\tkzDrawPolygon[line width = 0.7mm, color = black](a,b,c,d);
				
				\tkzDefPoint(1,0){a}
				\tkzDefPoint(1,-1){b}
				\tkzDefPoint(2,-1){c}
				\tkzDefPoint(2,0){d}
				\tkzDrawPolygon[line width = 0.7mm, color = black](a,b,c,d);
				
				\tkzDefPoint(2,0){a}
				\tkzDefPoint(2,-1){b}
				\tkzDefPoint(3,-1){c}
				\tkzDefPoint(3,0){d}
				\tkzDrawPolygon[line width = 0.7mm, color = black](a,b,c,d);
				
				\tkzDefPoint(0,-2){a}
				\tkzDefPoint(0,-1){b}
				\tkzDefPoint(1,-1){c}
				\tkzDefPoint(1,-2){d}
				\tkzDrawPolygon[line width = 0.7mm, color = black](a,b,c,d);
				
				\tkzDefPoint(1,-2){a}
				\tkzDefPoint(1,-1){b}
				\tkzDefPoint(2,-1){c}
				\tkzDefPoint(2,-2){d}
				\tkzDrawPolygon[line width = 0.7mm, color = black](a,b,c,d);
				
				\tkzLabelPoint(0.5,0.85){{$\ldiag{5,2}$}};
				\tkzLabelPoint(1.5,0.85){{$\ldiag{4,2}$}};
				\tkzLabelPoint(2.5,0.85){{$\ldiag{3,1}$}};
				\tkzLabelPoint(3.5,0.85){{$\ldiag{2,1}$}};
				\tkzLabelPoint(4.5,0.85){{$\ldiag{1,1}$}};
				\tkzLabelPoint(0.5,-0.15){{$\ldiag{6,2}$}};
				\tkzLabelPoint(1.5,-0.15){{$\ldiag{5,1}$}};
				\tkzLabelPoint(2.5,-0.15){{$\ldiag{4,1}$}};
				\tkzLabelPoint(0.5,-1.15){{$\ldiag{7,1}$}};
				\tkzLabelPoint(1.5,-1.15){{$\ldiag{6,1}$}};
				\tkzLabelPoint[below](4.5,-1.1){{\Huge $(b)$}};
			\end{scope}
	\end{tikzpicture}}
	\caption{\label{fig:lambdadiag} With $\lambda = (5,3,2)$, $(a)$: the Ferrers diagram of $\lambda$ with the usual coordinates for English convention; $(b):$ the $\lambda$-diagonal coordinates of the boxes in $\Fer(\lambda)$}
\end{figure}

Given a nonzero integer partition $\lambda$, we consider the \emph{$\lambda$-diagonal coordinates} for elements in $\Fer(\lambda)$ as follows. For $k\in \mathbb{Z}$, we define the \new{$k$th diagonal of $\lambda$} as the set $D_k(\lambda)$ of boxes such that $\lambda_1 + i -j = k$. Note that $D_k(\lambda) \neq \varnothing$ if and only if $k \in \{1,\ldots,h_\lambda(1,1)\}$. Assume that $\lambda \in \Hk_n$. For $k \in  \{1,\ldots,n\}$, we set $\delta_k = \max \left(\{\min(i,j) \mid (i,j) \in D_k(\lambda)\} \right)$. We define the \new{$\lambda$-diagonal coordinates} of a box $(i,j) \in \Fer(\lambda)$ to be the pair $\ldiag{k,\delta}_{\lambda}$ where $k \in \{1,\ldots,n\}$ and $\delta \in \{1,\ldots,\delta_k\}$ such that $(i,j) \in D_k(\lambda)$ and $\delta = \delta_k - \min(i,j) + 1$. See \cref{fig:lambdadiag} for an example with $\lambda = (5,3,2)$.
Given $k \in \{1,\ldots,n\}$, we define $\square_k(\lambda)$ the \new{$k$th square of $\lambda$} as the order ideal in $(\Fer(\lambda),\unlhd)$ generated by $D_k(\lambda)$. Note that $D_{\lambda_1}(\lambda)$ corresponds to the \emph{Durfee square of $\lambda$}.

A \new{filling} of shape $\lambda$ is an function $f :\Fer(\lambda) \longrightarrow \mathbb{N}$. Such a filling $f$ is a \new{(weak) reverse plane partition (of shape $\lambda$)} whenever $f$ is weakly increasing with respect to $\unlhd$. We call these reverse plane partitions ``weak'' because we allow $0$ as the value of a box, but we drop this adjective from now on. We denote by $\RPP(\lambda)$ the set of reverse plane partitions of shape $\lambda$.

A reverse plane partition $f$ (of shape $\lambda$) is a \new{(weak) semi-standard Young tableau (of shape $\lambda$)} if $f(i,j) > f(i',j) \geqslant 0$ for any $(i,j), (i',j) \in \Fer(\lambda)$ such that $i' < i$. Write $\SSYT(\lambda)$ for the set of such semi-standard Young tableau of shape $\lambda$, and, given an integer $m \in \mathbb{N}$, $\SSYT(\lambda,m)$ for those with values in $\{0,\ldots,m\}$. 

A filling $f$ is a \new{(weak) increasing tableau (of shape $\lambda$)} whenever $f$ is a reverse plane partition which is strictly increasing with respect to $\unlhd$. We denote by $\INC(\lambda)$ the set of increasing tableaux of shape $\lambda$. 

We define \new{standard Young tableaux (of shape $\lambda$)} as a bijective increasing tableau $f : \Fer(\lambda) \longrightarrow \{1,\ldots,|\lambda|\}$. We write $\SYT(\lambda)$ for the set of standard Young tableaux of shape $\lambda$.

\begin{figure}[!ht]
	\centering
	\scalebox{0.6}{\begin{tikzpicture}
			\begin{scope}[xshift=0cm]
				\tkzDefPoint(0,0){a}
				\tkzDefPoint(0,1){b}
				\tkzDefPoint(1,1){c}
				\tkzDefPoint(1,0){d}
				\tkzDrawPolygon[line width = 0.7mm, color = black](a,b,c,d);
				
				\tkzDefPoint(1,0){a}
				\tkzDefPoint(1,1){b}
				\tkzDefPoint(2,1){c}
				\tkzDefPoint(2,0){d}
				\tkzDrawPolygon[line width = 0.7mm, color = black](a,b,c,d);
				
				\tkzDefPoint(2,0){a}
				\tkzDefPoint(2,1){b}
				\tkzDefPoint(3,1){c}
				\tkzDefPoint(3,0){d}
				\tkzDrawPolygon[line width = 0.7mm, color = black](a,b,c,d);
				
				\tkzDefPoint(3,0){a}
				\tkzDefPoint(3,1){b}
				\tkzDefPoint(4,1){c}
				\tkzDefPoint(4,0){d}
				\tkzDrawPolygon[line width = 0.7mm, color = black](a,b,c,d);
				
				\tkzDefPoint(4,1){a}
				\tkzDefPoint(4,0){b}
				\tkzDefPoint(5,0){c}
				\tkzDefPoint(5,1){d}
				\tkzDrawPolygon[line width = 0.7mm, color = black](a,b,c,d);
				
				\tkzDefPoint(0,0){a}
				\tkzDefPoint(0,-1){b}
				\tkzDefPoint(1,-1){c}
				\tkzDefPoint(1,0){d}
				\tkzDrawPolygon[line width = 0.7mm, color = black](a,b,c,d);
				
				\tkzDefPoint(1,0){a}
				\tkzDefPoint(1,-1){b}
				\tkzDefPoint(2,-1){c}
				\tkzDefPoint(2,0){d}
				\tkzDrawPolygon[line width = 0.7mm, color = black](a,b,c,d);
				
				\tkzDefPoint(2,0){a}
				\tkzDefPoint(2,-1){b}
				\tkzDefPoint(3,-1){c}
				\tkzDefPoint(3,0){d}
				\tkzDrawPolygon[line width = 0.7mm, color = black](a,b,c,d);
				
				\tkzDefPoint(0,-2){a}
				\tkzDefPoint(0,-1){b}
				\tkzDefPoint(1,-1){c}
				\tkzDefPoint(1,-2){d}
				\tkzDrawPolygon[line width = 0.7mm, color = black](a,b,c,d);
				
				\tkzDefPoint(1,-2){a}
				\tkzDefPoint(1,-1){b}
				\tkzDefPoint(2,-1){c}
				\tkzDefPoint(2,-2){d}
				\tkzDrawPolygon[line width = 0.7mm, color = black](a,b,c,d);
				
				\tkzLabelPoint[below](0.5,0.9){{\Huge $0$}};
				\tkzLabelPoint[below](1.5,0.9){{\Huge $0$}};
				\tkzLabelPoint[below](2.5,0.9){{\Huge $2$}};
				\tkzLabelPoint[below](3.5,0.9){{\Huge $4$}};
				\tkzLabelPoint[below](4.5,0.9){{\Huge $6$}};
				\tkzLabelPoint[below](0.5,-0.1){{\Huge $0$}};
				\tkzLabelPoint[below](1.5,-0.1){{\Huge $2$}};
				\tkzLabelPoint[below](2.5,-0.1){{\Huge $2$}};
				\tkzLabelPoint[below](0.5,-1.1){{\Huge $5$}};
				\tkzLabelPoint[below](1.5,-1.1){{\Huge $6$}};
				\tkzLabelPoint[below](4.5,-1.1){{\Huge $(a)$}};
			\end{scope}
			
			\begin{scope}[xshift=7cm, yshift=0cm]
				\tkzDefPoint(0,0){a}
				\tkzDefPoint(0,1){b}
				\tkzDefPoint(1,1){c}
				\tkzDefPoint(1,0){d}
				\tkzDrawPolygon[line width = 0.7mm, color = black](a,b,c,d);
				
				\tkzDefPoint(1,0){a}
				\tkzDefPoint(1,1){b}
				\tkzDefPoint(2,1){c}
				\tkzDefPoint(2,0){d}
				\tkzDrawPolygon[line width = 0.7mm, color = black](a,b,c,d);
				
				\tkzDefPoint(2,0){a}
				\tkzDefPoint(2,1){b}
				\tkzDefPoint(3,1){c}
				\tkzDefPoint(3,0){d}
				\tkzDrawPolygon[line width = 0.7mm, color = black](a,b,c,d);
				
				\tkzDefPoint(3,0){a}
				\tkzDefPoint(3,1){b}
				\tkzDefPoint(4,1){c}
				\tkzDefPoint(4,0){d}
				\tkzDrawPolygon[line width = 0.7mm, color = black](a,b,c,d);
				
				\tkzDefPoint(4,1){a}
				\tkzDefPoint(4,0){b}
				\tkzDefPoint(5,0){c}
				\tkzDefPoint(5,1){d}
				\tkzDrawPolygon[line width = 0.7mm, color = black](a,b,c,d);
				
				\tkzDefPoint(0,0){a}
				\tkzDefPoint(0,-1){b}
				\tkzDefPoint(1,-1){c}
				\tkzDefPoint(1,0){d}
				\tkzDrawPolygon[line width = 0.7mm, color = black](a,b,c,d);
				
				\tkzDefPoint(1,0){a}
				\tkzDefPoint(1,-1){b}
				\tkzDefPoint(2,-1){c}
				\tkzDefPoint(2,0){d}
				\tkzDrawPolygon[line width = 0.7mm, color = black](a,b,c,d);
				
				\tkzDefPoint(2,0){a}
				\tkzDefPoint(2,-1){b}
				\tkzDefPoint(3,-1){c}
				\tkzDefPoint(3,0){d}
				\tkzDrawPolygon[line width = 0.7mm, color = black](a,b,c,d);
				
				\tkzDefPoint(0,-2){a}
				\tkzDefPoint(0,-1){b}
				\tkzDefPoint(1,-1){c}
				\tkzDefPoint(1,-2){d}
				\tkzDrawPolygon[line width = 0.7mm, color = black](a,b,c,d);
				
				\tkzDefPoint(1,-2){a}
				\tkzDefPoint(1,-1){b}
				\tkzDefPoint(2,-1){c}
				\tkzDefPoint(2,-2){d}
				\tkzDrawPolygon[line width = 0.7mm, color = black](a,b,c,d);
				
				\tkzLabelPoint[below](0.5,0.9){{\Huge $0$}};
				\tkzLabelPoint[below](1.5,0.9){{\Huge $0$}};
				\tkzLabelPoint[below](2.5,0.9){{\Huge $2$}};
				\tkzLabelPoint[below](3.5,0.9){{\Huge $4$}};
				\tkzLabelPoint[below](4.5,0.9){{\Huge $6$}};
				\tkzLabelPoint[below](0.5,-0.1){{\Huge $2$}};
				\tkzLabelPoint[below](1.5,-0.1){{\Huge $3$}};
				\tkzLabelPoint[below](2.5,-0.1){{\Huge $3$}};
				\tkzLabelPoint[below](0.5,-1.1){{\Huge $5$}};
				\tkzLabelPoint[below](1.5,-1.1){{\Huge $6$}};
				\tkzLabelPoint[below](4.5,-1.1){{\Huge $(b)$}};
			\end{scope}
			
			\begin{scope}[xshift=0cm,yshift=-4cm]
				\tkzDefPoint(0,0){a}
				\tkzDefPoint(0,1){b}
				\tkzDefPoint(1,1){c}
				\tkzDefPoint(1,0){d}
				\tkzDrawPolygon[line width = 0.7mm, color = black](a,b,c,d);
				
				\tkzDefPoint(1,0){a}
				\tkzDefPoint(1,1){b}
				\tkzDefPoint(2,1){c}
				\tkzDefPoint(2,0){d}
				\tkzDrawPolygon[line width = 0.7mm, color = black](a,b,c,d);
				
				\tkzDefPoint(2,0){a}
				\tkzDefPoint(2,1){b}
				\tkzDefPoint(3,1){c}
				\tkzDefPoint(3,0){d}
				\tkzDrawPolygon[line width = 0.7mm, color = black](a,b,c,d);
				
				\tkzDefPoint(3,0){a}
				\tkzDefPoint(3,1){b}
				\tkzDefPoint(4,1){c}
				\tkzDefPoint(4,0){d}
				\tkzDrawPolygon[line width = 0.7mm, color = black](a,b,c,d);
				
				\tkzDefPoint(4,1){a}
				\tkzDefPoint(4,0){b}
				\tkzDefPoint(5,0){c}
				\tkzDefPoint(5,1){d}
				\tkzDrawPolygon[line width = 0.7mm, color = black](a,b,c,d);
				
				\tkzDefPoint(0,0){a}
				\tkzDefPoint(0,-1){b}
				\tkzDefPoint(1,-1){c}
				\tkzDefPoint(1,0){d}
				\tkzDrawPolygon[line width = 0.7mm, color = black](a,b,c,d);
				
				\tkzDefPoint(1,0){a}
				\tkzDefPoint(1,-1){b}
				\tkzDefPoint(2,-1){c}
				\tkzDefPoint(2,0){d}
				\tkzDrawPolygon[line width = 0.7mm, color = black](a,b,c,d);
				
				\tkzDefPoint(2,0){a}
				\tkzDefPoint(2,-1){b}
				\tkzDefPoint(3,-1){c}
				\tkzDefPoint(3,0){d}
				\tkzDrawPolygon[line width = 0.7mm, color = black](a,b,c,d);
				
				\tkzDefPoint(0,-2){a}
				\tkzDefPoint(0,-1){b}
				\tkzDefPoint(1,-1){c}
				\tkzDefPoint(1,-2){d}
				\tkzDrawPolygon[line width = 0.7mm, color = black](a,b,c,d);
				
				\tkzDefPoint(1,-2){a}
				\tkzDefPoint(1,-1){b}
				\tkzDefPoint(2,-1){c}
				\tkzDefPoint(2,-2){d}
				\tkzDrawPolygon[line width = 0.7mm, color = black](a,b,c,d);
				
				\tkzLabelPoint[below](0.5,0.9){{\Huge $0$}};
				\tkzLabelPoint[below](1.5,0.9){{\Huge $1$}};
				\tkzLabelPoint[below](2.5,0.9){{\Huge $5$}};
				\tkzLabelPoint[below](3.5,0.9){{\Huge $6$}};
				\tkzLabelPoint[below](4.5,0.9){{\Huge $12$}};
				\tkzLabelPoint[below](0.5,-0.1){{\Huge $2$}};
				\tkzLabelPoint[below](1.5,-0.1){{\Huge $4$}};
				\tkzLabelPoint[below](2.5,-0.1){{\Huge $7$}};
				\tkzLabelPoint[below](0.5,-1.1){{\Huge $8$}};
				\tkzLabelPoint[below](1.5,-1.1){{\Huge $10$}};
				\tkzLabelPoint[below](4.5,-1.1){{\Huge $(c)$}};
			\end{scope}
			
			\begin{scope}[xshift=7cm,yshift=-4cm]
				\tkzDefPoint(0,0){a}
				\tkzDefPoint(0,1){b}
				\tkzDefPoint(1,1){c}
				\tkzDefPoint(1,0){d}
				\tkzDrawPolygon[line width = 0.7mm, color = black](a,b,c,d);
				
				\tkzDefPoint(1,0){a}
				\tkzDefPoint(1,1){b}
				\tkzDefPoint(2,1){c}
				\tkzDefPoint(2,0){d}
				\tkzDrawPolygon[line width = 0.7mm, color = black](a,b,c,d);
				
				\tkzDefPoint(2,0){a}
				\tkzDefPoint(2,1){b}
				\tkzDefPoint(3,1){c}
				\tkzDefPoint(3,0){d}
				\tkzDrawPolygon[line width = 0.7mm, color = black](a,b,c,d);
				
				\tkzDefPoint(3,0){a}
				\tkzDefPoint(3,1){b}
				\tkzDefPoint(4,1){c}
				\tkzDefPoint(4,0){d}
				\tkzDrawPolygon[line width = 0.7mm, color = black](a,b,c,d);
				
				\tkzDefPoint(4,1){a}
				\tkzDefPoint(4,0){b}
				\tkzDefPoint(5,0){c}
				\tkzDefPoint(5,1){d}
				\tkzDrawPolygon[line width = 0.7mm, color = black](a,b,c,d);
				
				\tkzDefPoint(0,0){a}
				\tkzDefPoint(0,-1){b}
				\tkzDefPoint(1,-1){c}
				\tkzDefPoint(1,0){d}
				\tkzDrawPolygon[line width = 0.7mm, color = black](a,b,c,d);
				
				\tkzDefPoint(1,0){a}
				\tkzDefPoint(1,-1){b}
				\tkzDefPoint(2,-1){c}
				\tkzDefPoint(2,0){d}
				\tkzDrawPolygon[line width = 0.7mm, color = black](a,b,c,d);
				
				\tkzDefPoint(2,0){a}
				\tkzDefPoint(2,-1){b}
				\tkzDefPoint(3,-1){c}
				\tkzDefPoint(3,0){d}
				\tkzDrawPolygon[line width = 0.7mm, color = black](a,b,c,d);
				
				\tkzDefPoint(0,-2){a}
				\tkzDefPoint(0,-1){b}
				\tkzDefPoint(1,-1){c}
				\tkzDefPoint(1,-2){d}
				\tkzDrawPolygon[line width = 0.7mm, color = black](a,b,c,d);
				
				\tkzDefPoint(1,-2){a}
				\tkzDefPoint(1,-1){b}
				\tkzDefPoint(2,-1){c}
				\tkzDefPoint(2,-2){d}
				\tkzDrawPolygon[line width = 0.7mm, color = black](a,b,c,d);
				
				\tkzLabelPoint[below](0.5,0.9){{\Huge $1$}};
				\tkzLabelPoint[below](1.5,0.9){{\Huge $2$}};
				\tkzLabelPoint[below](2.5,0.9){{\Huge $5$}};
				\tkzLabelPoint[below](3.5,0.9){{\Huge $8$}};
				\tkzLabelPoint[below](4.5,0.9){{\Huge $10$}};
				\tkzLabelPoint[below](0.5,-0.1){{\Huge $3$}};
				\tkzLabelPoint[below](1.5,-0.1){{\Huge $4$}};
				\tkzLabelPoint[below](2.5,-0.1){{\Huge $9$}};
				\tkzLabelPoint[below](0.5,-1.1){{\Huge $6$}};
				\tkzLabelPoint[below](1.5,-1.1){{\Huge $7$}};
				\tkzLabelPoint[below](4.5,-1.1){{\Huge $(d)$}};
			\end{scope}
	\end{tikzpicture}}
	\caption{\label{fig:FerandRPP} With $\lambda = (5,3,2)$, $(a)$: a  (weak) reverse plane partition of shape $\lambda$; $(b):$ a (weak) semi-standard Young tableau of shape $\lambda$; $(c):$ an increasing tableaux of shape $\lambda$; $(d)$ a standard Young tableaux of shape $\lambda$.}
\end{figure}
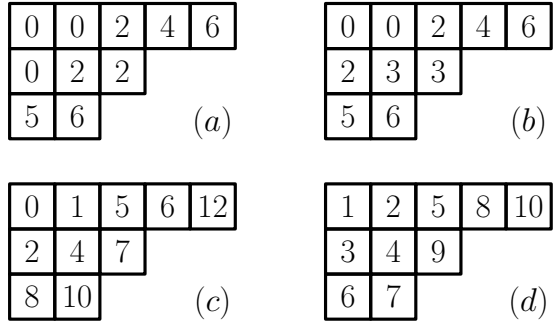

See \cref{fig:FerandRPP} for an example of each of the previous notions for $\lambda = (5,3,2)$.

\subsection{The classical story}
\label{ss:GenRSK} We recall the classical way to present the \emph{Robinson--Schensted--Knuth correspondence}. For more details, we invite the reader to look at the following references: \cite{St99}, \cite{F96}.

Let us first recall the Schensted row-insertion. 

Let $k \in \mathbb{N}$. The \new{Schensted row-insertion of $k$}, denoted $k \overset{\operatorname{Sch.}}{\longrightarrow} -$, is a function on semi-standard Young tableaux defined as follows. Let $f$ be a semi-standard Young tableau of shape $\lambda$. Then the filling $\left(k \overset{\operatorname{Sch.}}{\longrightarrow} f \right) = g$ is obtained thanks to the following procedure:
\begin{enumerate}[label = $\arabic*)$]
	\item Put $i'=1$ and $x = k$;
	\item If it exists, let $j'$ be the smallest index such that $x < f(i',j')$, otherwise, we put $j' = \lambda_{i'} + 1$;
	\item If $(i',j') \notin \Fer(\lambda)$, then we put $g(i',j') = x$ and for $(i,j) \in \Fer(\lambda)$ such that $i \geqslant i'$ we put $g(i,j) = f(i,j)$, and we are done;
	\item Otherwise, we put $g(i',j') = x$ and for $j \neq j'$ such that $(i',j) \in \Fer(\lambda)$, we put $g(i',j) = f(i',j)$; put $x = f(i',j')$, we increase $i'$ by $1$, and we come back to step $2)$.
\end{enumerate}
We end with a filling $g$ of the integer partition obtained from $\lambda$ by adding the box $(i',j')$ from step $3$) of the algorithm. We illustrate how Schensted row insertion works with the following example.
\begin{example} Consider $f$ to be the semi-standard Young tableau of shape $\lambda= (5,2)$ below (\cref{fig:exf}).
	\begin{figure}[!ht]
		\begin{center}
			\scalebox{0.6}{\begin{tikzpicture}
					\tkzDefPoint(0,0){a}
					\tkzDefPoint(0,1){b}
					\tkzDefPoint(1,1){c}
					\tkzDefPoint(1,0){d}
					\tkzDrawPolygon[line width = 0.7mm, color = black](a,b,c,d);
					
					\tkzDefPoint(1,0){a}
					\tkzDefPoint(1,1){b}
					\tkzDefPoint(2,1){c}
					\tkzDefPoint(2,0){d}
					\tkzDrawPolygon[line width = 0.7mm, color = black](a,b,c,d);
					
					\tkzDefPoint(2,0){a}
					\tkzDefPoint(2,1){b}
					\tkzDefPoint(3,1){c}
					\tkzDefPoint(3,0){d}
					\tkzDrawPolygon[line width = 0.7mm, color = black](a,b,c,d);
					
					\tkzDefPoint(3,0){a}
					\tkzDefPoint(3,1){b}
					\tkzDefPoint(4,1){c}
					\tkzDefPoint(4,0){d}
					\tkzDrawPolygon[line width = 0.7mm, color = black](a,b,c,d);
					
					\tkzDefPoint(4,1){a}
					\tkzDefPoint(4,0){b}
					\tkzDefPoint(5,0){c}
					\tkzDefPoint(5,1){d}
					\tkzDrawPolygon[line width = 0.7mm, color = black](a,b,c,d);
					
					\tkzDefPoint(0,0){a}
					\tkzDefPoint(0,-1){b}
					\tkzDefPoint(1,-1){c}
					\tkzDefPoint(1,0){d}
					\tkzDrawPolygon[line width = 0.7mm, color = black](a,b,c,d);
					
					\tkzDefPoint(1,0){a}
					\tkzDefPoint(1,-1){b}
					\tkzDefPoint(2,-1){c}
					\tkzDefPoint(2,0){d}
					\tkzDrawPolygon[line width = 0.7mm, color = black](a,b,c,d);
					
					\tkzLabelPoint[below](0.5,0.9){{\Huge $1$}};
					\tkzLabelPoint[below](1.5,0.9){{\Huge $2$}};
					\tkzLabelPoint[below](2.5,0.9){{\Huge $2$}};
					\tkzLabelPoint[below](3.5,0.9){{\Huge $3$}};
					\tkzLabelPoint[below](4.5,0.9){{\Huge $3$}};
					\tkzLabelPoint[below](0.5,-.1){{\Huge $3$}};
					\tkzLabelPoint[below](1.5,-.1){{\Huge $3$}};
			\end{tikzpicture}}
			\caption{\label{fig:exf} Example of a semi-standard Young tableau $f$.}
		\end{center}
	\end{figure}
	
	We obtain $1 \overset{\operatorname{Sch.}}{\longrightarrow} f$ by replacing the value in the box $(1,2)$ by $1$, the value in the box $(2,1)$ by $2$, by adding a box at $(3,1)$ and giving it the value $3$. The filling obtained following the Schensted row-insertion algorithm is illustrated below (\cref{fig:schrowins}). We greyed the modified boxes and framed the added box.
	\begin{figure}[!ht]
		\begin{center}
			\scalebox{0.6}{\begin{tikzpicture}
					\tkzDefPoint(0,0){a}
					\tkzDefPoint(0,1){b}
					\tkzDefPoint(1,1){c}
					\tkzDefPoint(1,0){d}
					\tkzDrawPolygon[line width = 0.7mm, color = black](a,b,c,d);
					
					\tkzDefPoint(1,0){a}
					\tkzDefPoint(1,1){b}
					\tkzDefPoint(2,1){c}
					\tkzDefPoint(2,0){d}
					\tkzDrawPolygon[line width = 0.7mm, color = black, fill=black!10](a,b,c,d);
					
					\tkzDefPoint(2,0){a}
					\tkzDefPoint(2,1){b}
					\tkzDefPoint(3,1){c}
					\tkzDefPoint(3,0){d}
					\tkzDrawPolygon[line width = 0.7mm, color = black](a,b,c,d);
					
					\tkzDefPoint(3,0){a}
					\tkzDefPoint(3,1){b}
					\tkzDefPoint(4,1){c}
					\tkzDefPoint(4,0){d}
					\tkzDrawPolygon[line width = 0.7mm, color = black](a,b,c,d);
					
					\tkzDefPoint(4,1){a}
					\tkzDefPoint(4,0){b}
					\tkzDefPoint(5,0){c}
					\tkzDefPoint(5,1){d}
					\tkzDrawPolygon[line width = 0.7mm, color = black](a,b,c,d);
					
					\tkzDefPoint(0,0){a}
					\tkzDefPoint(0,-1){b}
					\tkzDefPoint(1,-1){c}
					\tkzDefPoint(1,0){d}
					\tkzDrawPolygon[line width = 0.7mm, color = black,, fill=black!10](a,b,c,d);
					
					\tkzDefPoint(1,0){a}
					\tkzDefPoint(1,-1){b}
					\tkzDefPoint(2,-1){c}
					\tkzDefPoint(2,0){d}
					\tkzDrawPolygon[line width = 0.7mm, color = black](a,b,c,d);
					
					\tkzDefPoint(0,-2){a}
					\tkzDefPoint(0,-1){b}
					\tkzDefPoint(1,-1){c}
					\tkzDefPoint(1,-2){d}
					\tkzDrawPolygon[line width = 2mm, color = black, fill = black!10](a,b,c,d);
					
					\tkzLabelPoint[below](0.5,0.9){{\Huge $1$}};
					\tkzLabelPoint[below](1.5,0.9){{\Huge $\mathbf{1}$}};
					\tkzLabelPoint[below](2.5,0.9){{\Huge $2$}};
					\tkzLabelPoint[below](3.5,0.9){{\Huge $3$}};
					\tkzLabelPoint[below](4.5,0.9){{\Huge $3$}};
					\tkzLabelPoint[below](0.5,-0.1){{\Huge $\mathbf{2}$}};
					\tkzLabelPoint[below](1.5,-0.1){{\Huge $3$}};
					\tkzLabelPoint[below](0.5,-1.1){{\Huge $\mathbf{3}$}};
			\end{tikzpicture}}
			\caption{\label{fig:schrowins} The filling $1 \overset{\operatorname{Sch.}}{\longrightarrow} f$ obtained from $f$ in \cref{fig:exf}}
		\end{center}
	\end{figure}
	Remark that this new filling is a semi-standard Young tableau of shape $\mu = (5,2,1)$.
\end{example}

We can now present the RSK correspondence.

The \new{RSK correspondence} is a map from nonnegative integer matrices and pairs of semi-standard Young tableaux of the same shape, described as follows:
\begin{enumerate}[label=\arabic*)]
	
	\item From $A = (a_{i,j})$, an $n \times m$ nonnegative integer matrix, consider the associated two-line array, \[w_A = \left( \begin{matrix}
		i_1 & i_2 & \ldots & i_s \\
		j_1 & j_2 & \ldots & j_s
	\end{matrix} \right)\] such that, for any $(i,j) \in \{1, \ldots,n\}\times\{1,\ldots,m\}$, there are $a_{i,j}$ copies of the column $\left( \begin{matrix}
		i \\
		j
	\end{matrix}\right)$, and all the columns are in lexicographic order, meaning:
	\begin{enumerate}[label = $\bullet$]
		\item $i_1 \leqslant \ldots \leqslant i_s$, and;
		
		\item if $i_p = i_{p+1}$ then $j_p \leqslant j_{p+1}$.
	\end{enumerate}
	They are usually called \emph{biwords}.
	
	\item We construct two sequences of semi-standard Young tableaux,  $(P(k))_{0 \leqslant k \leqslant s}$ and $(Q(k))_{0 \leqslant k \leqslant s}$, as it follows:
	
	\begin{enumerate}[label = $\bullet$]
		
		\item we begin with $P(0) = Q(0) = \varnothing$;
		
		\item for all $k \in \{1, \ldots,s\}$, we put $P(k) = j_k \overset{\operatorname{Sch.}}{\longrightarrow} P(k-1)$.
		
		\item for all $k \in \{1, \ldots,s\}$, we get $Q(k)$ from $Q(k-1)$ by recording $i_k$ in box created when passing from $P(k-1)$ to $P(k)$.
		
	\end{enumerate}
	
	\item We define \new{$\RSK(A) = (P(s),Q(s))$}.
\end{enumerate} 

\begin{example} \label{ex:RSK1}
	Let us take \[A = \left( 
 \right).\] 
	
	Then we get $$ w_A = \left( %
 \right). $$
	
	In \cref{fig:RSKex1}, we make explicit the step-by-step calculations of $({\color{red}{P(k)}}, {\color{blue}{Q(k)}})_{1 \leqslant k \leqslant 9} $. Here $\RSK(A) = (P(9),Q(9))$. \qedhere
	
	\begin{figure}[!ht]
		\centering
		
		\scalebox{0.5}{
			%
}

		\caption{\label{fig:RSKex1} Illustration of the calculations for $\RSK(A)$. The framed box at the row $k$ is the one added by Schensted row-insertion from ${\color{red}{P(k-1)}}$ to ${\color{red}{P(k)}}$ and the colored boxes are the modified ones from ${\color{red}{P(k-1)}}$ to ${\color{red}{P(k)}}$ following the row-insertion process.}
	\end{figure}
\end{example}

\begin{theorem} \label{th:RSKbij}
	Let $n,m \in \mathbb{N}^*$. The RSK correspondence gives a bijection from $n \times m$ nonnegative integer matrices to pairs $(P,Q)$ of semi-standard Young tableaux of the same shape such that their entries are from $1$ to $m$ for $P$, and from $1$ to $n$ for $Q$. 
\end{theorem}

A well-known combinatorial consequence is the Cauchy identity. Before stating it, we recall Schur polynomials. 

Let $n \in \mathbb{N}^*$, and $x_1,\ldots,x_n$ be $n$ formal variables. Consider $\lambda$ to be a nonzero integer partition. We define the \new{Schur polynomial of $\lambda$} as follows:
\[s_\lambda(x_1,\ldots,x_n) = \sum_{f \in \SSYT(\lambda, n-1)} \prod_{b \in \Fer(\lambda)} x_{f(b)+1}.\] We write $x_{f(b)+1}$ instead of $x_{f(b)}$ because we are considering weak semi-standard Young tableaux. It is well-known that this is a homogeneous symmetric polynomial of degree $|\lambda|$. Moreover, for any $1 \leqslant m \leqslant n$, $(s_\lambda)_{\lambda \vdash m}$ gives a basis of the vector space of the homogeneous symmetric polynomials of degree $m$.

For $n=2$ and $\lambda = (2,1)$ we get $s_\lambda (x_1,x_2) = x_1^2 x_2 +  x_1 x_2^2.$
\begin{corollary}[\textbf{Cauchy identity} (see \cite{St99}, \cite{F96})] For any $n,m \in \mathbb{N}^*$, and for any $x_1,\ldots, x_n$ and $y_1,\ldots,y_m$ sets of formal variables, we have
	\[\sum_{\lambda} s_\lambda(x_1,\ldots,x_n) s_\lambda(y_1,\ldots,y_m) = \prod_{i=1}^n \prod_{j=1}^m \dfrac{1}{1 - x_i y_j}\]
	where the sum is over all the integer partitions $\lambda$.
\end{corollary}

\begin{remark} \label{rem:RSversion}
	The RSK correspondence induces a bijection from permutations of $\mathfrak{S}_n$ to pairs of standard Young tableaux $(P,Q)$ of size $n$. This recovers the Robinson--Schensted correspondence. It allows us to establish the following representation-theoretic identity combinatorially:
	\[ n! = \sum_{\lambda \vdash n} t_\lambda^2, \] where $t_\lambda$ is both the number of standard Young tableaux of shape $\lambda$, with values in $\{1,\ldots, n\}$, and the dimension of the irreducible representation of $\mathfrak{S}_n$ corresponding to the partition $\lambda$.
\end{remark}

\subsection{The Greene--Kleitman invariant}
\label{ss:GKinv}
Let $G=(G_0,G_1)$ be a directed graph. Assume that $G$ is acyclic. Consider a filling $f : G_0 \longrightarrow \mathbb{N}$ of $G$. We assign to any $\ell$-tuple of paths $\pmb{\gamma} \in \Pi(G)^\ell$ in $G$ a \new{$f$-weight} defined by \[\wt_f(\pmb{\gamma}) = \sum_{v \in \Supp(\pmb{\gamma})} f(v).\]
Set $M_0^G(f) = 0$, and for all integers $\ell \geqslant 1$, $M_\ell^G(f) = \max \left(\{\wt_f(\pmb{\gamma}) \mid \pmb{\gamma} \in \Pi(G)^\ell\} \right)$. We define the \new{Greene--Kleitman invariant} of $f$ in $G$ as \[\GK_G(f) = \left(M_\ell^G(f) - M_{\ell-1}^G(f) \right)_{\ell \geqslant 1}.\]

\begin{figure}[!ht]
	\centering
	\scalebox{0.65}{\begin{tikzpicture}[line width=0.4mm,>= angle 60]
			\node(0) at (-3,-1.5){\huge $G =$};
			\node(1) at (-1,0){$\bullet$};
			\node(2) at (1.5,0.5){$\bullet$};
			\node(3) at (-2,-1){$\bullet$};
			\node(4) at (-1,-1.5){$\bullet$};
			\node(5) at (0.75,-0.5){$\bullet$};
			\node(6) at (2.5,-0.5){$\bullet$};
			\node(7) at (-2,-2){$\bullet$};
			\node(8) at (0,-1.5){$\bullet$};
			\node(9) at (1.75,-1.75){$\bullet$};
			\node(10) at (-1,-3){$\bullet$};
			\node(11) at (1,-3){$\bullet$};
			\draw[->] (1) -- (3);
			\draw[->] (1) -- (4);
			\draw[->] (1) -- (8);
			\draw[->] (2) -- (5);
			\draw[->] (2) -- (6);
			\draw[->] (3) -- (7);
			\draw[->] (4) -- (10);
			\draw[->] (5) -- (9);
			\draw[->] (6) -- (9);
			\draw[->] (7) -- (10);
			\draw[->] (8) -- (10);
			\draw[->] (8) -- (11);
			\draw[->] (9) -- (11);
			\draw[->] (5) -- (8);
			\begin{scope}[xshift=7cm]
				\node(1) at (-3,-1.5){\huge $f =$};
				\node(1) at (-1,0){$1$};
				\node(2) at (1.5,0.5){$2$};
				\node(3) at (-2,-1){$3$};
				\node(4) at (-1,-1.5){$2$};
				\node(5) at (0.75,-0.5){$2$};
				\node(6) at (2.5,-0.5){$1$};
				\node(7) at (-2,-2){$0$};
				\node(8) at (0,-1.5){$4$};
				\node(9) at (1.75,-1.75){$2$};
				\node(10) at (-1,-3){$5$};
				\node(11) at (1,-3){$1$};
				\draw[->] (1) -- (3);
				\draw[->] (1) -- (4);
				\draw[->] (1) -- (8);
				\draw[->] (2) -- (5);
				\draw[->] (2) -- (6);
				\draw[->] (3) -- (7);
				\draw[->] (4) -- (10);
				\draw[->] (5) -- (9);
				\draw[->] (6) -- (9);
				\draw[->] (7) -- (10);
				\draw[->] (8) -- (10);
				\draw[->] (8) -- (11);
				\draw[->] (9) -- (11);
				\draw[->] (5) -- (8);
			\end{scope}
			
			\begin{scope}[scale = 0.8, xshift=-4.5cm, yshift=-5.5cm]
				\node(1) at (-1,0){$1$};
				\node(2) at (1.5,0.5){$2$};
				\node(3) at (-2,-1){$3$};
				\node(4) at (-1,-1.5){$2$};
				\node(5) at (0.75,-0.5){$2$};
				\node(6) at (2.5,-0.5){$1$};
				\node(7) at (-2,-2){$0$};
				\node(8) at (0,-1.5){$4$};
				\node(9) at (1.75,-1.75){$2$};
				\node(10) at (-1,-3){$5$};
				\node(11) at (1,-3){$1$};
				\draw[->] (1) -- (3);
				\draw[->] (1) -- (4);
				\draw[->] (1) -- (8);
				\draw[->] (2) -- (5);
				\draw[->] (2) -- (6);
				\draw[->] (3) -- (7);
				\draw[->] (4) -- (10);
				\draw[->] (5) -- (9);
				\draw[->] (6) -- (9);
				\draw[->] (7) -- (10);
				\draw[->] (8) -- (10);
				\draw[->] (8) -- (11);
				\draw[->] (9) -- (11);
				\draw[->] (5) -- (8);
				\draw[-,line width=4.5mm, black!80, draw opacity = 0.3] (1.7,0.7) -- (-1.2,-3.2);
			\end{scope}
			
			\begin{scope}[scale = 0.8, xshift=1cm, yshift=-5.5cm]
				\node(1) at (-1,0){$1$};
				\node(2) at (1.5,0.5){$2$};
				\node(3) at (-2,-1){$3$};
				\node(4) at (-1,-1.5){$2$};
				\node(5) at (0.75,-0.5){$2$};
				\node(6) at (2.5,-0.5){$1$};
				\node(7) at (-2,-2){$0$};
				\node(8) at (0,-1.5){$4$};
				\node(9) at (1.75,-1.75){$2$};
				\node(10) at (-1,-3){$5$};
				\node(11) at (1,-3){$1$};
				\draw[->] (1) -- (3);
				\draw[->] (1) -- (4);
				\draw[->] (1) -- (8);
				\draw[->] (2) -- (5);
				\draw[->] (2) -- (6);
				\draw[->] (3) -- (7);
				\draw[->] (4) -- (10);
				\draw[->] (5) -- (9);
				\draw[->] (6) -- (9);
				\draw[->] (7) -- (10);
				\draw[->] (8) -- (10);
				\draw[->] (8) -- (11);
				\draw[->] (9) -- (11);
				\draw[->] (5) -- (8);
				\draw[-,line width=4.5mm, black!80, draw opacity = 0.3] (1.7,0.7) -- (0,-1.5)-- (1.2,-3.2);
				\draw[-,line width=4.5mm, black!80, draw opacity = 0.3](-0.8,0.2) -- (-2,-1) -- (-2,-2) -- (-0.8,-3.2);
			\end{scope}
			
			\begin{scope}[scale = 0.8, xshift=6.5cm, yshift=-5.5cm]
				\node(1) at (-1,0){$1$};
				\node(2) at (1.5,0.5){$2$};
				\node(3) at (-2,-1){$3$};
				\node(4) at (-1,-1.5){$2$};
				\node(5) at (0.75,-0.5){$2$};
				\node(6) at (2.5,-0.5){$1$};
				\node(7) at (-2,-2){$0$};
				\node(8) at (0,-1.5){$4$};
				\node(9) at (1.75,-1.75){$2$};
				\node(10) at (-1,-3){$5$};
				\node(11) at (1,-3){$1$};
				\draw[->] (1) -- (3);
				\draw[->] (1) -- (4);
				\draw[->] (1) -- (8);
				\draw[->] (2) -- (5);
				\draw[->] (2) -- (6);
				\draw[->] (3) -- (7);
				\draw[->] (4) -- (10);
				\draw[->] (5) -- (9);
				\draw[->] (6) -- (9);
				\draw[->] (7) -- (10);
				\draw[->] (8) -- (10);
				\draw[->] (8) -- (11);
				\draw[->] (9) -- (11);
				\draw[->] (5) -- (8);
				\draw[-,line width=4.5mm, black!80, draw opacity = 0.3] (1.4,0.7) -- (2.5,-0.5) -- (0.8,-3.2);
				\draw[-,line width=4.5mm, black!80, draw opacity = 0.3] (1.7,0.7) -- (0,-1.5)-- (1.2,-3.2);
				\draw[-,line width=4.5mm, black!80, draw opacity = 0.3](-0.8,0.2) -- (-2,-1) -- (-2,-2) -- (-0.8,-3.2); 
			\end{scope}
			
			\begin{scope}[scale = 0.8, xshift=12cm, yshift=-5.5cm]
				\node(1) at (-1,0){$1$};
				\node(2) at (1.5,0.5){$2$};
				\node(3) at (-2,-1){$3$};
				\node(4) at (-1,-1.5){$2$};
				\node(5) at (0.75,-0.5){$2$};
				\node(6) at (2.5,-0.5){$1$};
				\node(7) at (-2,-2){$0$};
				\node(8) at (0,-1.5){$4$};
				\node(9) at (1.75,-1.75){$2$};
				\node(10) at (-1,-3){$5$};
				\node(11) at (1,-3){$1$};
				\draw[->] (1) -- (3);
				\draw[->] (1) -- (4);
				\draw[->] (1) -- (8);
				\draw[->] (2) -- (5);
				\draw[->] (2) -- (6);
				\draw[->] (3) -- (7);
				\draw[->] (4) -- (10);
				\draw[->] (5) -- (9);
				\draw[->] (6) -- (9);
				\draw[->] (7) -- (10);
				\draw[->] (8) -- (10);
				\draw[->] (8) -- (11);
				\draw[->] (9) -- (11);
				\draw[->] (5) -- (8);
				\draw[-,line width=4.5mm, black!80, draw opacity = 0.3] (1.4,0.7) -- (2.5,-0.5) -- (0.8,-3.2);
				\draw[-,line width=4.5mm, black!80, draw opacity = 0.3] (1.7,0.7) -- (0,-1.5)-- (1.2,-3.2);
				\draw[-,line width=4.5mm, black!80, draw opacity = 0.3](-0.8,0.2) -- (-2,-1) -- (-2,-2) -- (-0.8,-3.2); \draw[-,line width=4.5mm, black!80, draw opacity = 0.3](-1,0.2) -- (-1,-3.2);
			\end{scope}
			
			\node(1) at (3.5,-8){\huge $\GK_G(f) = (13,5,3,2)$};
	\end{tikzpicture}}
	\caption{\label{fig:GK} An example of the computation of $\GK_G$.}
\end{figure}
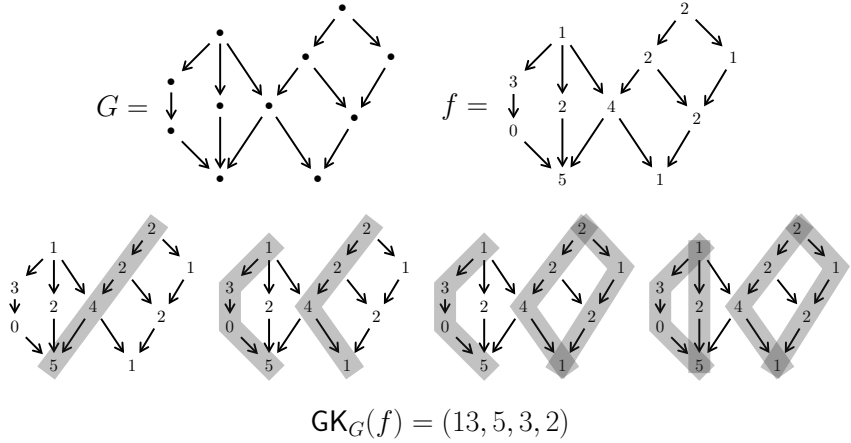

See \cref{fig:GK} for an explicit computation example.

\begin{proposition}[Greene--Kleitman \cite{GK76}] \label{GKprop} 
	Let $G$ be an acyclic directed graph and $f$ be a filling of $G$. The integer sequence $\GK_G(f)$ is an integer partition of length equal to the maximal cardinality of an antichain in $G$. Moreover, \[|\GK_G(f)| = \sum_{v \in G_0} f(v).\]
\end{proposition}

\subsection{The Gansner story}
\label{ss:Gansner} In the following, we present another way to realize the RSK correspondence. We refer the reader to \cite{Ga81Hi} for more details.

Let $A = (a_{i,j})_{1 \leqslant i,j \leqslant n}$ be a $n \times n$ integer matrix. We construct a directed graph $G_A$ where the vertices are labelled by $(i,j)$ for $1 \leqslant i,j \leqslant n$, and the arrows are given by $(i,j) \longrightarrow (i+1,j)$ and $(i,j) \longrightarrow (i,j+1)$.

For $1 \leqslant i \leqslant n$, consider $G_A^{[i,-]}$ the full subgraph of $G_A$ whose vertices are $(k,j)$ for $1 \leqslant k \leqslant i$ and $1 \leqslant j \leqslant n$. The entries of $A$ endows the graph  $G_A^{[i,-]}$ with a filling $f_A^{[i,-]}$. We define a sequence of integer partitions $(\nu^i)_{1 \leqslant i \leqslant n}$ by \[\forall i \in \{1,\ldots,n\},\ \nu^i = \GK_{G_A^{[i,-]}} \left(f_A^{[i,-]} \right).\]

Analogously, by considering, for $1 \leqslant j \leqslant n$, $G_A^{[-j]}$ the full subgraph of $G_A$ whose vertices are $(i,k)$ for $1 \leqslant i \leqslant n$ and $1 \leqslant j \leqslant k$,  we define a sequence of integer partitions $(\mu^j)_{1 \leqslant j \leqslant n}$ by
\[\forall j \in \{1,\ldots,n\},\ \mu^j = \GK_{G_A^{[-,j]}} \left(f_A^{[-,j]} \right).\]

Note $\nu^n = \mu^n$. Moreover, we can show that $\nu^{i} \subseteq \nu^{i+1}$ for all $1 \leqslant i < n$, and $\mu^j \subseteq \mu^{j+1}$ for all $1 \leqslant j < n$.

The pairs of integer partitions sequences $((\nu_i)_{1 \leqslant i \leqslant n}, (\mu_j)_{1 \leqslant j \leqslant n})$ are called the \new{Greene--Kleitman invariants of $A$}. They allow us to recover $\RSK(A)$. The semi-standard Young tableau $P(s)$ is the filling of $\mu^n$ obtained by labelling $j$ the boxes in $\mu^j \setminus \mu^{j-1}$. The other one, $Q(s)$, is constructed similarly but with the sequence $(\nu^i)$.

\begin{example} \label{ex:RSK2} In \cref{fig:RSKex2}, we give the explicit calculations of the sequences $(\mu^j)$ and $(\nu^i)$ for the matrix $A$ of \cref{ex:RSK1}.
	
	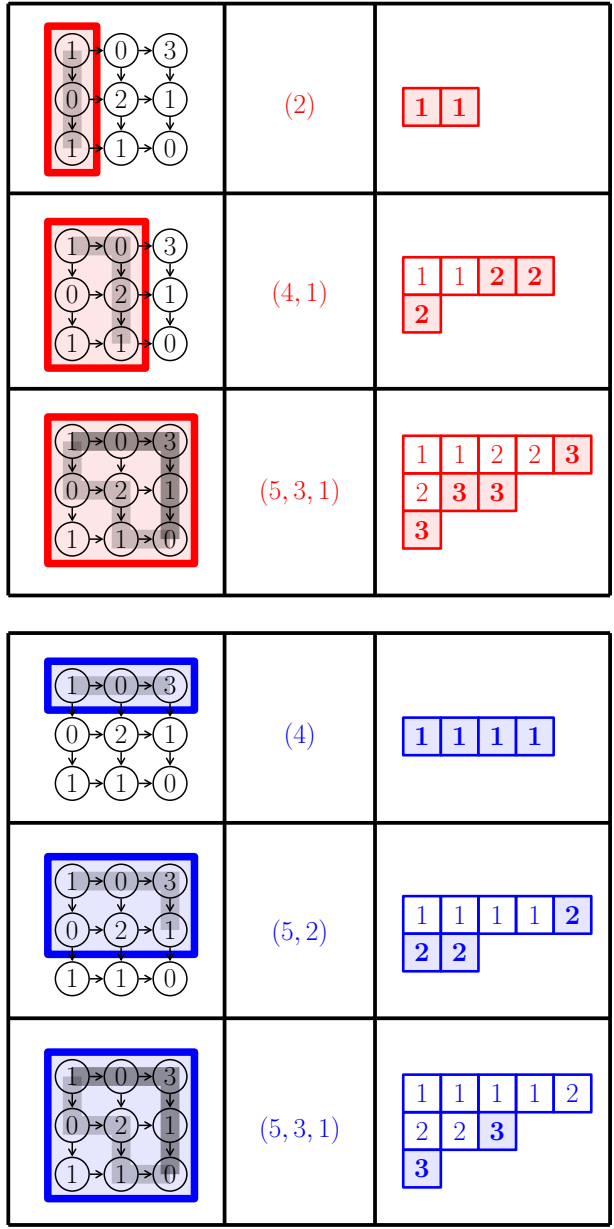
\begin{figure}[!ht]
		\centering
		
		\scalebox{0.5}{
			\begin{tikzpicture}[scale=1]
				
				\begin{scope}[xshift=1cm,yshift=4.5cm]
					\tkzDefPoint(0,0){a}
					\tkzDefPoint(0,1){b}
					\tkzDefPoint(1,1){c}
					\tkzDefPoint(1,0){d}
					\tkzDrawPolygon[line width = 0.7mm, color = red,fill=red!10](a,b,c,d);
					
					\tkzDefPoint(1,0){a}
					\tkzDefPoint(1,1){b}
					\tkzDefPoint(2,1){c}
					\tkzDefPoint(2,0){d}
					\tkzDrawPolygon[line width = 0.7mm, color = red, fill=red!10](a,b,c,d);
					
					\tkzLabelPoint[below,red](0.5,0.9){{\Huge $\mathbf{1}$}};
					\tkzLabelPoint[below,red](1.5,0.9){{\Huge $\mathbf{1}$}};
					
				\end{scope};
				
				\begin{scope}[xshift=1cm,yshift = 0cm]
					
					\tkzDefPoint(0,0){a}
					\tkzDefPoint(0,1){b}
					\tkzDefPoint(1,1){c}
					\tkzDefPoint(1,0){d}
					\tkzDrawPolygon[line width = 0.7mm, color = red](a,b,c,d);
					
					\tkzDefPoint(1,0){a}
					\tkzDefPoint(1,1){b}
					\tkzDefPoint(2,1){c}
					\tkzDefPoint(2,0){d}
					\tkzDrawPolygon[line width = 0.7mm, color = red](a,b,c,d);
					
					\tkzDefPoint(2,0){a}
					\tkzDefPoint(2,1){b}
					\tkzDefPoint(3,1){c}
					\tkzDefPoint(3,0){d}
					\tkzDrawPolygon[line width = 0.7mm, color = red,fill=red!10](a,b,c,d);
					
					\tkzDefPoint(3,0){a}
					\tkzDefPoint(3,1){b}
					\tkzDefPoint(4,1){c}
					\tkzDefPoint(4,0){d}
					\tkzDrawPolygon[line width = 0.7mm, color = red,fill=red!10](a,b,c,d);

					\tkzDefPoint(0,0){a}
					\tkzDefPoint(0,-1){b}
					\tkzDefPoint(1,-1){c}
					\tkzDefPoint(1,0){d}
					\tkzDrawPolygon[line width = 0.7mm, color = red,fill=red!10](a,b,c,d);

					\tkzLabelPoint[below,red](0.5,0.9){{\Huge $1$}};
					\tkzLabelPoint[below,red](1.5,0.9){{\Huge $1$}};
					\tkzLabelPoint[below,red](2.5,0.9){{\Huge $\mathbf{2}$}};
					\tkzLabelPoint[below,red](3.5,0.9){{\Huge $\mathbf{2}$}};
					\tkzLabelPoint[below,red](0.5,-0.1){{\Huge $\mathbf{2}$}};
					
				\end{scope};
				
				\begin{scope}[xshift=1cm,yshift = -4.75cm]
					
					\tkzDefPoint(0,0){a}
					\tkzDefPoint(0,1){b}
					\tkzDefPoint(1,1){c}
					\tkzDefPoint(1,0){d}
					\tkzDrawPolygon[line width = 0.7mm, color = red](a,b,c,d);
					
					\tkzDefPoint(1,0){a}
					\tkzDefPoint(1,1){b}
					\tkzDefPoint(2,1){c}
					\tkzDefPoint(2,0){d}
					\tkzDrawPolygon[line width = 0.7mm, color = red](a,b,c,d);
					
					\tkzDefPoint(2,0){a}
					\tkzDefPoint(2,1){b}
					\tkzDefPoint(3,1){c}
					\tkzDefPoint(3,0){d}
					\tkzDrawPolygon[line width = 0.7mm, color = red](a,b,c,d);
					
					\tkzDefPoint(3,0){a}
					\tkzDefPoint(3,1){b}
					\tkzDefPoint(4,1){c}
					\tkzDefPoint(4,0){d}
					\tkzDrawPolygon[line width = 0.7mm, color = red](a,b,c,d);
					
					\tkzDefPoint(4,1){a}
					\tkzDefPoint(4,0){b}
					\tkzDefPoint(5,0){c}
					\tkzDefPoint(5,1){d}
					\tkzDrawPolygon[line width = 0.7mm, color = red,fill=red!10](a,b,c,d);
					
					\tkzDefPoint(0,0){a}
					\tkzDefPoint(0,-1){b}
					\tkzDefPoint(1,-1){c}
					\tkzDefPoint(1,0){d}
					\tkzDrawPolygon[line width = 0.7mm, color = red](a,b,c,d);
					
					\tkzDefPoint(1,0){a}
					\tkzDefPoint(1,-1){b}
					\tkzDefPoint(2,-1){c}
					\tkzDefPoint(2,0){d}
					\tkzDrawPolygon[line width = 0.7mm, color = red,fill=red!10](a,b,c,d);
					
					\tkzDefPoint(2,0){a}
					\tkzDefPoint(2,-1){b}
					\tkzDefPoint(3,-1){c}
					\tkzDefPoint(3,0){d}
					\tkzDrawPolygon[line width = 0.7mm, color = red, fill = red!10](a,b,c,d);
					
					\tkzDefPoint(0,-2){a}
					\tkzDefPoint(0,-1){b}
					\tkzDefPoint(1,-1){c}
					\tkzDefPoint(1,-2){d}
					\tkzDrawPolygon[line width = 0.7mm, color = red,fill=red!10](a,b,c,d);
					
					\tkzLabelPoint[below,red](0.5,0.9){{\Huge $1$}};
					\tkzLabelPoint[below,red](1.5,0.9){{\Huge $1$}};
					\tkzLabelPoint[below,red](2.5,0.9){{\Huge $2$}};
					\tkzLabelPoint[below,red](3.5,0.9){{\Huge $2$}};
					\tkzLabelPoint[below,red](4.5,0.9){{\Huge $\mathbf{3}$}};
					\tkzLabelPoint[below,red](0.5,-0.1){{\Huge $2$}};
					\tkzLabelPoint[below,red](1.5,-0.1){{\Huge $\mathbf{3}$}};
					\tkzLabelPoint[below,red](2.5,-0.1){{\Huge $\mathbf{3}$}};
					\tkzLabelPoint[below,red](0.5,-1.1){{\Huge $\mathbf{3}$}};
					
				\end{scope};

				\begin{scope}[xshift=1cm,yshift=-12.25cm]
					
					\tkzDefPoint(0,0){a}
					\tkzDefPoint(0,1){b}
					\tkzDefPoint(1,1){c}
					\tkzDefPoint(1,0){d}
					\tkzDrawPolygon[line width = 0.7mm, color = blue,fill=blue!10](a,b,c,d);
					
					\tkzDefPoint(1,0){a}
					\tkzDefPoint(1,1){b}
					\tkzDefPoint(2,1){c}
					\tkzDefPoint(2,0){d}
					\tkzDrawPolygon[line width = 0.7mm, color = blue,fill=blue!10](a,b,c,d);
					
					\tkzDefPoint(2,0){a}
					\tkzDefPoint(2,1){b}
					\tkzDefPoint(3,1){c}
					\tkzDefPoint(3,0){d}
					\tkzDrawPolygon[line width = 0.7mm, color = blue,fill=blue!10](a,b,c,d);
					
					\tkzDefPoint(3,0){a}
					\tkzDefPoint(3,1){b}
					\tkzDefPoint(4,1){c}
					\tkzDefPoint(4,0){d}
					\tkzDrawPolygon[line width = 0.7mm, color = blue, fill=blue!10](a,b,c,d);

					\tkzLabelPoint[below,blue](0.5,0.9){{\Huge $\mathbf{1}$}};
					\tkzLabelPoint[below,blue](1.5,0.9){{\Huge $\mathbf{1}$}};
					\tkzLabelPoint[below,blue](2.5,0.9){{\Huge $\mathbf{1}$}};
					\tkzLabelPoint[below,blue](3.5,0.9){{\Huge $\mathbf{1}$}};

				\end{scope};

				\begin{scope}[yshift = -17cm,xshift=1cm]
					
					\tkzDefPoint(0,0){a}
					\tkzDefPoint(0,1){b}
					\tkzDefPoint(1,1){c}
					\tkzDefPoint(1,0){d}
					\tkzDrawPolygon[line width = 0.7mm, color = blue](a,b,c,d);
					
					\tkzDefPoint(1,0){a}
					\tkzDefPoint(1,1){b}
					\tkzDefPoint(2,1){c}
					\tkzDefPoint(2,0){d}
					\tkzDrawPolygon[line width = 0.7mm, color = blue](a,b,c,d);
					
					\tkzDefPoint(2,0){a}
					\tkzDefPoint(2,1){b}
					\tkzDefPoint(3,1){c}
					\tkzDefPoint(3,0){d}
					\tkzDrawPolygon[line width = 0.7mm, color = blue](a,b,c,d);
					
					\tkzDefPoint(3,0){a}
					\tkzDefPoint(3,1){b}
					\tkzDefPoint(4,1){c}
					\tkzDefPoint(4,0){d}
					\tkzDrawPolygon[line width = 0.7mm, color = blue](a,b,c,d);
					
					\tkzDefPoint(4,1){a}
					\tkzDefPoint(4,0){b}
					\tkzDefPoint(5,0){c}
					\tkzDefPoint(5,1){d}
					\tkzDrawPolygon[line width = 0.7mm, color = blue, fill = blue!10](a,b,c,d);
					
					\tkzDefPoint(0,0){a}
					\tkzDefPoint(0,-1){b}
					\tkzDefPoint(1,-1){c}
					\tkzDefPoint(1,0){d}
					\tkzDrawPolygon[line width = 0.7mm, color = blue,fill=blue!10](a,b,c,d);
					
					\tkzDefPoint(1,0){a}
					\tkzDefPoint(1,-1){b}
					\tkzDefPoint(2,-1){c}
					\tkzDefPoint(2,0){d}
					\tkzDrawPolygon[line width = 0.7mm, color = blue,fill=blue!10](a,b,c,d);
					
					\tkzLabelPoint[below,blue](0.5,0.9){{\Huge $1$}};
					\tkzLabelPoint[below,blue](1.5,0.9){{\Huge $1$}};
					\tkzLabelPoint[below,blue](2.5,0.9){{\Huge $1$}};
					\tkzLabelPoint[below,blue](3.5,0.9){{\Huge $1$}};
					\tkzLabelPoint[below,blue](4.5,0.9){{\Huge $\mathbf{2}$}};
					\tkzLabelPoint[below,blue](0.5,-0.1){{\Huge $\mathbf{2}$}};
					\tkzLabelPoint[below,blue](1.5,-0.1){{\Huge $\mathbf{2}$}};
					
				\end{scope};
				
				\begin{scope}[xshift=1cm, yshift = -21.75cm]
					
					\tkzDefPoint(0,0){a}
					\tkzDefPoint(0,1){b}
					\tkzDefPoint(1,1){c}
					\tkzDefPoint(1,0){d}
					\tkzDrawPolygon[line width = 0.7mm, color = blue](a,b,c,d);
					
					\tkzDefPoint(1,0){a}
					\tkzDefPoint(1,1){b}
					\tkzDefPoint(2,1){c}
					\tkzDefPoint(2,0){d}
					\tkzDrawPolygon[line width = 0.7mm, color = blue](a,b,c,d);
					
					\tkzDefPoint(2,0){a}
					\tkzDefPoint(2,1){b}
					\tkzDefPoint(3,1){c}
					\tkzDefPoint(3,0){d}
					\tkzDrawPolygon[line width = 0.7mm, color = blue](a,b,c,d);
					
					\tkzDefPoint(3,0){a}
					\tkzDefPoint(3,1){b}
					\tkzDefPoint(4,1){c}
					\tkzDefPoint(4,0){d}
					\tkzDrawPolygon[line width = 0.7mm, color = blue](a,b,c,d);
					
					\tkzDefPoint(4,1){a}
					\tkzDefPoint(4,0){b}
					\tkzDefPoint(5,0){c}
					\tkzDefPoint(5,1){d}
					\tkzDrawPolygon[line width = 0.7mm, color = blue](a,b,c,d);
					
					\tkzDefPoint(0,0){a}
					\tkzDefPoint(0,-1){b}
					\tkzDefPoint(1,-1){c}
					\tkzDefPoint(1,0){d}
					\tkzDrawPolygon[line width = 0.7mm, color = blue](a,b,c,d);
					
					\tkzDefPoint(1,0){a}
					\tkzDefPoint(1,-1){b}
					\tkzDefPoint(2,-1){c}
					\tkzDefPoint(2,0){d}
					\tkzDrawPolygon[line width = 0.7mm, color = blue](a,b,c,d);
					
					\tkzDefPoint(2,0){a}
					\tkzDefPoint(2,-1){b}
					\tkzDefPoint(3,-1){c}
					\tkzDefPoint(3,0){d}
					\tkzDrawPolygon[line width = 0.7mm, color = blue, fill = blue!10](a,b,c,d);
					
					\tkzDefPoint(0,-2){a}
					\tkzDefPoint(0,-1){b}
					\tkzDefPoint(1,-1){c}
					\tkzDefPoint(1,-2){d}
					\tkzDrawPolygon[line width = 0.7mm, color = blue,fill=blue!10](a,b,c,d);
					
					\tkzLabelPoint[below,blue](0.5,0.9){{\Huge $1$}};
					\tkzLabelPoint[below,blue](1.5,0.9){{\Huge $1$}};
					\tkzLabelPoint[below,blue](2.5,0.9){{\Huge $1$}};
					\tkzLabelPoint[below,blue](3.5,0.9){{\Huge $1$}};
					\tkzLabelPoint[below,blue](4.5,0.9){{\Huge $2$}};
					\tkzLabelPoint[below,blue](0.5,-0.1){{\Huge $2$}};
					\tkzLabelPoint[below,blue](1.5,-0.1){{\Huge $2$}};
					\tkzLabelPoint[below,blue](2.5,-0.1){{\Huge $\mathbf{3}$}};
					\tkzLabelPoint[below,blue](0.5,-1.1){{\Huge $\mathbf{3}$}};
					
				\end{scope};
				
				\draw[line width=1mm] (0.25,7.75) -- (0.25,-8);
				\draw[line width=1mm] (-3.75,7.75) -- (-3.75,-8);
				\draw[line width=1mm] (-9.5,7.75) -- (-9.5,-8);
				\draw[line width=1mm] (6.5,7.75) -- (6.5,-8);
				\draw[line width=1mm] (-9.5,7.75) -- (6.5,7.75);
				\draw[line width=1mm] (-9.5,2.7) -- (6.5,2.7);
				\draw[line width=1mm] (-9.5,-2.5) -- (6.5,-2.5);
				\draw[line width=1mm] (-9.5,-8) -- (6.5,-8);
				
				\begin{scope}[yshift=-16.75cm]
					\draw[line width=1mm] (0.25,7.75) -- (0.25,-8);
					\draw[line width=1mm] (-3.75,7.75) -- (-3.75,-8);
					\draw[line width=1mm] (-9.5,7.75) -- (-9.5,-8);
					\draw[line width=1mm] (6.5,7.75) -- (6.5,-8);
					\draw[line width=1mm] (-9.5,7.75) -- (6.5,7.75);
					\draw[line width=1mm] (-9.5,2.7) -- (6.5,2.7);
					\draw[line width=1mm] (-9.5,-2.5) -- (6.5,-2.5);
					\draw[line width=1mm] (-9.5,-8) -- (6.5,-8);
				\end{scope}
				
				\begin{scope}[->,line width=0.4mm,>= angle 60, scale =1.3,xshift=-6cm,yshift=5cm]
					
					\tkzDefPoint(-0.5,0.5){a}
					\tkzDefPoint(0.5,0.5){b}
					\tkzDefPoint(0.5,-2.5){c}
					\tkzDefPoint(-0.5,-2.5){d}
					\tkzDrawPolygon[line width = 2mm, color = red, fill = red!10](a,b,c,d);
					
					\draw[-,line width=5mm, black,opacity=0.2](0,0) edge (0,-2);
					
					\node[circle,draw] (a) at (0,0){{\Huge $1$}};
					\node[circle,draw] (b) at (1,0){{\Huge $0$}};
					\node[circle,draw] (c) at (2,0){{\Huge $3$}};
					
					\node[circle,draw](f) at (0,-1){{\Huge $0$}};
					\node[circle,draw] (g) at (1,-1){{\Huge $2$}};
					\node[circle,draw] (h) at (2,-1){{\Huge $1$}};
					
					\node[circle,draw] (k) at (0,-2){{\Huge $1$}};
					\node[circle,draw] (l) at (1,-2){{\Huge $1$}};
					\node[circle,draw] (m) at (2,-2){{\Huge $0$}};
					
					\draw (a) -- (b);
					\draw (a) -- (f);
					\draw (b) -- (c);
					\draw (b) -- (g);
					\draw (c) -- (h);
					\draw (f) -- (g);
					\draw (f) -- (k);
					\draw (g) -- (h);
					\draw (g) -- (l);
					\draw (h) -- (m);
					\draw (k) -- (l);
					\draw (l) -- (m);
					
				\end{scope};
				
				\begin{scope}[->,line width=0.4mm,>= angle 60, scale =1.3,xshift=-6cm,yshift=1cm]
					
					\tkzDefPoint(-0.5,0.5){a}
					\tkzDefPoint(1.5,0.5){b}
					\tkzDefPoint(1.5,-2.5){c}
					\tkzDefPoint(-0.5,-2.5){d}
					\tkzDrawPolygon[line width = 2mm, color = red, fill = red!10](a,b,c,d);
					
					\draw[-,line width=5mm, black, opacity=0.2](0,0) -- (1,0) -- (1,-2);
					
					\node[circle,draw] (a) at (0,0){{\Huge $1$}};
					\node[circle,draw] (b) at (1,0){{\Huge $0$}};
					\node[circle,draw] (c) at (2,0){{\Huge $3$}};
					
					\node[circle,draw](f) at (0,-1){{\Huge $0$}};
					\node[circle,draw] (g) at (1,-1){{\Huge $2$}};
					\node[circle,draw] (h) at (2,-1){{\Huge $1$}};
					
					\node[circle,draw] (k) at (0,-2){{\Huge $1$}};
					\node[circle,draw] (l) at (1,-2){{\Huge $1$}};
					\node[circle,draw] (m) at (2,-2){{\Huge $0$}};
					
					\draw (a) -- (b);
					\draw (a) -- (f);
					\draw (b) -- (c);
					\draw (b) -- (g);
					\draw (c) -- (h);
					\draw (f) -- (g);
					\draw (f) -- (k);
					\draw (g) -- (h);
					\draw (g) -- (l);
					\draw (h) -- (m);
					\draw (k) -- (l);
					\draw (l) -- (m);
					
				\end{scope};
				
				\begin{scope}[->,line width=0.4mm,>= angle 60, scale =1.3,xshift=-6cm,yshift=-8cm]
					
					\tkzDefPoint(-0.5,-0.5){a}
					\tkzDefPoint(-0.5,0.5){b}
					\tkzDefPoint(2.5,0.5){c}
					\tkzDefPoint(2.5,-0.5){d}
					\tkzDrawPolygon[line width = 2mm, color = blue, fill = blue!10](a,b,c,d);
					
					\draw[-,line width=5mm, black,opacity=0.2](0,0) edge (2,0);
					
					\node[circle,draw] (a) at (0,0){{\Huge $1$}};
					\node[circle,draw] (b) at (1,0){{\Huge $0$}};
					\node[circle,draw] (c) at (2,0){{\Huge $3$}};
					
					\node[circle,draw](f) at (0,-1){{\Huge $0$}};
					\node[circle,draw] (g) at (1,-1){{\Huge $2$}};
					\node[circle,draw] (h) at (2,-1){{\Huge $1$}};
					
					\node[circle,draw] (k) at (0,-2){{\Huge $1$}};
					\node[circle,draw] (l) at (1,-2){{\Huge $1$}};
					\node[circle,draw] (m) at (2,-2){{\Huge $0$}};
					
					\draw (a) -- (b);
					\draw (a) -- (f);
					\draw (b) -- (c);
					\draw (b) -- (g);
					\draw (c) -- (h);
					\draw (f) -- (g);
					\draw (f) -- (k);
					\draw (g) -- (h);
					\draw (g) -- (l);
					\draw (h) -- (m);
					\draw (k) -- (l);
					\draw (l) -- (m);
					
				\end{scope};
				
				\begin{scope}[->,line width=0.4mm,>= angle 60, scale =1.3,xshift=-6cm,yshift=-12cm]
					
					\tkzDefPoint(-0.5,-1.5){a}
					\tkzDefPoint(-0.5,0.5){b}
					\tkzDefPoint(2.5,0.5){c}
					\tkzDefPoint(2.5,-1.5){d}
					\tkzDrawPolygon[line width = 2mm, color = blue, fill = blue!10](a,b,c,d);
					
					\draw[-,line width=5mm, black,opacity=0.2](0,0) -- (2,0) -- (2,-1);
					
					\node[circle,draw] (a) at (0,0){{\Huge $1$}};
					\node[circle,draw] (b) at (1,0){{\Huge $0$}};
					\node[circle,draw] (c) at (2,0){{\Huge $3$}};
					
					\node[circle,draw](f) at (0,-1){{\Huge $0$}};
					\node[circle,draw] (g) at (1,-1){{\Huge $2$}};
					\node[circle,draw] (h) at (2,-1){{\Huge $1$}};
					
					\node[circle,draw] (k) at (0,-2){{\Huge $1$}};
					\node[circle,draw] (l) at (1,-2){{\Huge $1$}};
					\node[circle,draw] (m) at (2,-2){{\Huge $0$}};
					
					\draw (a) -- (b);
					\draw (a) -- (f);
					\draw (b) -- (c);
					\draw (b) -- (g);
					\draw (c) -- (h);
					\draw (f) -- (g);
					\draw (f) -- (k);
					\draw (g) -- (h);
					\draw (g) -- (l);
					\draw (h) -- (m);
					\draw (k) -- (l);
					\draw (l) -- (m);
					
				\end{scope};
				
				\begin{scope}[->,line width=0.4mm,>= angle 60, scale =1.3,xshift=-6cm,yshift=-3cm]
					
					\tkzDefPoint(-0.5,-2.5){a}
					\tkzDefPoint(-0.5,0.5){b}
					\tkzDefPoint(2.5,0.5){c}
					\tkzDefPoint(2.5,-2.5){d}
					\tkzDrawPolygon[line width = 2mm, color = red, fill = red!10](a,b,c,d);
					
					\draw[-,line width=5mm, black,opacity=0.4](0,0) -- (2,0) -- (2,-2);
					
					\draw[-,line width=5mm, black,opacity=0.2](0,0) -- (0,-1) -- (1,-1) -- (1,-2) -- (2,-2);
					
					\node[circle,draw] (a) at (0,0){{\Huge $1$}};
					\node[circle,draw] (b) at (1,0){{\Huge $0$}};
					\node[circle,draw] (c) at (2,0){{\Huge $3$}};
					
					\node[circle,draw](f) at (0,-1){{\Huge $0$}};
					\node[circle,draw] (g) at (1,-1){{\Huge $2$}};
					\node[circle,draw] (h) at (2,-1){{\Huge $1$}};
					
					\node[circle,draw] (k) at (0,-2){{\Huge $1$}};
					\node[circle,draw] (l) at (1,-2){{\Huge $1$}};
					\node[circle,draw] (m) at (2,-2){{\Huge $0$}};
					
					\draw (a) -- (b);
					\draw (a) -- (f);
					\draw (b) -- (c);
					\draw (b) -- (g);
					\draw (c) -- (h);
					\draw (f) -- (g);
					\draw (f) -- (k);
					\draw (g) -- (h);
					\draw (g) -- (l);
					\draw (h) -- (m);
					\draw (k) -- (l);
					\draw (l) -- (m);
					
				\end{scope};
				
				\begin{scope}[->,line width=0.4mm,>= angle 60, scale =1.3,xshift=-6cm,yshift=-16cm]
					
					\tkzDefPoint(-0.5,-2.5){a}
					\tkzDefPoint(-0.5,0.5){b}
					\tkzDefPoint(2.5,0.5){c}
					\tkzDefPoint(2.5,-2.5){d}
					\tkzDrawPolygon[line width = 2mm, color = blue, fill = blue!10](a,b,c,d);
					
					\draw[-,line width=5mm, black, opacity=0.4](0,0) -- (2,0) -- (2,-2);
					
					\draw[-,line width=5mm, black,opacity=0.2](0,0) -- (0,-1) -- (1,-1) -- (1,-2) -- (2,-2);
					
					\node[circle,draw] (a) at (0,0){{\Huge $1$}};
					\node[circle,draw] (b) at (1,0){{\Huge $0$}};
					\node[circle,draw] (c) at (2,0){{\Huge $3$}};
					
					\node[circle,draw](f) at (0,-1){{\Huge $0$}};
					\node[circle,draw] (g) at (1,-1){{\Huge $2$}};
					\node[circle,draw] (h) at (2,-1){{\Huge $1$}};
					
					\node[circle,draw] (k) at (0,-2){{\Huge $1$}};
					\node[circle,draw] (l) at (1,-2){{\Huge $1$}};
					\node[circle,draw] (m) at (2,-2){{\Huge $0$}};
					
					\draw (a) -- (b);
					\draw (a) -- (f);
					\draw (b) -- (c);
					\draw (b) -- (g);
					\draw (c) -- (h);
					\draw (f) -- (g);
					\draw (f) -- (k);
					\draw (g) -- (h);
					\draw (g) -- (l);
					\draw (h) -- (m);
					\draw (k) -- (l);
					\draw (l) -- (m);
					
				\end{scope};

				\tkzLabelPoint[below,red](-1.75,5.5){{\Huge $(2)$}};
				\tkzLabelPoint[below,red](-1.75,0.5){{\Huge $(4,1)$}};
				\tkzLabelPoint[below,red](-1.75,-4.75){{\Huge $(5,3,1)$}};

				\tkzLabelPoint[below,blue](-1.75,-11.25){{\Huge $(4)$}};
				\tkzLabelPoint[below,blue](-1.75,-16.5){{\Huge $(5,2)$}};
				\tkzLabelPoint[below,blue](-1.75,-21.75){{\Huge $(5,3,1)$}};
		\end{tikzpicture}}
		
		\caption{\label{fig:RSKex2} Illustration of the calculations for $\RSK(A)$ in the ``Greene-Kleitman'' way. (Top): The sequence  {\color{red}{$(\mu^j)$}}; (Bottom): The sequence {\color{blue}{$(\nu^i)$}}.}
	\end{figure}
\end{example}

\subsection{The generalized Gansner story}
\label{ss:GeneralizedRSK} 
\begin{figure}[!ht]
	\centering
	\scalebox{0.5}{
		\begin{tikzpicture}[scale=1]
			
			\tkzDefPoint(0,0){a}
			\tkzDefPoint(0,1){b}
			\tkzDefPoint(1,1){c}
			\tkzDefPoint(1,0){d}
			\tkzDrawPolygon[line width = 0.7mm, color = black](a,b,c,d);
			
			\tkzDefPoint(0,-1){a}
			\tkzDefPoint(0,0){b}
			\tkzDefPoint(1,0){c}
			\tkzDefPoint(1,-1){d}
			\tkzDrawPolygon[line width = 0.7mm, color = black](a,b,c,d);
			
			\tkzDefPoint(0,-2){a}
			\tkzDefPoint(0,-1){b}
			\tkzDefPoint(1,-1){c}
			\tkzDefPoint(1,-2){d}
			\tkzDrawPolygon[line width = 0.7mm, color = black](a,b,c,d);
			
			\tkzDefPoint(1,0){a}
			\tkzDefPoint(1,1){b}
			\tkzDefPoint(2,1){c}
			\tkzDefPoint(2,0){d}
			\tkzDrawPolygon[line width = 0.7mm, color = black](a,b,c,d);
			
			\tkzDefPoint(1,-1){a}
			\tkzDefPoint(1,0){b}
			\tkzDefPoint(2,0){c}
			\tkzDefPoint(2,-1){d}
			\tkzDrawPolygon[line width = 0.7mm, color = black](a,b,c,d);
			
			\tkzDefPoint(1,-2){a}
			\tkzDefPoint(1,-1){b}
			\tkzDefPoint(2,-1){c}
			\tkzDefPoint(2,-2){d}
			\tkzDrawPolygon[line width = 0.7mm, color = black](a,b,c,d);
			
			\tkzDefPoint(2,0){a}
			\tkzDefPoint(2,1){b}
			\tkzDefPoint(3,1){c}
			\tkzDefPoint(3,0){d}
			\tkzDrawPolygon[line width = 0.7mm, color = black](a,b,c,d);
			
			\tkzDefPoint(2,-1){a}
			\tkzDefPoint(2,0){b}
			\tkzDefPoint(3,0){c}
			\tkzDefPoint(3,-1){d}
			\tkzDrawPolygon[line width = 0.7mm, color = black](a,b,c,d);
			
			\tkzDefPoint(2,-2){a}
			\tkzDefPoint(2,-1){b}
			\tkzDefPoint(3,-1){c}
			\tkzDefPoint(3,-2){d}
			\tkzDrawPolygon[line width = 0.7mm, color = black](a,b,c,d);
			
			\tkzLabelPoint[below](0.5,0.9){{\Huge $1$}};
			\tkzLabelPoint[below](1.5,0.9){{\Huge $0$}};
			\tkzLabelPoint[below](2.5,0.9){{\Huge $3$}};
			\tkzLabelPoint[below](0.5,-0.1){{\Huge $0$}};
			\tkzLabelPoint[below](1.5,-0.1){{\Huge $2$}};
			\tkzLabelPoint[below](2.5,-0.1){{\Huge $1$}};
			\tkzLabelPoint[below](0.5,-1.1){{\Huge $1$}};
			\tkzLabelPoint[below](1.5,-1.1){{\Huge $1$}};
			\tkzLabelPoint[below](2.5,-1.1){{\Huge $0$}};
			
			\begin{scope}[xshift=8cm]
				
				\tkzDefPoint(0,0){a}
				\tkzDefPoint(0,1){b}
				\tkzDefPoint(1,1){c}
				\tkzDefPoint(1,0){d}
				\tkzDrawPolygon[line width = 0.7mm, color = black](a,b,c,d);
				
				\tkzDefPoint(0,-1){a}
				\tkzDefPoint(0,0){b}
				\tkzDefPoint(1,0){c}
				\tkzDefPoint(1,-1){d}
				\tkzDrawPolygon[line width = 0.7mm, color = black](a,b,c,d);
				
				\tkzDefPoint(0,-2){a}
				\tkzDefPoint(0,-1){b}
				\tkzDefPoint(1,-1){c}
				\tkzDefPoint(1,-2){d}
				\tkzDrawPolygon[line width = 0.7mm, color = black](a,b,c,d);
				
				\tkzDefPoint(1,0){a}
				\tkzDefPoint(1,1){b}
				\tkzDefPoint(2,1){c}
				\tkzDefPoint(2,0){d}
				\tkzDrawPolygon[line width = 0.7mm, color = black](a,b,c,d);
				
				\tkzDefPoint(1,-1){a}
				\tkzDefPoint(1,0){b}
				\tkzDefPoint(2,0){c}
				\tkzDefPoint(2,-1){d}
				\tkzDrawPolygon[line width = 0.7mm, color = black](a,b,c,d);
				
				\tkzDefPoint(1,-2){a}
				\tkzDefPoint(1,-1){b}
				\tkzDefPoint(2,-1){c}
				\tkzDefPoint(2,-2){d}
				\tkzDrawPolygon[line width = 0.7mm, color = black](a,b,c,d);
				
				\tkzDefPoint(2,0){a}
				\tkzDefPoint(2,1){b}
				\tkzDefPoint(3,1){c}
				\tkzDefPoint(3,0){d}
				\tkzDrawPolygon[line width = 0.7mm, color = black](a,b,c,d);
				
				\tkzDefPoint(2,-1){a}
				\tkzDefPoint(2,0){b}
				\tkzDefPoint(3,0){c}
				\tkzDefPoint(3,-1){d}
				\tkzDrawPolygon[line width = 0.7mm, color = black](a,b,c,d);
				
				\tkzDefPoint(2,-2){a}
				\tkzDefPoint(2,-1){b}
				\tkzDefPoint(3,-1){c}
				\tkzDefPoint(3,-2){d}
				\tkzDrawPolygon[line width = 0.7mm, color = black](a,b,c,d);
				
				\tkzLabelPoint[below,Plum](0.5,0.9){{\Huge $1$}};
				\tkzLabelPoint[below,blue](1.5,0.9){{\Huge $2$}};
				\tkzLabelPoint[below,blue](2.5,0.9){{\Huge $4$}};
				\tkzLabelPoint[below,red](0.5,-0.1){{\Huge $1$}};
				\tkzLabelPoint[below,Plum](1.5,-0.1){{\Huge $3$}};
				\tkzLabelPoint[below,blue](2.5,-0.1){{\Huge $5$}};
				\tkzLabelPoint[below,red](0.5,-1.1){{\Huge $2$}};
				\tkzLabelPoint[below,red](1.5,-1.1){{\Huge $4$}};
				\tkzLabelPoint[below,Plum](2.5,-1.1){{\Huge $5$}};
				
			\end{scope}
			
			\begin{scope}[xshift=2cm]
				\draw [line width=0.7mm, blue!30, dashed] (7.5,1.5) --  (9.5,-0.5);
				\node[blue] at (10,-0.5){{\Huge $\nu^1$}};
				\draw [line width=0.7mm, blue!30, dashed] (6.5,1.5) --  (9.5,-1.5);
				\node[blue] at (10,-1.5){{\Huge $\nu^2$}};
				\draw [line width=0.7mm, Plum!30, dashed] (5.5,1.5) --  (9.5,-2.5);
				\node[Plum] at (10.5,-3){{\Huge $\nu^3 = \mu^3$}};
				\draw [line width=0.7mm, red!30, dashed] (5.5,0.5) --  (8.5,-2.5);
				\node[red] at (8.5,-3){{\Huge $\mu^2$}};
				\draw [line width=0.7mm, red!30, dashed] (5.5,-0.5) --  (7.5,-2.5);
				\node[red] at (7.5,-3){{\Huge $\mu^1$}};
			\end{scope}
			
			\draw [|->,line width=1.5mm] (4,-0.5) -- node[above]{{\Huge $\RSK $}} (7,-0.5);
			
	\end{tikzpicture}}
	\caption{\label{fig:GenRSKex} Displaying the results from the application of $\RSK $ seen in \cref{fig:RSKex2} in the ``Greene-Kleitman'' way.}
\end{figure}
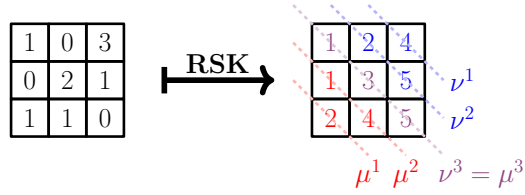

We can display these sequences  $(\mu^j)_{1 \leqslant j \leqslant n}$ and $(\nu^i)_{1 \leqslant j \leqslant n}$ as a reverse plane partition of shape the $n\times n$ box partition. We first locate the sink of the subgraph used to calculate $\nu^j$ or $\mu^i$. This sink corresponds to a box in the Young tableau we want. We fill the boxes in the same diagonal with the parts of $\nu^j$ or $\mu^i$ from the bottom left to the top right. We can look at \cref{fig:GenRSKex} to see how it goes with the results of \cref{ex:RSK2}.

We can generalize this way of applying $\RSK$ to a correspondence from fillings of a Ferrers diagram of any fixed integer partition $\lambda$ to reverse plane partitions of shape $\lambda$.

Let $f : \Fer(\lambda) \longrightarrow \mathbb{N}$ be a filling of $\Fer(\lambda)$
We consider the directed graph $G_\lambda$ associated to $\lambda$ whose vertex set is $\Fer(\lambda)$ and the arrows are given by $(i,j) \longrightarrow (i,j+1)$ and $(i,j) \longrightarrow (i+1,j)$. Indeed $G_\lambda$ corresponds to the Hasse diagram of $(\Fer(\lambda), \unlhd)$.  

\begin{figure}[!ht]
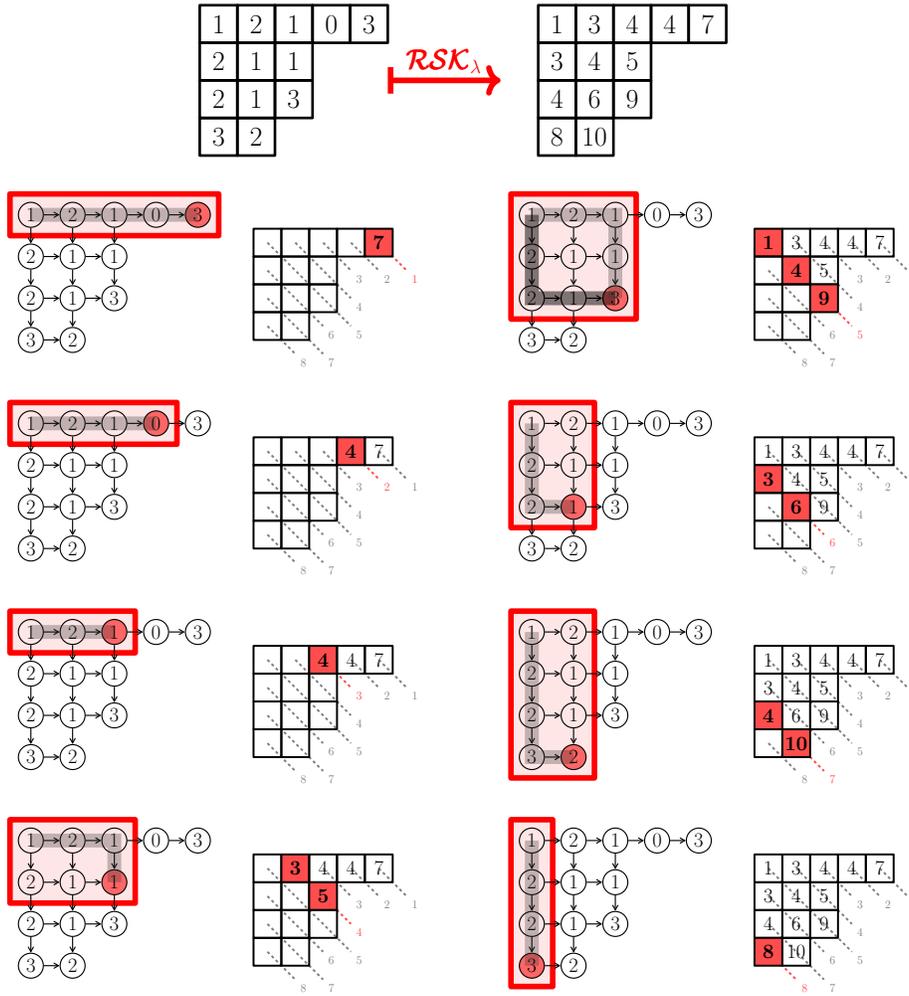

	\centering
	\scalebox{0.5}{
}
	\caption{\label{fig:genRSK} Explicit calculations of $\GRSK_\lambda(f)$ for a given filling $f$ of shape $\lambda = (5,3,3,2)$. For $1 \leqslant m \leqslant 8$, each framed subgraph corresponds to the subgraph $G_\lambda(m)$, and each filled diagonal colored in red correspond to $\GK_{G_\lambda(m)}(f)$.}
\end{figure}

Assume that $\lambda \in \Hk_n$. For each $k \in \{1,\ldots,n\}$, we consider $G_\lambda^{[k]}$ the full subgraph of $G_\lambda$ whose vertices are boxes in $\square_k(\lambda)$. Write $f^{[k]}$ for the induced filling of $G_\lambda^{[k]}$ from $f$ of $G_\lambda$. We define an integer partition $\pi^k$ by:
\[\pi^k = \GK_{G_\lambda^{[k]}}\left(f^{[k]}\right).\]
We can show that $\ell(\pi^k) \leqslant \# D_k(\lambda)$, as there exists $\pmb{\gamma} \in \Pi(G_\lambda^{[k]})^{\#D_k(\lambda)}$ such that $\Supp(\pmb{\gamma}) = (G_\lambda^{[k]})_0$. So we can place the values of $\pi^k$ on the diagonal $D_k(\lambda)$. We define a filling $\GRSK_\lambda(f): \Fer(\lambda) \longrightarrow \mathbb{N}$ as follows. Using the $\lambda$-diagonal coordinates for the boxes of $\lambda$, we compute $\GRSK_\lambda(f)$ by: \[\ \forall \ldiag{k,\delta}_\lambda \in \Fer(\lambda),\ \GRSK_\lambda(f)(\ldiag{k,\delta}_\lambda) = \pi^{k}_{\delta}.\]
See \cref{fig:genRSK} for a detailed example.

\begin{theorem}[{\cite{Ga81Hi,Ga81Ma}}] \label{thm:fillingsRPP}
	For any nonempty integer partition $\lambda$, The map $\GRSK_\lambda$ is a one-to-one correspondence from fillings of $\Fer(\lambda)$ to $\RPP(\lambda)$.
\end{theorem}

\begin{remark}
	By reversing the arrows $(i,j) \longrightarrow (i+1,j)$ in $G_\lambda$, and by proceeding to the same calculations that defined $\GRSK_\lambda$, it realizes the Hillman--Grassl correspondence. See \cite{Ga81Hi} for more details.
\end{remark}

We introduce some notations before stating a combinatorial identity, consequence of  \cref{thm:fillingsRPP}. Fix $n \in \mathbb{N}^*$. Consider $\lambda \in \Hk_n$. We assign a weight to the boxes of $\lambda$, using the $\lambda$-diagonal coordinates, as follows: \[\forall b = \ldiag{k,\delta}_\lambda \in \Fer(\lambda),\ w_{\lambda,b}(x_1,\ldots,x_n) = \prod_{1 \leqslant \ell \leqslant n,\ b \in \square_\ell(\lambda)} x_{\ell}.\]
We define the \new{trace generating function for $\RPP(\lambda)$}:
\[\rho_\lambda(x_1,\ldots,x_n) = \sum_{f \in \RPP(\lambda)} \prod_{\ldiag{\ell,\varepsilon}_\lambda \in \Fer(\lambda)} x_\ell^{f(\ldiag{\ell,\varepsilon}_\lambda)}.\] 

\begin{corollary}[\cite{Ga81Hi}] \label{cor:tracegenformula}
	Let $n \geqslant 1$ and $\lambda \in \Hk_n$. Let $x_1,\ldots,x_n$ be $n$ formal variables. We have \[\rho_\lambda(x_1,\ldots,x_n) = \prod_{b \in \Fer(\lambda)} \dfrac{1}{1 - w_{\lambda,b}(x_1,\ldots,x_n)}.\]
\end{corollary}

This result induces a well-known equality involving the norm-generating function of $\RPP(\lambda)$, previously proved by R. Stanley \cite{St72} using $\mathcal{P}$-partition theory. Here, by setting $\sigma(f) = \sum_{b \in \Fer(\lambda)} f(b)$ and mapping all the $x_i$ to $x$, we have \[\sum_{f \in \RPP(\lambda)} x^{\sigma(f)} = \prod_{b \in \Fer(\lambda)} \dfrac{1}{1 - x^{h_\lambda(b)}}.\]

	\section{Tools from Coxeter elements}
	\label{sec:Coxeter}
	In this section, we define combinatorial objects related to Coxeter elements that help present and study our extended version of Gansner's RSK correspondence.

\subsection{(Type \texorpdfstring{$A$}{A}) Coxeter elements}
\label{ss:typeACox}

For any $n \geqslant 1$, let $\mathfrak{S}_{n+1}$ be the symmetric group on $n+1$ letters. For $1 \leqslant i < j \leqslant n+1$, write $(i,j)$ for the transposition exchanging $i$ and $j$. For $1 \leqslant i \leqslant n$, let $s_i$ be the adjacent transposition $(i,i+1)$. Let $\Sigma_n$ be the set of adjacent transpositions of $\mathfrak{S}_{n+1}$. Recall that $\mathfrak{S}_{n+1}$ admits a presentation in terms of generators and relations using $\Sigma_n$ as follows:
\[\mathfrak{S}_{n+1} = \left\langle \Sigma_n \left|  \begin{matrix*}[l]
	s_i^2 = 1 & \text{for } i \in \{1,\ldots,n\}  \\
	s_i s_{i+1} s_i = s_{i+1} s_i s_{i+1}  &\text{for } i \in \{1,\ldots,n-1\} \\
	s_i s_j = s_j s_i & \text{for } i,j \in \{1,\ldots,n\} \text{ such that } |i-j| > 1 \}
\end{matrix*}\right.\right\rangle\]

For any $w \in \mathfrak{S}_{n+1}$, call an \textit{expression of} $w$ a way to write $w$ as a product of transpositions in $\Sigma_n$. The \emph{length of $w$}, denoted by $\ell(w)$, is the minimal number of transpositions in $\Sigma_n$ needed to express $w$. Whenever  $\ell(s w) < \ell(w)$ for some $s \in \Sigma_n$, we say that $s$ is \emph{initial in} $w$. Similarly, we call $s \in \Sigma_n$ \emph{final in} $w$ whenever $\ell(w s) < \ell (w)$.

A \new{Coxeter element (of $\mathfrak{S}_{n+1}$)} is an element $c \in \mathfrak{S}_{n+1}$ which can be written as a product of all the transpositions of $\Sigma_n$, in some order, where each of them appears precisely once.

\begin{example} \label{ex:CoxA} The permutation $c = s_2 s_1 s_3 s_6 s_5 s_4 s_8 s_7 = (1,3,4,7,9,8,6,5,2)$ is a Coxeter element of $\mathfrak{S}_9$. Note that $s_1,s_3,s_6$ and $s_8$ are initial in $c$, and $s_7$ and $s_4$ are final in $c$.
\end{example}

First, we state the result on conjugating Coxeter elements with their initial or final adjacent transpositions. 

\begin{lemma} \label{lem:Coxinitfin} $ $
	\begin{enumerate}[label=$(\roman*)$,itemsep=1mm]
		\item Let $c \in \mathfrak{S}_{n+1}$ be a Coxeter element. For any $s \in \Sigma_n$, either initial or final in $c$, the permutation $scs$ is a Coxeter element in $\mathfrak{S}_{n+1}$.
		\item {\cite[Lemma 1.7]{R07}} Given two Coxeter elements $c,c' \in \mathfrak{S}_{n+1}$, there exist $k \in \mathbb{N}$, and a sequence $c=c_0,c_1,\ldots,c_k = c'$ of Coxeter elements such that for any $i \in \{0,\ldots,k-1\}$ there exists an initial adjacent transposition $t_i$ such that $c_{i+1} = t_i c_i t_i$.
	\end{enumerate}
\end{lemma}

Then, as observed in \cref{ex:CoxA}, the following lemma allows us to write any Coxeter element of $\mathfrak{S}_{n+1}$ as a long cycle of a precise form. It is a direct consequence of \cref{lem:Coxinitfin}.

\begin{lemma} \label{lem:Coxcyc}
	An element $c \in \mathfrak{S}_{n+1}$ is a Coxeter element if and only if $c$ is a long cycle which can be written as follows \[c = (c_1,c_2, \ldots, c_m, c_{m+1}, \ldots, c_{n+1})\]
	where $c_1 =1 < c_2 <\ldots <c_m = n+1 > c_{m+1} > \ldots > c_{n+1} > c_1 = 1$.
\end{lemma}

Consider a Coxeter element $c \in \mathfrak{S}_{n+1}$. Write it $c = (c_1, \ldots, c_{n+1})$ as said in the previous lemma. We define the \new{left part of $c$} as $\L_c= \{c_2,\ldots,c_{m-1}\}$ and the \new{right part of $c$} as $\R_c = \{c_{m+1},\ldots,c_{n+1}\}$.

The following lemma characterizes initial and final adjacent transpositions thanks to $\L_c$ and $\R_c$.

\begin{lemma} \label{lem:initialfinalandleftright}
	Let $c \in \mathfrak{S}_{n+1}$ be a Coxeter element. For any $k \in \{2,\ldots,n-1\}$,
	\begin{enumerate}[label=$\bullet$,itemsep=1mm] 
		\item $s_k$ is final in $c$ if and only if $k \in \L_c$ and $k+1 \in \R_c$; and,
		\item $s_k$ is initial in $c$ if and only if $k \in \R_c$ and $k+1 \in \L_c$.
	\end{enumerate} 
	In special cases,
	\begin{enumerate}[label=$\bullet$,itemsep=1mm] 
		\item $s_1$ is initial in $c$ if and only if $2 \in \L_c$; and,
		\item $s_n$ is final in $c$ if and only if then $n \in \L_c$.
	\end{enumerate} 
\end{lemma}

We recall, via the following definition, that we can associate to any Coxeter element $c \in \mathfrak{S}_{n+1}$ a unique $A_n$ type quiver.

\begin{definition} \label{def:CoxQuiver}
	Let $c \in \mathfrak{S}_{n+1}$ be a Coxeter element.  We define the quiver \new{$Q(c)$} as follows:
	\begin{enumerate}[label = $\bullet$, itemsep=1mm]
		\item its set of vertices is $Q(c)_0 = \{1,\ldots,n\}$;
		
		\item its set of arrows is given by arrows between $i$ and $i+1$, for all $i \in \{1,\ldots,n\}$:
		\begin{enumerate}[label=$\bullet$,itemsep=0.5mm]
			\item we have $i \longrightarrow i+1$ if $s_i$ precedes $s_{i+1}$ in a reduced expression of $c$;
			\item we have $i \longleftarrow i+1$ otherwise.
		\end{enumerate}
	\end{enumerate} 
\end{definition}

\begin{example} For $c = (1,3,4,7,9,8,6,5,2)$, we obtain
	\begin{center}
		\begin{tikzpicture}[->,line width=0.6mm,>= angle 60,color=black,scale=0.5]
			\node (Q) at (-2,0){$Q(c)=$};
			\node (1) at (0,0){$1$};
			\node (2) at (2,0){$2$};
			\node (3) at (4,0){$3$};
			\node (4) at (6,0){$4$};
			\node (5) at (8,0){$5$};
			\node (6) at (10,0){$6$};
			\node (7) at (12,0){$7$};
			\node (8) at (14,0){$8$};
			\draw (2) -- (1);
			\draw (2) -- (3);
			\draw (3) -- (4);
			\draw (5) -- (4);
			\draw (6) -- (5);
			\draw (6) -- (7);
			\draw (8) -- (7);
		\end{tikzpicture}
	\end{center}
\end{example}

\begin{proposition} \label{prop:CoxQuiv}
	Let $n \in \mathbb{N}^*$. The map $c \longmapsto Q(c)$ realizes a one-to-one correspondence from Coxeter elements of $\mathfrak{S}_{n+1}$ to $A_n$ type quivers. Moreover:
	\begin{enumerate}[label=$\bullet$,itemsep=1mm]
		\item $v$ is a source of $Q(c)$ if and only if $s_v$ is initial in $c$;
		\item $v$ is a sink of $Q(c)$ if and only if $s_v$ is final in $c$; 
	\end{enumerate}
\end{proposition}

This map is crucial for connecting to representation-theoretic results.

Finally, we give a tiny result that links the inverse operation on Coxeter elements in $\mathfrak{S}_{n+1}$ and the opposite action on $A_n$ type quivers.

\begin{lemma} \label{lem:inverseandop}
	Let $c \in \mathfrak{S}_{n+1}$ be a Coxeter element. Then $Q(c^{-1}) = Q(c)^{\op}$.
\end{lemma}

\subsection{Interval bipartitions}
\label{ss:intbipartitions}

Let $\A \subset \mathbb{N}^*$. A \new{bipartition} of $\A$ is a pair $(\L,\R)$ sucht that $\L \cup \R = \A$ and $\L \cap \R = \varnothing$. We do not identify the pair $(\L,\R)$ with the pair $(\R,\L)$. The following result is a direct consequence of \cref{lem:Coxcyc}.

\begin{corollary} \label{cor:Coxbipart} Let $n \geqslant 1$. The map
	\[\Psi_{n+1} : \left\{\begin{matrix}
		\{\text{Coxeter elements of } \mathfrak{S}_{n+1}\} & \longrightarrow & \{\text{Bipartitions of } \{2, \ldots, n\}\} \\
		c & \longmapsto & (\L_c, \R_c)
	\end{matrix}\right.\]
	is bijective.
\end{corollary}
In the following, we focus on bipartitions of intervals in $\mathbb{N}^*$.

An \emph{interval (in $\mathbb{N}^*$)} is a set $ \llrr{i,j} = \{i,i+1,\ldots,j\}$ for some $i,j \in \mathbb{N}^*$ with $i \leqslant j$. For all $i \in \mathbb{N}^*$, we set $\llrr{i,i} = \llrr{i}$. We denote by $\mathcal{I}$ the set of all the intervals in $\mathbb{N}^*$, and by $\mathcal{I}_n$ the subset of those included in $\{1,\ldots,n+1\}$.  An \new{interval bipartition} is a bipartition $(\L,\R)$ of an interval in $\mathbb{N}^*$. Call it \new{elementary} whenever either $\L = \R = \varnothing$, or both $1 \in \L$ and $\max(\L \cup \R) \in \R$.

Fix $(\L,\R)$ as a bipartition of some finite set $\A \subset \mathbb{N}^*$. If $\L$ is nonempty, write $\L = \{\ell_1 < \ell_2 < \ldots < \ell_p\}$.  We define the integer partition $\oplamb(\L,\R)$, for all $i \in \{1,\ldots,p\}$, by $\oplamb(\L,\R)_i = \#\{r \in \R\ \mid \ell_i < r\}$ if $\L$ is not empty, and $\oplamb(\L,\R) = (0)$ otherwise. 

\begin{lemma}\label{lem:elemintervbipartition}
	For any bipartition $(\L,\R)$ of some finite set $\A \subset \mathbb{N}^*$, there exists an elementary interval bipartition $(\L',\R')$ such that $\oplamb(\L,\R) = \oplamb(\L',\R')$ 
\end{lemma}

\begin{proof}
	If $\oplamb(\L,\R) = (0)$, then we set $\L' = \R' = \varnothing$ and we are done.
	
	Otherwise consider $\mathbf{M} = \{\ell \in \L \mid \ell < \max(\R)\}$ and $\mathbf{S} = \{r \in \R \mid r > \min(\L)\}$. By construction, we easily check that $\oplamb(\L,\R) = \oplamb(\mathbf{M},\mathbf{S})$. Let $p = \#(\A \cup \B)$, and consider the stricly increasing map from $\varphi: \{1,\ldots, p\} \longrightarrow \mathbf{M} \cup \mathbf{S}$. By setting $\L' = \varphi^{-1}(\mathbf{M})$ and $\R' = \varphi^{-1}(\mathbf{S})$, we can check that $(\L',\R')$ is an elementary interval bipartition of $\llrr{1,p}$, and $\oplamb(\L',\R') = \oplamb(\mathbf{M},\mathbf{S}) = \oplamb(\L,\R)$.
\end{proof}

From now on, we assume that $(\L,\R)$ is an elementary interval bipartition. By also writing $\R = \{r_1 < \ldots < r_q\}$, we can picture $\oplamb(\L,\R)$ by its Ferrers diagram: we have $(i,j) \in \Fer(\oplamb(\L,\R))$ whenever $\ell_i < r_{q-j+1}$. It allows us to label the $i$th row of $\Fer(\oplamb(\L,\R))$ by $\ell_i$ and the $j$th column by $r_{q-j+1}$. 

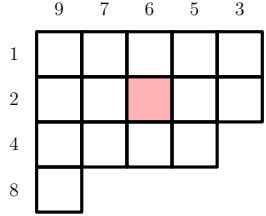
\begin{figure}[!ht]
	\centering
	\scalebox{0.6}{\begin{tikzpicture}
			\tkzDefPoint(0,0){a}
			\tkzDefPoint(0,1){b}
			\tkzDefPoint(1,1){c}
			\tkzDefPoint(1,0){d}
			\tkzDrawPolygon[line width = 0.7mm, color = black](a,b,c,d);
			
			\tkzDefPoint(1,0){a}
			\tkzDefPoint(1,1){b}
			\tkzDefPoint(2,1){c}
			\tkzDefPoint(2,0){d}
			\tkzDrawPolygon[line width = 0.7mm, color = black](a,b,c,d);
			
			\tkzDefPoint(2,0){a}
			\tkzDefPoint(2,1){b}
			\tkzDefPoint(3,1){c}
			\tkzDefPoint(3,0){d}
			\tkzDrawPolygon[line width = 0.7mm, color = black](a,b,c,d);
			
			\tkzDefPoint(3,0){a}
			\tkzDefPoint(3,1){b}
			\tkzDefPoint(4,1){c}
			\tkzDefPoint(4,0){d}
			\tkzDrawPolygon[line width = 0.7mm, color = black](a,b,c,d);
			
			\tkzDefPoint(4,1){a}
			\tkzDefPoint(4,0){b}
			\tkzDefPoint(5,0){c}
			\tkzDefPoint(5,1){d}
			\tkzDrawPolygon[line width = 0.7mm, color = black](a,b,c,d);
			
			\tkzDefPoint(0,0){a}
			\tkzDefPoint(0,-1){b}
			\tkzDefPoint(1,-1){c}
			\tkzDefPoint(1,0){d}
			\tkzDrawPolygon[line width = 0.7mm, color = black](a,b,c,d);
			
			\tkzDefPoint(1,0){a}
			\tkzDefPoint(1,-1){b}
			\tkzDefPoint(2,-1){c}
			\tkzDefPoint(2,0){d}
			\tkzDrawPolygon[line width = 0.7mm, color = black](a,b,c,d);
			
			\tkzDefPoint(2,0){a}
			\tkzDefPoint(2,-1){b}
			\tkzDefPoint(3,-1){c}
			\tkzDefPoint(3,0){d}
			\tkzDrawPolygon[line width = 0.7mm, color = black, fill=red!30](a,b,c,d);
			
			\tkzDefPoint(3,0){a}
			\tkzDefPoint(3,-1){b}
			\tkzDefPoint(4,-1){c}
			\tkzDefPoint(4,0){d}
			\tkzDrawPolygon[line width = 0.7mm, color = black](a,b,c,d);
			
			\tkzDefPoint(4,-1){a}
			\tkzDefPoint(4,0){b}
			\tkzDefPoint(5,0){c}
			\tkzDefPoint(5,-1){d}
			\tkzDrawPolygon[line width = 0.7mm, color = black](a,b,c,d);
			
			\tkzDefPoint(0,-2){a}
			\tkzDefPoint(0,-1){b}
			\tkzDefPoint(1,-1){c}
			\tkzDefPoint(1,-2){d}
			\tkzDrawPolygon[line width = 0.7mm, color = black](a,b,c,d);
			
			\tkzDefPoint(1,-2){a}
			\tkzDefPoint(1,-1){b}
			\tkzDefPoint(2,-1){c}
			\tkzDefPoint(2,-2){d}
			\tkzDrawPolygon[line width = 0.7mm, color = black](a,b,c,d);
			
			\tkzDefPoint(2,-2){a}
			\tkzDefPoint(2,-1){b}
			\tkzDefPoint(3,-1){c}
			\tkzDefPoint(3,-2){d}
			\tkzDrawPolygon[line width = 0.7mm, color = black](a,b,c,d);
			
			\tkzDefPoint(3,-2){a}
			\tkzDefPoint(3,-1){b}
			\tkzDefPoint(4,-1){c}
			\tkzDefPoint(4,-2){d}
			\tkzDrawPolygon[line width = 0.7mm, color = black](a,b,c,d);
			
			\tkzDefPoint(0,-3){a}
			\tkzDefPoint(0,-2){b}
			\tkzDefPoint(1,-2){c}
			\tkzDefPoint(1,-3){d}
			\tkzDrawPolygon[line width = 0.7mm, color = black](a,b,c,d);
			
			\node at (0.5,1.5){{\Large $9$}};
			\node at (1.5,1.5){{\Large $7$}};
			\node at (2.5,1.5){{\Large $6$}};
			\node at (3.5,1.5){{\Large $5$}};
			\node at (4.5,1.5){{\Large $3$}};
			\node at (-0.5,0.5){{\Large $1$}};
			\node at (-0.5,-0.5){{\Large $2$}};
			\node at (-0.5,-1.5){{\Large $4$}};
			\node at (-0.5,-2.5){{\Large $8$}};
	\end{tikzpicture}}
	\caption{\label{fig:BElambda} The (labelled) integer partition $\oplamb(\L,\R)$ with $\L = \{1,2,4,8\}$ and $\R = \{3,5,6,7,9\}$. It also corresponds to $\lambda(c)$ with $c = (1,2,4,8,9,7,6,5,3)$. The red box has $c$-coordinates $\LR{2,6}_c$.}
\end{figure}

Given a Coxeter element $c \in \mathfrak{S}_{n+1}$, we write $\oplamb(c)$ for the integer partition $\oplamb(\{1\} \cup \L_c, \R_c \cup \{n+1\})$. Thanks to the observation above, we introduce the
\new{$c$-coordinates} of any box in $\Fer(\oplamb(c))$ as it follows. By setting $\L= \L_c \cup \{1\} = \{\ell_1<\ldots<\ell_p\}$ and $\R = \R_c \cup \{n+1\} = \{r_1 < \ldots < r_q\}$, we write $\LR{\ell_i, r_{q-j+1}}_c $ for the box $(i,j) \in \Fer(\oplamb(c))$ whenever $\ell_i < r_{q-j+1}$. See \cref{fig:BElambda} for an example of such an object.

\begin{proposition}\label{prop:elemintbipartintpartitions}
	For any integer partition $\lambda$, there exists a unique elementary interval bipartition $(\L,\R)$ such that $\lambda = \oplamb(\L,\R)$.
\end{proposition}

\begin{proof}
	If $\lambda = (0)$, then we set $\L = \R = \varnothing$ and we are done.
	Otherwise, we label the segments of the southeast border of the shape of $\Fer(\lambda)$ from $1$ to its length, going from the top-right to the bottom-left. This defines a label for each row and each column of $\Fer(\lambda)$. We set $\L$, the set of labels assigned to the rows, and $\R$, the set of labels assigned to the columns. We can easily check that $(\L,\R)$ is an elementary interval bipartition (of the interval $\llrr{1,h_\lambda(1,1)+1}$). By construction, it is unique.
\end{proof}

\begin{corollary}
	Let $n \in \mathbb{N}^*$. For any $\lambda \in \Hk_n$, there exists a unique Coxeter element $c$ of $\mathfrak{S}_{n+1}$ such that $\oplamb(c) = \lambda$.
\end{corollary}

\begin{proof}
	It follows automatically from \cref{cor:Coxbipart} and \cref{prop:elemintbipartintpartitions}.
\end{proof}

Given a $\lambda \in \Hk_n$, we denote by $\opc(\lambda)$ the unique Coxeter element of $\mathfrak{S}_n$ such that $\oplamb(\opc(\lambda)) = \lambda$. 

\begin{lemma} \label{lem:diagonalboxes} Let $c \in \mathfrak{S}_{n+1}$. For $k \in \{1,\ldots,n\}$, we have:
	\[\#D_k(\oplamb(c)) = \min(\#\{\ell \in \L_{c} \cup \{1\} \mid \ell \leqslant k\}, \#\{r \in \R_c \cup \{n+1\} \mid r > k\})\]
\end{lemma}

\begin{proof}
	This result follows by interpreting $D_k(\lambda)$ as "the diagonal" of the rectangle made of boxes $(i,j)$ such that $\ell_i \leqslant k < r_{\lambda_1 -j+1}$. 
\end{proof}

\subsection{Auslander--Reiten quiver}

Let $c \in \mathfrak{S}_{n+1}$ be a Coxeter element. The \new{Auslander--Reiten quiver of $c$}, denoted by $\AR(c)$, is the oriented graph satisfying the following conditions:
\begin{enumerate}[label = $\bullet$]
	\item The vertices of $\AR(c)$ are the transpositions $(i,j)$, with $i<j$, in $\mathfrak{S}_{n+1}$;
	\item The arrows of $\AR(c)$ are given, for all $i < j$, by
	\begin{enumerate}[label = $\bullet$]
		\item $(i,j) \longrightarrow (i,c(j))$ whenever $i < c(j)$;
		
		\item $(i,j) \longrightarrow (c(i),j)$ whenever  $c(i) < j$.
	\end{enumerate}
\end{enumerate}

Let us state an obvious and useful proposition about these quivers.

\begin{proposition} \label{prop:AcycAR}
	For any Coxeter element $c \in \mathfrak{S}_{n+1}$, the Auslander--Reiten quiver $\AR(c)$ is an acyclic connected directed graph. Moreover:
	\begin{enumerate}[label = $\bullet$, itemsep=1mm]
		\item  its sources are the initial adjacent transpositions in $c$, and
		\item its sinks are the final adjacent transpositions in $c$.
	\end{enumerate}
\end{proposition}

\begin{figure}[!ht]
	\centering
	\[
	\scalebox{0.9}{
		\begin{tikzpicture}

			\node (67) at (0,2) {(67)};
			\node (57) at (1,3) {(57)};
			\node (27) at (2,4) {(27)};
			\node (17) at (3,5) {(17)};
			\node (37) at (4,6) {(37)};
			\node (47) at (5,7) {(47)};
			
			\node (69) at (1,1) {(69)};
			\node (59) at (2,2) {(59)};
			\node (29) at (3,3) {(29)};
			\node (19) at (4,4) {(19)};
			\node (39) at (5,5) {(39)};
			
			\node (49) at (6,6) {(49)};
			
			\node (68) at (2,0) {(68)};
			\node (58) at (3,1) {(58)};
			\node (28) at (4,2) {(28)};
			\node (18) at (5,3) {(18)};
			\node (38) at (6,4) {(38)};
			\node (48) at (7,5) {(48)};
			\node (78) at (8,6) {(78)};
			\node (89) at (0,0) {(89)};
			\node (79) at (7,7) {(79)};
			
			\node (56) at (4,0) {(56)};
			\node (25) at (6,0) {(25)};
			\node (12) at (8,0) {(12)};
			\node (13) at (1,7) {(13)};
			\node (34) at (3,7) {(34)};
			
			\node (26) at (5,1) {(26)};
			\node (15) at (7,1) {(15)};
			\node (23) at (0,6) {(23)};
			\node (14) at (2,6) {(14)};
			
			\node (16) at (6,2) {(16)};
			\node (35) at (8,2) {(35)};
			\node (24) at (1,5) {(24)};
			
			\node (36) at (7,3) {(36)};
			\node (45) at (9,3) {(45)};
			
			\node (46) at (8,4) {(46)};
			
			\draw[->] (13)--(14);\draw[->] (14)--(34);\draw[->] (14)--(17);\draw[->] (15)--(35);\draw[->] (15)--(12);\draw[->] (16)--(36);\draw[->] (16)--(15);\draw[->] (17)--(37);\draw[->] (17)--(19);\draw[->] (18)--(38);\draw[->] (18)--(16);\draw[->] (19)--(39);\draw[->] (19)--(18);\draw[->] (23)--(13);\draw[->] (23)--(24);\draw[->] (24)--(14);\draw[->] (24)--(27);\draw[->] (25)--(15);\draw[->] (26)--(16);\draw[->] (26)--(25);\draw[->] (27)--(17);\draw[->] (27)--(29);\draw[->] (28)--(18);\draw[->] (28)--(26);\draw[->] (29)--(19);\draw[->] (29)--(28);\draw[->] (34)--(37);\draw[->] (35)--(45);\draw[->] (36)--(46);\draw[->] (36)--(35);\draw[->] (37)--(47);\draw[->] (37)--(39);\draw[->] (38)--(48);\draw[->] (38)--(36);\draw[->] (39)--(49);\draw[->] (39)--(38);\draw[->] (46)--(45);\draw[->] (47)--(49);\draw[->] (48)--(78);\draw[->] (48)--(46);\draw[->] (49)--(79);\draw[->] (49)--(48);\draw[->] (56)--(26);\draw[->] (57)--(27);\draw[->] (57)--(59);\draw[->] (58)--(28);\draw[->] (58)--(56);\draw[->] (59)--(29);\draw[->] (59)--(58);\draw[->] (67)--(57);\draw[->] (67)--(69);\draw[->] (68)--(58);\draw[->] (69)--(59);\draw[->] (69)--(68);\draw[->] (79)--(78);\draw[->] (89)--(69);
	\end{tikzpicture} }\]
	\caption{\label{fig:ARc} The Auslander--Reiten quiver of $c = (1,3,4,7,9,8,6,5,2)=s_2s_1s_3s_6s_5s_4s_8s_7$.}
\end{figure}
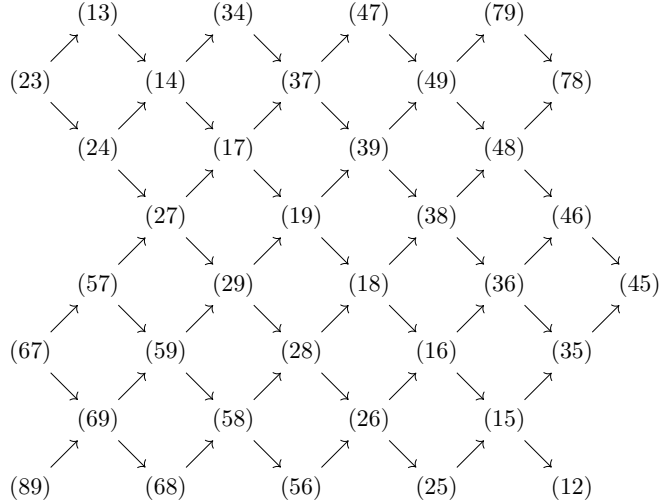

To construct such a graph recursively, we can first find the initial adjacent transpositions of $c$, which are all the sources,  and step by step, using the second rule, construct the arrows and the vertices of $\AR(c)$ until we reach all the transpositions of $\mathfrak{S}_{n+1}$. Note that the sinks of $\AR(c)$ are given by the final adjacent transpositions of $c$.
See \cref{fig:ARc} for an explicit example.

Moreover, one can notice that we can construct $\AR(c)$ from any transposition $(i,j) \in \mathfrak{S}_{n+1}$ using the second rule.

\begin{proposition} \label{prop:ARandinverse}
	Let $c \in \mathfrak{S}_{n+1}$ be a Coxeter element. Then $\AR(c^{-1}) = \AR(c)^{\op}$.
\end{proposition}

\begin{remark}
	The Auslander--Reiten quiver of a Coxeter element has a representation-theoretic meaning for the quiver $Q(c)^{\op}$ (see \cref{ss:quiver}).
\end{remark}

	\section{Storability}
	\label{sec:stor}
	In this section, we first recall the notion of storability, introduced in \cite{Deq23}, and we enumerate a few primary results. Then, given a positive integer $n \in \mathbb{N}^*$ and a Coxeter element $c \in \mathfrak{S}_{n+1}$, we introduce the notion of $c$-storabilty for $n$-tuples of integer partitions. We highlight their bijective link with the reverse plane partitions of $\oplamb(c)$.

\subsection{Storable pairs and storable triplets}
\label{ss:storpairstriplets} 

\begin{definition}\label{def:storpairs}
	Let $\lambda$ and $\mu$ be two integer partitions. The pair $(\lambda, \mu)$ is  \new{storable} whenever for all $i \in \mathbb{N}^*$, $ \lambda_i \geqslant \mu_i \geqslant \lambda_{i+1}$ (we can add zero parts if needed). Such a pair is \new{strongly storable} if we have $\lambda_1 = \mu_1$. 
\end{definition}

We can picture storable pairs as follows. See a partition $\lambda$ as a right-infinite row of forty-five-degree rotated squares filled with the parts of $\lambda$ from left to right. We can add infinitely many zeros to the right. See $\mu$ in the same way. We say that two such rows of squares are \emph{intertwining} if, for all $i \geqslant 1$, the $i$th square of one row is placed between the $i$th and the $(i+1)$th squares of the other row.

\begin{figure}[!ht]
	\centering
	\begin{tikzpicture}[scale=0.7]
		\draw[line width=0.7mm](-1,1) edge (1,3);
		\draw[line width=0.7mm](0,0) edge (3,3);
		\draw[line width=0.7mm](2,0) edge (5,3);
		\draw[line width=0.7mm](4,0) edge (6,2);
		
		\draw[line width=0.7mm](-1,1) edge (0,0);
		\draw[line width=0.7mm](0,2) edge (2,0);
		\draw[line width=0.7mm](1,3) edge (4,0);
		\draw[line width=0.7mm](3,3) edge (5,1);
		\draw[line width=0.7mm](5,3) edge (6,2);
		
		\draw[line width=0.7mm,dotted](6.5,2) edge (8,2);
		\draw[line width=0.7mm,dotted](5.5,1) edge (7,1);
		
		\tkzLabelPoint[below](0.05,1.55){{\large $\color{blue}{\lambda_1}$}}
		\tkzLabelPoint[below](2.05,1.55){{\large $\color{blue}{\lambda_2}$}}
		\tkzLabelPoint[below](4.05,1.55){{\large $\color{blue}{\lambda_3}$}}
		
		\tkzLabelPoint[below](1.05,2.45){{\large $\color{red}{\mu_1}$}}
		\tkzLabelPoint[below](3.05,2.45){{\large $\color{red}{\mu_2}$}}
		\tkzLabelPoint[below](5.05,2.45){{\large $\color{red}{\mu_3}$}}
	\end{tikzpicture}
	\caption{\label{fig:store} Illustration of storability of $(\lambda, \mu)$.}
\end{figure}
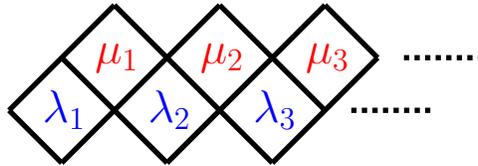

Then $(\lambda, \mu)$ is a storable pair if and only if we can intertwine the two rows of filled squares such that, when we read the two rows together from left to right, the values remain decreasing (\cref{fig:store}). In other words, the square $\mu_i$ intertwines the squares $\lambda_i$ and $\lambda_{i+1}$ whenever $\lambda_i \geqslant \mu_i \geqslant \lambda_{i+1}$.

We give two results that arise from the definition.

\begin{lemma} \label{lem:elementstorable} Let $\lambda$ and $\mu$ be two integer partitions. 
	\begin{enumerate}[label = \arabic*)]
		\item \label{estor1} If $(\lambda, \mu)$ and $(\mu, \lambda)$ are both storable, then $\lambda = \mu$;
		
		\item \label{estor2} If $(\lambda, \mu)$ is storable, then $\ell(\lambda) \in \{\ell(\mu), \ell(\mu) + 1\}$.
	\end{enumerate}
\end{lemma}

\begin{definition} \label{def:storabletriplets} Let $\lambda, \mu$ and $\nu$ be three integer partitions. The triplet $(\lambda, \mu, \nu)$ is \new{storable} if the two following conditions are satisfied: \begin{enumerate}[label = $\bullet$]
		
		\item either $(\lambda, \mu)$ or $(\mu, \lambda)$ is a storable pair;
		
		\item either $(\mu, \nu)$ or $(\nu, \mu)$ is a storable pair.
		
	\end{enumerate}
	More precisely, we say that $(\lambda, \mu, \nu)$ is:
	\begin{enumerate}
		
		\item[$(\boxplus \boxplus)$]  $(\boxplus,\boxplus)$-storable if $(\lambda, \mu)$ and $(\nu, \mu)$ are storable pairs;
		
		\item[$(\boxplus \boxminus)$]  $(\boxplus,\boxminus)$-storable if $(\lambda, \mu)$ and $(\mu, \nu)$ are storable pairs;
		
		\item[$(\boxminus \boxplus)$]  $(\boxminus,\boxplus)$-storable if $(\mu, \lambda)$ and $(\nu, \mu)$ are storable pairs;
		
		\item[$(\boxminus \boxminus)$] $(\boxminus,\boxminus)$-storable if $(\mu, \lambda)$ and $(\mu, \nu)$ are storable pairs.
	\end{enumerate}
	Such a triplet is \new{strongly storable} whenever $\lambda_1 = \mu_1$ or $\mu_1 = \nu_1$. 
\end{definition}

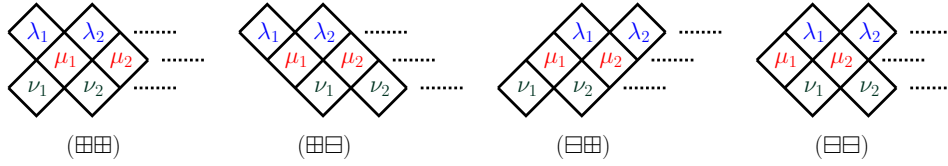
\begin{figure}[!ht]
	\centering
	
	\scalebox{0.53}{
		\begin{tikzpicture}[scale=0.7]
			\draw[line width=0.7mm](0,2) -- (1,3) -- (4,0) ;
			\draw[line width=0.7mm](0,2) -- (3,-1) -- (5,1) ;
			\draw[line width=0.7mm](0,0) -- (1,-1) -- (4,2) ;
			\draw[line width=0.7mm](0,0) -- (3,3) -- (5,1);
			
			\draw[line width=0.7mm,dotted](4.5,2) edge (6,2);
			\draw[line width=0.7mm,dotted](5.5,1) edge (7,1);
			\draw[line width=0.7mm,dotted](4.5,0) edge (6,0);
			
			\tkzLabelPoint[below](2.05,1.45){{\large $\color{red}{\mu_1}$}}
			\tkzLabelPoint[below](4.05,1.45){{\large $\color{red}{\mu_2}$}}

			\tkzLabelPoint[below](1.05,2.55){{\large $\color{blue}{\lambda_1}$}}
			\tkzLabelPoint[below](3.05,2.55){{\large $\color{blue}{\lambda_2}$}}
			
			\tkzLabelPoint[below](1.05,0.4){{\large $\color{darkgreen}{\nu_1}$}}
			\tkzLabelPoint[below](3.05,0.4){{\large $\color{darkgreen}{\nu_2}$}}
			
			\tkzLabelPoint[below](3.05,-1.5){{\large $(\boxplus \boxplus)$}}
		\end{tikzpicture} \qquad 
		\begin{tikzpicture}[scale=0.7]
			\draw[line width=0.7mm](-1,1) -- (1,3) -- (4,0) ;
			\draw[line width=0.7mm](-2,2) -- (-1,3) -- (3,-1) -- (4,0) ;
			\draw[line width=0.7mm](-2,2) -- (1,-1) -- (3,1) ;
			\draw[line width=0.7mm](0,0) -- (2,2);
			
			\draw[line width=0.7mm,dotted](2.5,2) edge (4,2);
			\draw[line width=0.7mm,dotted](3.5,1) edge (5,1);
			\draw[line width=0.7mm,dotted](4.5,0) edge (6,0);
			
			\tkzLabelPoint[below](0.05,1.45){{\large $\color{red}{\mu_1}$}}
			\tkzLabelPoint[below](2.05,1.45){{\large $\color{red}{\mu_2}$}}

			\tkzLabelPoint[below](-.95,2.55){{\large $\color{blue}{\lambda_1}$}}
			\tkzLabelPoint[below](1.05,2.55){{\large $\color{blue}{\lambda_2}$}}
			
			\tkzLabelPoint[below](1.05,0.4){{\large $\color{darkgreen}{\nu_1}$}}
			\tkzLabelPoint[below](3.05,0.4){{\large $\color{darkgreen}{\nu_2}$}}
			
			\tkzLabelPoint[below](1.05,-1.5){{\large $(\boxplus \boxminus)$}}
		\end{tikzpicture} \qquad \begin{tikzpicture}[scale=0.7]
			\draw[line width=0.7mm](-2,0) -- (1,3) -- (3,1) ;
			\draw[line width=0.7mm](0,2) -- (2,0);
			\draw[line width=0.7mm](-1,1) -- (1,-1) -- (4,2) ;
			\draw[line width=0.7mm](-2,0) --(-1,-1) -- (3,3) -- (4,2);
			
			\draw[line width=0.7mm,dotted](4.5,2) edge (6,2);
			\draw[line width=0.7mm,dotted](3.5,1) edge (5,1);
			\draw[line width=0.7mm,dotted](2.5,0) edge (4,0);
			
			\tkzLabelPoint[below](0.05,1.45){{\large $\color{red}{\mu_1}$}}
			\tkzLabelPoint[below](2.05,1.45){{\large $\color{red}{\mu_2}$}}

			\tkzLabelPoint[below](1.05,2.55){{\large $\color{blue}{\lambda_1}$}}
			\tkzLabelPoint[below](3.05,2.55){{\large $\color{blue}{\lambda_2}$}}
			
			\tkzLabelPoint[below](-0.95,0.4){{\large $\color{darkgreen}{\nu_1}$}}
			\tkzLabelPoint[below](1.05,0.4){{\large $\color{darkgreen}{\nu_2}$}}
			
			\tkzLabelPoint[below](1.05,-1.5){{\large $(\boxminus \boxplus)$}}
		\end{tikzpicture}  \qquad  \begin{tikzpicture}[scale=0.7]
			\draw[line width=0.7mm](-1,1) -- (1,3) -- (4,0) ;
			\draw[line width=0.7mm](0,2) -- (3,-1) -- (4,0) ;
			\draw[line width=0.7mm](-1,1) -- (1,-1) -- (4,2) ;
			\draw[line width=0.7mm](0,0) -- (3,3) -- (4,2);
			
			\draw[line width=0.7mm,dotted](4.5,2) edge (6,2);
			\draw[line width=0.7mm,dotted](3.5,1) edge (5,1);
			\draw[line width=0.7mm,dotted](4.5,0) edge (6,0);
			
			\tkzLabelPoint[below](0.05,1.45){{\large $\color{red}{\mu_1}$}}
			\tkzLabelPoint[below](2.05,1.45){{\large $\color{red}{\mu_2}$}}

			\tkzLabelPoint[below](1.05,2.55){{\large $\color{blue}{\lambda_1}$}}
			\tkzLabelPoint[below](3.05,2.55){{\large $\color{blue}{\lambda_2}$}}
			
			\tkzLabelPoint[below](1.05,0.4){{\large $\color{darkgreen}{\nu_1}$}}
			\tkzLabelPoint[below](3.05,0.4){{\large $\color{darkgreen}{\nu_2}$}}
			
			\tkzLabelPoint[below](2.05,-1.5){{\large $(\boxminus \boxminus)$}}
	\end{tikzpicture}}

	\caption{\label{fig:AllConf} Illustration of the four storability configurations of $(\lambda, \mu, \nu)$.}
\end{figure}

We illustrate the four storability configurations in \cref{fig:AllConf}.

We introduce the following notion via the representation-theoretic interpretation introduced in \cref{ss:refl}.

\begin{definition}\label{def:diag}
	Let $\lambda, \mu$, and $\nu$ be three integer partitions. Assume that $(\lambda, \mu, \nu)$ is a storable triplet. We define the \new{diagonal transformation of $\mu$ in $(\lambda, \mu, \nu)$}, denoted $\diag(\lambda, \mu, \nu)$, to be the integer partition $\theta = (\theta_1, \theta_2, \ldots)$ such that:
	
	\begin{enumerate}[label = $\bullet$]
		
		\item if $(\lambda, \mu, \nu)$ is $(\boxplus, \boxplus)$-storable, then we define, for all $i \geqslant 1$, $$\theta_i =  \begin{cases}
			\max(\lambda_1, \nu_1) & \text{if } i = 1 \\
			\min(\lambda_{i-1},\nu_{i-1}) + \max(\lambda_i, \nu_i) - \mu_{i-1} & \text{otherwise;}\\
		\end{cases} $$
		
		\item if $(\lambda, \mu, \nu)$ is $(\boxplus, \boxminus)$-storable, then we define, for all $i \geqslant 1$, $$ \theta_i = \begin{cases} \lambda_1 + \max(\lambda_2, \nu_1) - \mu_1 & \text{if } i = 1\\
			\min(\lambda_i,\nu_{i-1}) + \max(\lambda_{i+1}, \nu_i) - \mu_{i}& \text{otherwise;}
		\end{cases}$$
		
		\item if $(\lambda, \mu, \nu)$ is $(\boxminus, \boxplus)$-storable, then we define, for all $i \geqslant 1$, $$ \theta_i = \begin{cases} \nu_1 + \max(\lambda_1, \nu_2) - \mu_1 & \text{if } i = 1\\
			\min(\lambda_{i-1},\nu_{i}) + \max(\lambda_{i}, \nu_{i+1}) - \mu_{i}& \text{otherwise;}
		\end{cases}$$
		
		\item if $(\lambda, \mu, \nu)$ is  $(\boxminus, \boxminus)$-storable, then we define, for all $i \geqslant 1$, $$\theta_i =  \min(\lambda_i, \nu_i) + \max(\lambda_{i+1}, \nu_{i+1}) - \mu_{i+1}.$$
	\end{enumerate}
\end{definition}

\begin{figure}[!ht]
	\centering
	
	\scalebox{0.6}{
		\begin{tikzpicture}[scale=1.3]
			
			\tkzDefPoint(1,1){a}
			\tkzDefPoint(2,2){b}
			\tkzDefPoint(3,1){c}
			\tkzDefPoint(2,0){d}
			\tkzDrawPolygon[line width = 2mm, color = black, fill = black!10](a,b,c,d);
			
			\draw[line width=0.7mm](0,2) -- (1,3) -- (4,0) ;
			\draw[line width=0.7mm](0,2) -- (3,-1) -- (4,0) ;
			\draw[line width=0.7mm](0,0) -- (1,-1) -- (4,2) ;
			\draw[line width=0.7mm](0,0) -- (3,3) -- (4,2);
			
			\draw[line width=0.7mm,dotted](4.5,2) edge (5.5,2);
			\draw[line width=0.7mm,dotted](3.5,1) edge (4.5,1);
			\draw[line width=0.7mm,dotted](4.5,0) edge (5.5,0);
			
			\draw[line width=0.7mm,dotted](-0.5,2) edge node[above]{\Large $\color{blue}{\lambda}$} (-1.5,2);
			\draw[line width=0.7mm,dotted](0.5,1) edge node[above]{\Large $\color{red}{\mu}$}(-0.5,1);
			\draw[line width=0.7mm,dotted](-0.5,0) edge node[above]{\Large $\color{darkgreen}{\nu}$} (-1.5,0);
			
			\draw[->,>= angle 60, line width=2mm](5,1) -- node[midway, xshift=-9pt, yshift=15pt]{\Huge $\diag$} (7,1);
			
			\tkzLabelPoint[below](2,1.27){{\Huge $\color{red}{e}$}}

			\tkzLabelPoint[below](1,2.3){{\Huge $\color{blue}{a}$}}
			\tkzLabelPoint[below](3,2.4){{\Huge $\color{blue}{b}$}}
			
			\tkzLabelPoint[below](1,0.27){{\Huge $\color{darkgreen}{c}$}}
			\tkzLabelPoint[below](3,0.4){{\Huge $\color{darkgreen}{d}$}}
			
			\begin{scope}[xshift = 8cm]
				
				\tkzDefPoint(1,1){a}
				\tkzDefPoint(2,2){b}
				\tkzDefPoint(3,1){c}
				\tkzDefPoint(2,0){d}
				\tkzDrawPolygon[line width = 2mm, color = black, fill = black!10](a,b,c,d);
				
				\draw[line width=0.7mm](0,2) -- (1,3) -- (4,0) ;
				\draw[line width=0.7mm](0,2) -- (3,-1) -- (4,0) ;
				\draw[line width=0.7mm](0,0) -- (1,-1) -- (4,2) ;
				\draw[line width=0.7mm](0,0) -- (3,3) -- (4,2);
				
				\draw[line width=0.7mm,dotted](4.5,2) edge (5.5,2);
				\draw[line width=0.7mm,dotted](3.5,1) edge (4.5,1);
				\draw[line width=0.7mm,dotted](4.5,0) edge (5.5,0);
				
				\draw[line width=0.7mm,dotted](-0.5,2) edge (-1.5,2);
				\draw[line width=0.7mm,dotted](0.5,1) edge (-0.5,1);
				\draw[line width=0.7mm,dotted](-0.5,0) edge (-1.5,0);

				\tkzLabelPoint[below](2,1.5){{\small$\color{red}{\min(a,c)}$}}
				\tkzLabelPoint[below](2,1.15){{\small$\color{red}{+\max(b,d)}$}}
				\tkzLabelPoint[below](2,0.8){{\small$\color{red}{-e}$}}

				\tkzLabelPoint[below](1,2.3){{\Huge $\color{blue}{a}$}}
				\tkzLabelPoint[below](3,2.4){{\Huge $\color{blue}{b}$}}
				
				\tkzLabelPoint[below](1,0.27){{\Huge $\color{darkgreen}{c}$}}
				\tkzLabelPoint[below](3,0.4){{\Huge $\color{darkgreen}{d}$}}
				
			\end{scope};
	\end{tikzpicture}}

	\caption{\label{fig:GenDiagop} Illustration of the local operations to calculate $\diag(\lambda, \mu, \nu)$.}
\end{figure}
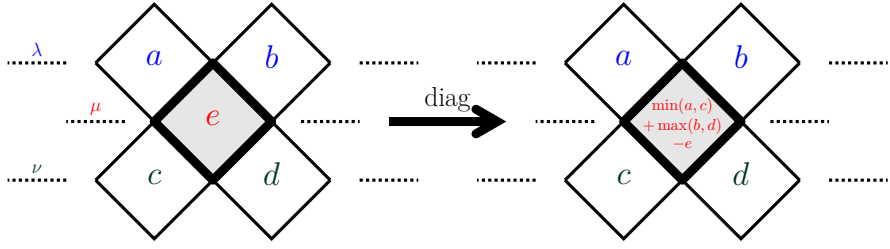

We can picture the diagonal operation as doing local operations for each square of $\mu$ in the diagram representing the storable triple $(\lambda, \mu, \nu)$ (\cref{fig:GenDiagop}).

Remark that $\lambda$ and $\nu$ play symmetric roles: $\diag(\lambda, \mu, \nu) = \diag(\nu, \mu, \lambda)$. Here are some elementary statements about the diagonal transformation.

\begin{lemma}\label{lem:elementary diag store} Let $\lambda, \mu$ and $\nu$ be three integer partitions. When it is well-defined,  consider $\theta = \diag(\lambda, \nu, \mu)$.
	\begin{enumerate}[label = \arabic*),itemsep=1mm]
		
		\item \label{ediag0} If $(\lambda, \mu)$ is a storable pair, then $\diag(\lambda,\mu, \mu) = \lambda$.
		
		\item \label{ediag1}  If $(\lambda, \mu, \nu)$ is $(\boxplus, \boxplus)$-storable, then $(\lambda, \theta, \nu)$ is strongly $(\boxminus, \boxminus)$-storable.
		
		\item \label{ediag2}  If $(\lambda, \mu, \nu)$ is $(\boxplus,\boxminus)$-storable, then $(\lambda, \theta, \nu)$ is $(\boxplus, \boxminus)$-storable. 
		
		\item \label{ediag3} If $(\lambda, \mu, \nu)$ is $(\boxminus,\boxplus)$-storable, then $(\lambda, \theta, \nu)$ is $(\boxminus, \boxplus)$-storable.  
		
		\item \label{ediag4} If $(\lambda, \mu, \nu)$ is $(\boxminus,\boxminus)$-storable, then $(\lambda, \theta, \nu)$ is $(\boxplus, \boxplus)$-storable. 
		
		\item \label{ediag5} If $(\lambda, \mu, \nu)$ is $(\boxplus, \boxplus)$-storable, $(\boxplus, \boxminus)$-storable, $(\boxminus, \boxplus)$-storable or strongly $(\boxminus, \boxminus)$-storable, then $\diag(\lambda, \theta, \nu) = \mu$.
	\end{enumerate}
\end{lemma}

\subsection{\texorpdfstring{$c$}{c}-storability}
\label{ss:cstor}

For $\A \subset \mathbb{Z}$ and $j \in \mathbb{Z}$, we set $\A[j] = \{a+j \mid a \in A\}$.

Let $n \in \mathbb{N}^*$.  Let $\pmb{\pi} = (\pi^k)_{1 \leqslant k \leqslant n}$ a $n+2$-tuple of integer partitions. In the following, we set $\pi^0 = \pi^{n+1} = (0)$ to make things more convenient with the following definitions. We say that $\pmb{\pi}$ is \new{(strongly) $(\boxminus, \boxminus)$-storable at $k$} if  $(\pi^{k-1}, \pi^k, \pi^{k+1})$ is a (strongly) $(\boxminus, \boxminus)$-storable triplet. We use the same formulation for the three other storability configurations.

\begin{definition} \label{def:BEstore}
	Let $(\L,\R)$ be a pair of disjoint subsets of $\{1,\ldots,n+1\}$ such that $\min(\L,\R) \in \L$ and $\max(\L,\R) \in \R$. A $n$-tuple of integer partitions $\pmb{\pi} = (\pi^k)_{1 \leqslant k \leqslant n}$ is said to be $(\L,\R)$-storable whenever the following assertions hold:
	\begin{enumerate}[label = $(\alph*)$]
		\item $\pmb{\pi}$ is $(\boxplus, \boxplus)$-storable at $i$ for $i \notin \L\cup \R[-1]$;
		
		\item $\pmb{\pi}$ is $(\boxplus, \boxminus)$-storable at $i$ for  $i \in \R[-1] \setminus \L$;
		
		\item $\pmb{\pi}$ is $(\boxminus, \boxplus)$-storable at $i$ if $i\in \L \setminus \R[-1]$;
		
		\item $\pmb{\pi}$ is $(\boxminus, \boxminus)$-storable at $i$ for $i \in \L \cap \R[-1]$.
		
		\item $\pi^k=(0)$ for $k \notin \{\min(\L), \ldots, \max(\R)\}$.
	\end{enumerate}
	Let $c \in \mathfrak{S}_{n+1}$ be a Coxeter element.  A $n$-tuple of integer partitions $\pmb{\pi}$ is \new{$c$-storable} whenever $\pmb{\pi}$ is $(\L_c\cup\{1\},\R_c \cup\{n+1\})$-storable.
\end{definition}

Note that $(a)$ and $(e)$ are not exclusive.

\begin{example} \label{ex:BEstorex}
	Fix $n=5$ and the Coxeter element $c = (1,2,4,6,5,3) \in \mathfrak{S}_n$. The tuple $\pmb{\pi} = ((2), (4,2),(2,1),(3,2),(3))$ is $c$-storable. From the drawing (\cref{fig:BEstorex}), we can determine the $c$-storability property.
	
	\begin{figure}[!ht]
		\centering
		\begin{tikzpicture}[scale=0.6]
			
			\foreach \y in {-1,...,3}{
				\FPeval{\res}{clip(4-\y)}
				\draw [line width=0.7mm, black!10, dashed] (-0.5,\y) --  (5,\y);
				\node[black] at (-0.9,\y){$\mathbf{\res}$};
			}
			
			\draw[line width=0.7mm](0,2) -- (1,3) -- (4,0) ;
			\draw[line width=0.7mm](0,2) -- (3,-1) -- (4,0) ;
			\draw[line width=0.7mm](0,0) -- (1,-1) -- (4,2) ;
			\draw[line width=0.7mm](0,0) -- (3,3) -- (4,2);
			\draw[line width=0.7mm](1,-1) -- (2,-2) -- (3,-1);
			\draw[line width=0.7mm](1,3) -- (2,4) -- (3,3);
			\draw[line width=0.7mm](4,2) -- (5,1) -- (4,0);

			\tkzLabelPoint[below](2,3.5){{\Large $2$}}
			
			\tkzLabelPoint[below](1,2.5){{\Large $4$}}
			\tkzLabelPoint[below](3,2.5){{\Large $2$}}
			
			\tkzLabelPoint[below](2,1.5){{\Large $2$}}
			\tkzLabelPoint[below](4,1.5){{\Large $1$}}
			
			\tkzLabelPoint[below](1,0.5){{\Large $3$}}
			\tkzLabelPoint[below](3,0.5){{\Large $2$}}
			
			\tkzLabelPoint[below](2,-0.5){{\Large $3$}}
			
		\end{tikzpicture}
		\caption{\label{fig:BEstorex} Illustration of the $c$-storability of $\pmb{\pi}$ in \cref{ex:BEstorex}.}
	\end{figure}
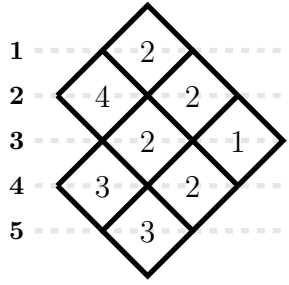
\end{example}

Given a Coxeter element $c \in \mathfrak{S}_{n+1}$, we denote by $\st(c)$ the set of $n$-tuples of integer partitions that are $c$-storable. 

For $k \in \{1,\ldots,n\}$, we define the \new{ diagonal transformation of $\pmb{\pi}$ at $k$}, denoted by $\cdiag_k(\pmb{\pi})$, to be the $n$-tuple of integer partitions obtained from $\pmb{\pi}$ by replacing $\pi^k$ with $\diag(\pi^{k-1},\pi^k, \pi^{k+1})$. The following lemma follows from \cref{lem:elementary diag store}.
\begin{lemma} \label{lem:diagc-stor}
	Let $c \in \mathfrak{S}_{n+1}$ and $\pmb{\pi} \in \st(c)$. For any $k \in \{1,\ldots,n\}$:
	\begin{enumerate}[label=\arabic*),itemsep = 1mm]
		
		\item  if $s_k$ is initial in $c$, then then $\cdiag_k(\pmb{\pi}) \in \st(s_k c s_k)$;
		
		\item if $s_k$ is final in $c$, then $\cdiag_k(\pmb{\pi}) \in \st(s_k c s_k)$, and $\cdiag_k(\pmb{\pi})$ is strongly $(\boxminus, \boxminus)$-storable at $k$;
		
		\item otherwise $\cdiag_k(\pmb{\pi}) \in \st(c)$.
	\end{enumerate}
\end{lemma}
This result will help explore ways to realize the extended generalization of RSK we present in \cref{sec:Extension} via local transformations. See \ref{sec:Further} for additional details. In the remainder of this section, we show that we have a bijection from $\st(c)$ to  $\RPP(\oplamb(c))$. 

We first exhibit a symmetric operation. We define the \emph{reverse map} $\rev_n$ on    $\{1,\ldots,n\}$ by $\rev_n(k) = n+1-k$ for all $k \in \{1,\ldots,n\}$. It induces an action on:
\begin{enumerate}[label = $\bullet$,itemsep=1mm]
	\item subsets of $\{1,\ldots,n\}$: for any $\A \subseteq \{1,\ldots,n\}$, we set $\rev_n(\A) = \{\rev_n(a) \mid a \in \A\}$;
	\item  Coxeter elements of $\mathfrak{S}_{n+1}$: for any Coxeter element $c \in \mathfrak{S}_{n+1}$, write $\rev_{n+1}(c)$ the Coxeter element obtained by conjugating $c$ with $\rev_{n+1}$ seen as a permutation;
	\item $n$-tuples of integer partitions: given $\pmb{\pi} = (\pi^k)_{1 \leqslant k \leqslant n}$, we set $\rev_n(\pmb{\pi}) = (\pi^{\rev_n(k)})_{1 \leqslant k \leqslant n}$. 
\end{enumerate}

\begin{lemma} \label{lem:symrev}
	Consider $c \in \mathfrak{S}_{n+1}$ a Coxeter element. For any $\pmb{\pi} \in \st(c)$, we have $\rev_n(\pmb{\pi}) \in \st(\rev_{n+1}(c))$.
\end{lemma}

The following result allows us to control the length of the integer partitions that constitute any $c$-storable $n$-tuple of integer partitions.

\begin{proposition} \label{prop:boundedlength}
	Let $c \in \mathfrak{S}_{n+1}$ be a Coxeter element. For all $\pmb{\pi} = (\pi^k)_{1 \leqslant k \leqslant n} \in \st(c)$, and for all $1 \leqslant k \leqslant n$, \[\ell(\pi^k) \leqslant \min \left(\#\{i \in \L_c \cup \{1\} \mid i \leqslant k\}, \#\{j \in \R_c \cup \{n+1\} \mid j > k\} \right).\]
\end{proposition}

\begin{proof} Let $c \in \mathfrak{S}_{n+1}$ and $\pmb{\pi} \in \st(c)$. We set $\L = \L_c \cup \{1\}$ and $\R = \R_c \cup \{n+1\}$.  
	
	We will prove that $\ell(\pi^k) \leqslant
	\#\{i \in \L, i \leqslant k\}$ by induction over $1 \leqslant k \leqslant n$. By hypothesis, as $1 \in \L$, we have that $(\pi^1,\pi^0=(0))$ is a storable pair. Via \cref{lem:elementstorable}, we get that $\ell(\pi^1) \leqslant \ell(\pi^0) + 1 = 1$. In addition, $\#\{i \in \L \mid i \leqslant 1\} = 1$.
	
	Assume that, for a fixed $1 \leqslant k < n$, $\ell(\pi^k) \leqslant
	\#\{i \in \L, i \leqslant k\}$. We can distinguish two cases.
	\begin{enumerate}[label = $\bullet$]
		\item If $k+1 \notin \L$, then $\pmb{\pi}$ is $(\boxplus, \boxplus)$-storable or $(\boxplus, \boxminus)$-storable at $k+1$. Thus $(\pi^k, \pi^{k+1})$ is a storable pair. By \cref{lem:elementstorable}, either $\ell(\pi^k) = \ell(\pi^{k+1})+1$ or $\ell(\pi^k) = \ell(\pi^{k+1})$. In either way, \[\ell(\pi^{i+1}) \leqslant \ell(\pi^k) \leqslant \#\{i \in \L \mid i \leqslant k \} = \#\{i \in \L \mid i \leqslant k+1\}.\]
		
		\item If $k+1 \in \L$, then $\pmb{\pi}$ is $(\boxminus, \boxplus)$-storable or $(\boxminus, \boxminus)$-storable at $k+1$. Thus $(\pi^{k+1}, \pi^k)$ is a storable pair. By \cref{lem:elementstorable}, either $\ell(\pi^{k+1}) = \ell(\pi^{i})+1$ or $\ell(\pi^{k+1}) = \ell(\pi^{i})$. Either way, 
		\[\ell(\pi^{k+1}) \leqslant \ell(\pi^k) +1 \leqslant \#\{i \in \L \mid i \leqslant k \} +1  = \#\{i \in \L \mid i \leqslant k+1\}.\]
	\end{enumerate}
	This completes the proof of the first desired inequality.
	
	To get the other one, we play with the symmetry the reverse map gives. By \cref{lem:symrev}, we know that $\rev_n(\pmb{\pi}) \in \st(\rev_{n+1}(c))$. Therefore we have, for $1 \leqslant k \leqslant n$, \[ \begin{split} \ell(\pi^k) = \ell(\rev(\pmb{\pi})^{n+1-k}) & \leqslant \{j' \in \rev_{n+1}(\R_c) \mid j' \leqslant n+2-k\} \\
		& = \{j \in \R_c \mid j > k\}. \end{split}\]
	
	We can, therefore, conclude the desired result.
\end{proof}

Now, thanks to \cref{lem:diagonalboxes}, we can define a map from $\st(c)$ to fillings of $\oplamb(c)$ as follows:
\[ \Phi_c : \left\{ \begin{matrix}
	\st(c) & \longrightarrow &\{\text{fillings of } \oplamb(c)\} \\
	\pmb{\pi} = (\pi_k)_{1 \leqslant k \leqslant n} & \longmapsto & (\Phi_c(\pmb{\pi}) : \ldiag{k,r}_{\oplamb(c)} \longmapsto \pi^{k}_{r} )
\end{matrix} \right.\] 

The map $\Phi_c$ consists on, for any $\pmb{\pi} \in \st(c)$, and $k \in \{1,\ldots,n\}$, placing the $k$th integer partition $\pi^k$ of $\pmb{\pi}$ in the boxes of $D_k(\oplamb(c))$ in the decreasing order, from right to left. We complete the diagonal with zeros if necessary. Therefore, one can check that the map $\Phi_c$ consists of a $135^\circ$-degree counterclockwise rotation of the picture of the $c$-storable $n$-tuples of integer partitions. Then, the $c$-storability conditions allow us to check that $\Phi_c(\pmb{\pi})$ is a reverse plane partition of shape $\oplamb(c)$.

\begin{theorem} \label{thm:bijLRstorRPP}
	For any $n \in \mathbb{N}^*$, and for any Coxeter element $c \in \mathfrak{S}_{n+1}$, The map $\Phi_c$ induces a bijection from $\st(c)$ to  $\RPP(\oplamb(c))$.
\end{theorem}

	\section{Overview of Jordan recoverability}
	\label{sec:JRstory}
	\subsection{Quiver representations}
\label{ss:quiver}

Let $Q=(Q_0,Q_1,s,t)$ be a finite connected quiver. Consider $\mathbb{K}$ an algebraically closed field. A \new{(finite-dimensionnal) representation of $Q$ over $\mathbb{K}$} is a pair \[E = ((E_q)_{q \in Q_0}, (E_\alpha)_{\alpha \in Q_1})\] seen as an assignement of:
\begin{enumerate}[label = $\bullet$]
	\item a finite-dimensional $\mathbb{K}$ -vector space $E_q$ to each $q \in Q_0$; 
	\item a linear transformation $E_\alpha : E_{s(\alpha)} \longrightarrow E_{t(\alpha)}$ to each $\alpha \in Q_1$.
\end{enumerate}
Write $\rep_\mathbb{K}(Q)$ for the set of the  representations of $Q$ over $\mathbb{K}$. For $E \in \rep_\mathbb{K}(Q)$, we denote by $\vdim(E) = (\dim(E_q))_{q \in Q_0}$ the \emph{dimension vector} of $E$.

Given $E,F \in \rep_\mathbb{K}(Q)$, a \new{morphism of representations} $\phi : E \longrightarrow F$ is a collection $(\phi_q)_{q \in Q_0}$ seen as an assignement of a linear transformation $\phi_q : E_q \longrightarrow F_q$ to each $q \in Q_0$ such that  for any $\alpha \in Q_1$, $F_\alpha \phi_{s(\alpha)} = \phi_{t(\alpha)} E_\alpha$. Such a morphism of representations $\phi$ is an isomorphism whenever, for all $q \in Q_0$, $\phi_q$ is an isomorphism of vector spaces. We say that two representations, $E$ and $F$, are isomorphic. We denote it by $E \cong F$ whenever an isomorphism exists from $E$ to $F$  Write $\Hom_\mathbb{K}(E,F)$ for the $\mathbb{K}$-vector space of morphisms of representations from $E$ to $F$. In particular, we set $\End_\mathbb{K}(E) = \Hom_{\mathbb{K}}(E,E)$ to refer to the endomorphisms of a representation $E$. 

We define the \emph{path algebra} $\mathbb{K}Q$ to be the $\mathbb{K}$-vector space generated as a basis by the paths over $Q$ endowed with multiplication acting on the basis as a concatenation of paths  Recall that $\rep_\mathbb{K}(Q)$ endowed with the morphisms of representations between them is a category. This category is equivalent to the category of finitely generated (right) $\mathbb{K}Q$-modules. See \cite{ASS06} for more details. 

Given $F,G \in \rep_\mathbb{K}(Q)$, we write $F \oplus G$ the direct sum of $F$ and $G$. A representation $E$ is said to be \new{indecomposable} whenever, if $E \cong F \oplus G$, then $F = 0$ or $G = 0$.  By $\ind_\mathbb{K}(Q)$, we denote the isomorphism classes of indecomposable representations of $Q$.

In our setting, we can describe indecomposable representations thanks to intervals in $\{1,\ldots,n\}$. 
Let $n \geqslant 1$ be an integer  Assume that $Q$ is an $A_n$ type quiver  Given an interval $K \subset \{1,\ldots,n\}$, we denote by $X_K$ the indecomposable representation such that:
\begin{enumerate}[label = $\bullet$]
	\item $(X_K)_q = \mathbb{K}$ if $q \in K$, and $(X_K)_q = 0$ otherwise;
	\item  $(X_K)_\alpha = \Id_\mathbb{K}$ if $s(\alpha), t(\alpha) \in K$, and $(X_K)_\alpha = 0$ otherwise.
\end{enumerate}

\begin{proposition}
	Let $Q$ be an $A_n$ type quiver  Any indecomposable representation $E \in \rep_\mathbb{K}(Q)$ is isomorphic to $X_K$ for some interval $K \subset \{1, \ldots, n\}$.
\end{proposition}

Note that knowing the indecomposable representations of $\rep_\mathbb{K}(Q)$ and knowing the morphisms between them is enough to describe the entire category $\rep_\mathbb{K}(Q)$  This data is contained in the so-called \new{Auslander--Reiten quiver of $Q$ (over $\mathbb{K}$)}, denoted by $\AR_\mathbb{K}(Q)$. The vertices of $\AR_\mathbb{K}(Q)$ are the isomorphism classes of the indecomposable representations of $\rep_\mathbb{K}(Q)$, and the arrows correspond to the irreducible representations between them. 

\begin{proposition} \label{prop:ARCoxandRep}
	Let $c \in \mathfrak{S}_{n+1}$ be a Coxeter element. Then the quivers $\AR(c)$ and $\AR_\mathbb{K}(Q(c))^{\op}$ are isomorphic. More precisely, the map $\Psi : \AR(c) \longrightarrow \AR_\mathbb{K}(Q(c))^{\op}$ defined by:
	\begin{enumerate}[label=$\bullet$,itemsep=1mm]
		\item $\Psi_0((i,j)) = X_{\llrr{i,j-1}}$ for all $1 \leqslant i < j \leqslant n+1$, and,
		
		\item $\Psi_1(((i,j),(k,\ell))) = (X_{\llrr{i,j-1}},X_{\llrr{k,\ell-1}})$ for all arrow $((i,j),(k,\ell))$ in $\AR(c)$,
	\end{enumerate}
	is a quiver (directed graph) isomorphism.
\end{proposition}

\begin{remark} \label{rem:ARCoxandRep}
	To see further details about Auslander-Reiten quivers of $A_n$ type quivers, we refer the reader to \cite[Section 3.1]{Sch14}. To learn more about quiver representation theory and for more in-depth knowledge on the notion of Auslander--Reiten quivers, we invite the reader to look at \cite{ASS06}.
\end{remark}

In the following, we focus on full subcategories of $\rep_\mathbb{K}(Q)$ that are closed under sums and summands. Thus, those categories are additively generated by indecomposable representations, up to isomorphism, and are characterized by the sets of indecomposable representations. 

From now one, we write $\mathcal{I}_n$ the set of intervals in $\{1,\ldots,n\}$. Given such a category $\mathscr{C}$, we write $\Int(\mathscr{C})$ for the interval set $\mathscr{J} \subseteq \mathcal{I}_n$ such that $\mathscr{C}$ is additively generated by the indecomposable representations $X_K$ for $K \in \mathscr{J}$. Given an interval set $\mathscr{J} \subseteq \mathcal{I}_n$, we denote by $\Cat_Q(\mathscr{J})$ for the category additively generated by the $X_K$ for $K \in \mathscr{J}$.

\subsection{Reflection functors}
\label{ss:refl}

In this subsection, we recall the definition of \emph{reflection functors} for any quiver $Q$. For our purposes in this paper, defining these functors only on objects is sufficient. We refer the reader to \cite{BGP73} and \cite[Section VII.5]{ASS06} for more details.

Let $Q$ be an arbitrary quiver and $v$ be a vertex $Q$  Denote $\sigma_v(Q)$ the quiver obtained from $Q$ by reversing the directions of the arrows incident to $v$  If $\alpha \in Q_1$ such that $v \in \{s(\alpha),t(\alpha)\}$, denote $\tilde{\alpha}$ the reversed arrow of $\alpha$ in $\sigma_v(Q)$. 

Now assume that $v$ is a sink of $Q$  Consider $\Xi = \sigma_v(Q)$. The \new{reflection functor} $$\mathcal{R}_v^+ : \rep(Q) \longrightarrow \rep(\Xi)$$ is defined as follows. Let $E = ((E_q)_{q \in Q_0}, (E_\beta)_{\beta \in Q_1}) \in \rep(Q)$. We set $\mathcal{R}_v^+(E) = ((F_q)_{q \in \Xi_0}, (F_\beta)_{\beta \in \Xi_1}) \in \rep(\Xi)$ where \begin{enumerate}[label = $\bullet$]
	\item $F_q = E_q$ for $q \neq v$ and $$\displaystyle F_v = \Ker \left(\bigoplus_{\alpha \in Q_1,\ t(\alpha) = v} E_\alpha : \bigoplus _{\alpha \in Q_1,\ t(\alpha) = v} E_{s(\alpha)} \longrightarrow E_v \right);$$
	\item $F_\beta= E_\beta$ if $\beta \in Q_1$ such that $t(\beta) \neq v$, otherwise $F_{\tilde{\beta}} : Y_v \longrightarrow E_{s(\beta)}$ is the composition of the kernel inclusion of $F_v$ to $\displaystyle \bigoplus_{\alpha \in Q_1,\ t(\alpha) = v} E_{s(\alpha)}$ with the projection onto the direct summand $E_{s(\beta)}$.
\end{enumerate}
If $v$ is a source of $Q$, the \new{reflection functor} $$\mathcal{R}_v^- : \rep(Q) \longrightarrow \rep(\sigma_v(Q))$$ is defined dually. We refer the reader to \cite{Deq23} for an explicit example of the calculation of $\mathcal{R}^\pm_v(E)$ given an $A$ type quiver $Q$, and a representation $E \in \rep_\mathbb{K}(Q)$.

The reflection functors are additive, so we can understand their actions on objects by knowing their actions on indecomposable objects.

By the following propositions, we recall the interpretation of reflection functors for Coxeter elements of $\mathfrak{S}_{n+1}$, and their action on $\ind_\mathbb{K}(Q)$, for $Q$ an $A_n$ type quiver.

\begin{proposition}\label{prop:reflandconj} Let $c \in \mathfrak{S}_{n+1}$ be a Coxeter element. Then if $v$ is either a sink or a source of $Q(c)$, then $\sigma_v(Q(c)) = Q(s_v c s_v)$.
\end{proposition}

\begin{proposition}\label{prop:reflonAntypequivers} Let $Q$ be an $A_n$ type quiver, $v \in Q_0$ and $\llrr{v} \neq K \in \mathcal{I}_n$. Write $\Xi = \sigma_v(Q)$. If $v$ is a sink of $Q$, then $\mathcal{R}_v^+(X_K) \cong X_{K'} \in \rep(\Xi)$ where $$K' = \begin{cases}
		K \cup \{v\} & \text{if either } e(K) = v-1 \text{ or } b(k) = v+1; \\
		K \setminus \{v\} & \text{if either } e(K)=v \text{ or } b(K) = v;\\
		K & \text{otherwise.}
	\end{cases}.$$ If $v$ is a source of $Q$, then $\mathcal{R}_v^-(X_K) = X_{K'}$ where $K'$ is defined as above.
\end{proposition}
Note that, if $v$ is a sink of $Q$, $\mathcal{R}_v^+(X_{\llrr{v}}) = 0$, and if $v$ is a source, $\mathcal{R}_v^-(X_{\llrr{v}}) = 0$. We also recall the following result, which highlights one of the main algebraic uses of those functors.
\begin{theorem} \label{thm:eqofcatrefl}
	Let $Q$ be a quiver and $v$ be one of its sinks. Write $\Xi = \sigma_v(Q)$. The reflection functor $\mathcal{R}_v^+: \rep(Q) \longrightarrow \rep(\Xi)$ induces a category equivalence between the full subcategory of $\rep(Q)$ additively generated by the indecomposable representations of $Q$ except the simple projective representation at $v$ and the full subcategory of $\rep(\Xi)$ additively generated by indecomposable representations of $\Xi$ except the simple injective representation at $v$  The quasi-inverse is induced by the reflection functor $\mathcal{R}_v^-:\rep(\Xi) \longrightarrow \rep(Q)$.
\end{theorem}
For more details, see \cite[Theorem VII.5.3]{ASS06}.

\subsection{Jordan recoverability and canonical Jordan recoverability}
\label{ss:JRandCJR}

Consider $E \in \rep_\mathbb{K}(Q)$. An endomorphism $N: E \longrightarrow E$ is \new{nilpotent} if $N^k = 0$ for some integer $k \geqslant 1$. One can see a nilpotent endomorphism as a collection of nilpotent transformations $(N_q)_{q \in Q_0}$ with some additional compatibility relations. Write $\NEnd(E)$ for the set of nilpotent endomorphisms. 

Assume that $\vdim(E)= \pmb{d} = (d_q)_{q \in Q_0}$. Given $N \in \NEnd(E)$, we consider the Jordan form of $N_q$ at each $q \in Q_0$  It gives us a sequence of integer partitions $\lambda^q \vdash d_q$  We refer to $\pmb{\lambda} = (\lambda^q)_{q \in Q_0}$ as the \new{Jordan form of $N$}, and we set $\JF(N) = \pmb{\lambda}$. On integer partitions, we consider the \new{dominance order}, denote by $\unlhd$, defined as follows: for $\lambda$ and $\mu$ two integer partitions of $d \in \mathbb{N}$, $\lambda \unlhd \mu$ whenever, for all $k \geqslant 1$, $\lambda_1 + \ldots + \lambda_k \leqslant \mu_1+ \ldots + \mu_k$. We extend this order to $n$-tuples of partitions  First, let us introduce a notation  Given $\pmb{d} = (d_k)_{1 \leqslant k \leqslant n} \in \mathbb{N}^n$ and $\pmb{\lambda} = (\lambda^k)_{1 \leqslant k \leqslant n}$ a $n$-tuple of integer, we write $\pmb{\lambda} \pmb{\vdash} \pmb{d}$ whenever $\lambda^k \vdash d_k$ for all $k \in \{1,\ldots,n\}$. Now fix $\pmb{d} \in \mathbb{N}^n$, and $\pmb{\lambda},\pmb{\mu} \pmb{\vdash} \pmb{d}$. We write $\pmb{\lambda} \pmb{\unlhd} \pmb{\mu}$  whenever for all $k \in \{1,\ldots, n\}$, $\lambda^k \unlhd \mu^k$.

\begin{theorem}[\cite{GPT19},\cite{Deq23}]
	\label{thm:defGenJF}
	Let $Q$ be an $A_n$ type quiver, and $\mathbb{K}$ be an algebraically closed field  Consider $E \in \rep_\mathbb{K}(Q)$. 
	\begin{enumerate}[label = $(\roman*)$]
		\item  The set $\NEnd(E)$ is an irreducible algebraic variety.
		\item  There exists a maximal value of $\JF$ in $\NEnd(E)$, with respect to $\pmb{\unlhd}$, and it is attained in a dense open set (in the Zariski topology) of $\NEnd(E)$.
	\end{enumerate}
\end{theorem}

Note that the previous result can be generalized to finitely generated modules over an arbitrary $\mathbb{K}$-algebra (see \cite{GPT19}).

\begin{definition}
	Let $Q$ be an $A_n$ type quiver, and $\mathbb{K}$ be an algebraically closed field  For all $E \in \rep_\mathbb{K}(Q)$, we call the \new{generic Jordan form data of $E$}, denoted by $\GenJF(E)$, the maximal value of $\JF$ in $\NEnd(E)$
\end{definition}

This definition is entirely algebraic. We introduce a combinatorial method using Greene--Kleitman invariant calculations to compute $\GenJF(E)$ effectively.

Fix an $A_n$ type quiver $Q$  Let $E \in \rep_\mathbb{K}(Q)$. We decompose $E$ as below
\[E \cong \bigoplus_{K \in \mathcal{I}_n} X_K^{m_K},\]
with $m_K \in \mathbb{N}$. Using this decomposition above allows us to see a representation $E$, up to isomorphism, as a filling $\opf_Q(E)$ of $\AR_\mathbb{K}(Q)$ as follows: \[\forall K \in \mathcal{I}_n,\ \opf_Q(E)(K) = m_K.\] In the following, for any filling $f$ of $\AR_\mathbb{K}(Q)$, we denote by $\opE_Q(f)$ its associated representation in $\rep_\mathbb{K}(Q)$, defined up to isomorphism.

For each $k \in \{1,\ldots,n\}$, write $\AR_\mathbb{K}^{[k]}(Q)$ for the complete subquiver of $\AR_\mathbb{K}(Q)$ whose vertices are given by $X_K$ such that $k \in K$. It induces a filling $\opf_Q(E)^{[k]}$ of $\AR_\mathbb{K}(Q)$. Using \cref{prop:AcycAR} and \cref{prop:ARCoxandRep}, the direct graph $\AR_\mathbb{K}(Q)$ is acyclic, which allows us to calculate the Greene--Kleitman invariant of $\opf_Q(E)^{[k]}$. 

\begin{theorem}[{\cite[Theorem 3.6]{Deq23}}] \label{thm:GenJFGK}
	Let $Q$ be an $A_n$ type quiver  For any $X \in \rep_\mathbb{K}(Q)$, we have \[\GenJF(E) = \left(\GK_{\AR_\mathbb{K}^{(k)}(Q)}(\opf_Q(E)^{[k]}) \right)_{k \in Q_0}.\] 
\end{theorem}

\begin{remark} \label{rem:acyclicityARviareptheory}
	For any $A_n$ type quiver $Q$, the acyclicity of $\AR_\mathbb{K}(Q)$, needed to calculate the Greene--Kleitman invariants, can be justified using \cref{prop:ARCoxandRep} and \cref{prop:AcycAR}. There exists a representation-theoretic explanation: $Q$ is a \emph{representation-directed quiver}. We refer the reader to \cite[Section IX.3.]{ASS06} for more details.
\end{remark}

The generic Jordan form data is an invariant on $\rep_\mathbb{K}(Q)$, but not complete. Nevertheless, we remain interested in determining the subcategories of $\rep_\mathbb{K}(Q)$ such that $\GenJF$ is complete, which brings us to the following definition.

\begin{definition} \label{def:JRcats}
	Let $Q$ be an $A_n$ quiver  A subcategory $\mathscr{C} \subset \rep_\mathbb{K}(Q)$  is \new{Jordan recoverable} if, for any tuple of integer partitions $\pmb{\lambda}$, there exists at most one $E \in \mathscr{C}$ such that $\GenJF(E) \cong \pmb{\lambda}$.
\end{definition}

\cite{Deq23} and \cite{GPT19}.

It is still difficult to characterize all the Jordan recoverable subcategories of $\rep_\mathbb{K}(Q)$.  In \cref{sec:Further}, we recall and restate a conjecture originally stated in \cite{Deq23}. However, we can focus on subcategories in which we can make explicit an algebraic inverse to $\GenJF$.

Fix $Q$ an $A_n$ type quiver  Consider $\pmb{\lambda}$ to be an $n$-tuple of integer partitions  We can consider the set $\rep_\mathbb{K}(Q,\pmb{\lambda})$ of representations $F \in \rep_\mathbb{K}(Q)$ such that it admits a nilpotent endomorphism $N$ with Jordan form $\JF(N) = \pmb{\lambda}$. We can ask for the existence of a dense open set (for Zariski's topology) $\Omega \subset \rep_\mathbb{K}(Q,\pmb{\lambda})$ such that all the representations in $\Omega$ are isomorphic to each other. Via \cite[Corollary 2.5]{GPT19}, such a set $\Omega$ exists. 

\begin{definition} \label{def:GenRep}
	Given $Q$ and $A_n$ type quiver, and any $n$-tuple of integer partitions $\pmb{\lambda}$, we define the \new{generic representation of $\pmb{\lambda}$}, denoted by $\GenRep(\pmb{\lambda})$, to be the representation such that there exists a dense open set $\Omega$ in $\rep_\mathbb{K}(Q,\pmb{\lambda})$ such that, for all $F \in \Omega$, $F \cong \GenRep(\pmb{\lambda})$.
\end{definition}

Note that, in general, for any $E \in \rep_\mathbb{K}(Q)$, $\GenRep(\GenJF(E)) \ncong E$, and for any $n$-tuple of partitions $\pmb{\lambda}$, $\GenJF(\GenRep(\pmb{\lambda})) \neq \pmb{\lambda}$.

It makes sense that whenever $\pmb{\lambda}$ does not correspond to a generic Jordan form data, applying $\GenJF \circ \GenRep$ does not give us back $\pmb{\lambda}$  Similarly, whenever we have two representations $E, F \in \rep_\mathbb{K}(Q)$ such that $\pmb{\lambda} = \GenJF(E) = \GenJF(F)$ but $E \ncong F$, then $\GenRep(\pmb{\lambda}) \ncong E$ or $\GenRep(\pmb{\lambda}) \ncong F$. Thus, we can focus on categories in which, at least, we can recover any representation from its generic Jordan form data by applying $\GenRep$.

\begin{definition} \label{def:CJRcats}
	A subcategory $\mathscr{C} \subset \rep_\mathbb{K}(Q)$ is said to be \new{canonically Jordan recoverable} whenever for all $E \in \mathscr{C}$, $\GenRep(\GenJF(E)) \cong E$.
\end{definition}

Note that any canonically Jordan recoverable category is Jordan recoverable, but the converse is false. We refer the reader to \cite{Deq23,GPT19} for examples.

\subsection{Operations preserving canonical Jordan recoverability}
\label{ss:CJRpreserving}

In this section, we restate some results about operations that preserve the canonical Jordan recoverability property from \cite{Deq23}.

First, we recall results involving the reflection functors. 

\begin{proposition}[{\cite[Section 4.3]{Deq23}}] \label{prop:GenfJFRepRefl}
	Let $v$ be a vertex of a quiver of $A_n$ type. Let $\pmb{\pi}$ be a $n$-tuple of
	integer partitions such that $\pmb{\pi}$ is either $(\boxplus, \boxplus)$-storable, $(\boxplus, \boxminus)$-storable, $(\boxminus, \boxplus)$
	storable or strongly $(\boxminus, \boxminus)$-storable at $v$. Consider $E = \GenRep(\pmb{\pi})$ and assume
	that $\GenJF(E) = \pmb{\pi}$. The following assertions hold: 
	\begin{enumerate}[label=$(\roman*)$,itemsep=1mm]
		\item If $v$ is a source, $\GenJF(\mathcal{R}^-_v (E)) = \cdiag_v(\pmb{\pi})$, and $\mathcal{R}^-_v(E) \cong \GenRep(\cdiag_v(\pmb{\pi}))$.
		
		\item If $v$ is a sink, $\GenJF(\mathcal{R}^+_v (E)) = \cdiag_v(\pmb{\pi})$, and $\mathcal{R}^+_v(E) \cong \GenRep(\cdiag_v(\pmb{\pi}))$.
	\end{enumerate}
\end{proposition}

\begin{theorem}[{\cite[Section 5.2]{Deq23}}] \label{thm:CJRefl} Let $v$ be a source of an $A_n$-type quiver $Q$. Let $\mathscr{C} \subset \rep_\mathbb{K}(Q)$
	be a canonically Jordan recoverable category such that \begin{enumerate}[label = $(\nabla)$]
		\item \label{nabla} For all $E \in \mathscr{C}$, $\GenJF(E)$ is either $(\boxplus, \boxplus)$-storable, $(\boxminus,\boxplus)$-storable, $(\boxplus,\boxminus)$-storable or strongly $(\boxminus,\boxminus)$-storable at $v$.
	\end{enumerate} Then $R^{-}_v (\mathscr{C})$ is a
	canonically Jordan recoverable subcategory of $\rep_\mathbb{K}(\sigma_v(Q))$. Similarly, if $v$ is a sink,  
	$R^{+}_v (\mathscr{C})$ is a canonically Jordan recoverable subcategory of 
	$\rep_\mathbb{K}(\sigma_v(Q))$.
\end{theorem}

Then, we recall some results about the \emph{adding simple operations}. Let $Q$ be an $A_n$ type quiver. For any $v \in Q_0$, we define the \new{adding simple at $v$ operation}, denoted by $\Adds_v$, on every subcategory $\mathscr{C}$ of $\rep_\mathbb{K}(Q)$  by $\Adds_v(\mathscr{C}) = \add(\mathscr{C}, X_{\llrr{v}})$.

\begin{lemma}[{\cite[Section 5.1]{Deq23}}] \label{lem:GenJFAddS}
	Let $v$ be a source or a sink of an $A_n$ type quiver $Q$. Consider $a \in \mathbb{N}$
	and $E \in rep_\mathbb{K}(Q)$. Write $\pmb{\pi} = \GenJF(E)$. Then $\GenJF(X_{\llrr{v}}^a
	\oplus E) = \pmb{\xi}$ where $\xi^q = \pi^q$
	if $q \neq v$ and $\xi^v = (\pi^v_1 + a, \pi^v_2, \pi^v_3,\ldots)$.
\end{lemma}

\begin{theorem}[{\cite[Section 5.1]{Deq23}}] \label{thm:CJRAddS} 
	Let $Q$ be a quiver of $A_n$ type, and $v$ be a source or a sink of $Q$.
	Let $\mathscr{C} \subset \rep_\mathbb{K}(Q)$ be a canonically Jordan recoverable category such that
	\begin{enumerate}[label = $(\ast)$]
		\item \label{ast} For any $E \in \mathscr{C}$, $\GenJF(E)$ is strongly $(\boxminus,\boxminus)$-storable at $v$.
	\end{enumerate} Then
	$ \Adds_v(\mathscr{C})$ is a canonically Jordan recoverable subcategory of $\rep_\mathbb{K}(Q)$.
\end{theorem} 

\subsection{Results from previous work} \label{ss:previously}

In this section, we resume some results from \cite{Deq23}, and by going through the main outlines of the proof, we establish some valuable consequences in the next section.

First, we state the main result of \cite{Deq23}. Let $Q$ be an $A_n$ type quiver. We recall that the subcategories we consider are characterized by their indecomposable representations and, therefore, by intervals of $\{1, \ldots, n\}$. 

For $K = \llbracket i;j \rrbracket \in \mathcal{I}_n$, we set $b(K) = i$ and $e(K)=j$. Two intervals $K,L \in \mathcal{I}_n$ are \new{adjacent} whenever either $b(L) = e(K)+1$ or $b(K) = e(L)+1$. Hence, $K$ and $L$ are not adjacent if either $K \cap L \neq \varnothing$, or $b(K) > e(L)+1$, or $b(L) > e(K)+1$.

\begin{definition} \label{def:AdjAvoidInt}
	An interval set $\mathscr{J} \subset \mathcal{I}_n$ is said to be \new{adjacency-avoiding} whenever there is no pair of intervals $K,L \in \mathscr{J}$ that are adjacent.
\end{definition}

\begin{example}\label{ex:MaxAdjAvoidInt}
	Below are given all the maximal adjacency-avoiding interval sets of $\mathcal{I}_3$:
	\begin{enumerate}[label = $\bullet$, itemsep=1mm]
		\item $\{\llrr{1}, \llrr{1,2}, \llrr{1,3}\}$;
		\item $\{\llrr{2}, \llrr{1,2}, \llrr{2,3},\llrr{1,3}\}$;
		\item $\{\llrr{3},\llrr{2,3},\llrr{1,3}\}$;
		\item $\{\llrr{1}, \llrr{3}, \llrr{1,3}\}$. \qedhere
	\end{enumerate}
\end{example}

Below is the main result of \cite{Deq23}, which gives a combinatorial characterization of all the subcategories of $\rep_\mathbb{K}(Q)$.

\begin{theorem} \label{thm:CJR}
	Let $Q$ be an $A_n$ type quiver  A subcategory $\mathscr{C} \subset \rep_\mathbb{K}(Q)$ is canonically Jordan recoverable if, and only if, $\Int(\mathscr{C})$ is adjacency-avoiding. 
\end{theorem}

\begin{remark} $ $
	\begin{enumerate}[label = $\bullet$, itemsep=1mm]
		\item This result generalizes a previous result of \cite{GPT19} for $A_n$ type quivers.
		
		\item For any interval set $\mathscr{J}$, the canonical Jordan recoverability of  $\Cat_Q(\mathscr{J})$ does not depend of $Q$.
	\end{enumerate}
\end{remark}

The first step of the proof was to highlight that the adjacency-avoiding property of $\Int(\mathscr{C})$ is a necessary condition for some category $\mathscr{C}$ to be canonically Jordan recoverable.

Then, we study the maximal adjacency-avoiding interval sets of $\{1,\ldots,n\}$.  Before stating the following result, let us introduce a relevant family of interval sets. For any pair of subsets $(\B,\E)$ of $\{1,\ldots, n\}$, we define the interval set $\mathscr{J}(\B,\E)$ as follows: \[\mathscr{J}(\B,\E) = \{ \llrr{b,e} \mid b \in \B,\ e \in \E\}.\] 

\begin{theorem}[\cite{Deq23}]
	Any adjacency-avoiding interval set is a subset of some $\mathscr{J}(\B,\E)$ such that $(\B,\E[1])$ is an elementary interval bipartition of $\{1,\ldots,n+1\}$. 
\end{theorem}

By \cref{cor:Coxbipart}, we can parametrize maximal adjacency-avoiding interval subsets of $\mathcal{I}_n$ via a Coxeter element $c \in \mathfrak{S}_{n+1}$. Write $\mathscr{J}(c) = \mathscr{J}(\L_c \cup \{1\}, \R_c[-1]\cup\{n\})$. 

We set $\Cat_Q(c) = \Cat_Q(\mathscr{J}(c))$. Finally, we show that, given an $A_n$ quiver, all the subcategories $\Cat_Q(c) \subset \rep_{\mathbb{K}}(Q)$ arising from maximal adjacency-avoiding interval subsets of $\mathcal{I}_n$ are canonically Jordan recoverable.  

\begin{theorem}[{\cite[Theorem 7.5]{Deq23}}] \label{thm:CJRprecise} Fix $n \in \mathbb{N}^*$. Let $\mathbb{K}$ be an algebraically closed field, and $Q$ be an $A_n$ type quiver. Consider $c \in \mathfrak{S}_{n+1}$ a Coxeter element. Then $\Cat_Q(c)$ is canonically Jordan recoverable. Moreover, $\GenJF$ induces a one-to-one correspondence from isomorphism classes of representations of $\Cat_Q(c)$ to $\st(c)$.
\end{theorem}

\begin{remark} \label{rem:morethaninj}
	In \cite[Theorem 7.5]{Deq23}, we state that $\GenJF$ induces an injective map from isomorphism classes of representations of $\Cat_Q(c)$ to $\st(c)$. The steps used in the proofs of \cite[Theorem 7.4 and 7.5]{Deq23} can be reversed. Thus $\GenJF$ indeed induces a one-to-one correspondence.
\end{remark}

The proof uses technical tools from combinatorics and quiver representation theory, particularly results presented in \cref{ss:CJRpreserving}. The proof of \cite[Theorem 7.5]{Deq23} uses induction on the orientations of the quiver $Q$, which can be seen in terms of Coxeter elements. We prove it for $Q(c_1)$ for $c_1 = (1,\ldots,n+1)$, and then we prove that if the statement holds for $Q(c_2)$ given by some Coxeter element $c_2 \in \mathfrak{S}_{n+1}$, then it holds for $Q(sc_2s)$ for any $s \in \Sigma_n$ that is either initial or final in $c_2$. 

Let $c \in \mathfrak{S}_{n+1}$ be a Coxeter element, and $k \in \{1,\ldots,n\}$. We define $\sigma_k(c)$ by:
\begin{enumerate}[label=$\bullet$,itemsep=1mm]
	\item $\sigma_k(c) = s_kcs_k$ whenever $s_k$ is either initial or final in $c$, and
	\item $\sigma_k(c) = c$ otherwise.
\end{enumerate}
Note that we use the same notation for mutations on quivers, and we saw that in our setting, they match (see Proposition \ref{prop:reflandconj}). We also consider a tiny deformation of this mutation, denoted by $\widetilde{\sigma_k}$, defined by:
\begin{enumerate}[label=$\bullet$,itemsep=1mm]
	\item $\widetilde{\sigma_k}(c) = c$ whenever $k=1$ and $s_1$ initial in $c$, or $k=n$ and $s_n$ final in $c$;
	\item $\widetilde{\sigma_k}(c) = \sigma_k(c)$ otherwise.
\end{enumerate}

\begin{proposition} \label{prop:reflonCJR} Let $c_1,c_2 \in \mathfrak{S}_{n+1}$ be two Coxeter elements. Fix $k \in \{1,\ldots,n\}$ such that $s_k$ is either initial or final in $c_1$. Then $\mathcal{R}_k^{\pm}(\Cat_{Q(c_1)}(c_2))$ is equivalent to a subcategory of $\Cat_{Q(c_1)}(\widetilde{\sigma_k}(c_2))$, as subcategory of $\rep_\mathbb{K}(Q(\sigma_k(c_1)))$.
\end{proposition}

\begin{proof}
	This is a consequence of \cref{prop:reflandconj} for $c_1$, and \cref{prop:reflonAntypequivers} for $c_2$.
\end{proof}

In the following proposition, we rephrase the precise result hidden in the proof of \cite[Theorem 7.5]{Deq23}.

\begin{theorem}[\cite{Deq23}] \label{thm:togRep}
	Let $\mathbb{K}$ be an algebraically closed field. Fix $n \in \mathbb{N}^*$. Let $c_1,c_2 \in \mathfrak{S}_{n+1}$ be two Coxeter elements.   Let $\Cat_{Q(c_1)}(c_2) \subset \rep_{\mathbb{K}}(Q(c_1))$. Fix $k \in \{1,\ldots,n\}$ such that $s_k \in \Sigma_n$ is either initial or final in $c_1$. For all $E \in \Cat_{Q(c_1)}(c_2)$,
	\[\GenJF(\mathcal{R}^{\pm}_k(E)) = \cdiag_k(\GenJF(E)).\] Moreover, $\GenJF(\mathcal{R}^{\pm}_k(E)) \in \st(\widetilde{\sigma_{k}}(c_2))$, and if $s_k$ is final in $c_2$, then $\GenJF(\mathcal{R}^{\pm}_k(E))$ is strongly $(\boxminus,\boxminus)$-storable at $k$.
\end{theorem}

	\section{Extended RSK via type $A$ quiver representation}
	\label{sec:Extension}
	\subsection{Construction of the extended generalization}

Next, we describe an extended version of RSK using (type $A$) Coxeter elements and state the main result.

Let $n \geqslant 1$ and $\lambda \in \Hk_n$.  Let $c \in \mathfrak{S}_{n+1}$, and consider $\AR(c)$ its Auslander--Reiten quiver. Recall that $\opc(\lambda) \in \mathfrak{S}_{n+1}$ is the unique Coxeter element such that $\oplamb(\opc(\lambda)) = \lambda$ by \cref{cor:Coxbipart}. We also recall that $(\L = \L_c \cup \{1\}, \R = \R_c \cup \{n+1\})$ is the unique elementary bipartition of $\{1,\ldots,n+1\}$ such that $\oplamb(\L,\R) = \lambda$ by \cref{prop:elemintbipartintpartitions}.

Using $\opc(\lambda)$-coordinates on $\Fer(\lambda)$, any box labelled $\LR{\ell,r}_{\opc(\lambda)}$ can be associated to the transposition $(\ell,r) \in \mathfrak{S}_{n+1}$. We construct a one-to-one correspondence \new{$\crep_{\lambda, c}$} from fillings of shape $\lambda$ to those of the Auslander--Reiten quiver $\AR(c)$ which are supported on vertices $(\ell,r) \in \L \times \R$ such that $\ell < r$. Precisely, for any filling $f$ of shape $\lambda$, $\crep_{\lambda,c} (f)(\ell,r) = f(\LR{\ell,r}_{\opc(\lambda)})$ whenever $(\ell,r) \in \L \times \R$ such that $\ell < r$,  and $\crep_{\lambda,c}(\ell,r) = 0$ otherwise.

As in \cref{ss:Gansner}, we also use the $\lambda$-diagonal coordinates. 
For each $k \in \{1,\ldots,n\}$, we consider the subgraph $\AR^{[k]}(c)$ of $\AR(c)$ where the vertices are the transpositions $(\ell,r)$ with $\ell \leqslant k < r$. This subgraph has only one source and only one sink. We denote by $\crep_{\lambda,c}(f)^{[k]}$ the filling or $\AR^{[k]}(c)$ induced by $\crep_{\lambda,c}(f)$.

We define \new{$\GRSK_{\lambda, c}(f)$} to be the fillings of shape $\lambda$ defined, for $k \in \{1,\ldots,n\}$, and for $\delta \in \{1,\ldots,\delta_k\}$, by:
\[\GRSK_{\lambda,c}(f)(\ldiag{k,\delta}_{\lambda}) = \GK_{\AR^{[k]}(c)} \left(\crep_{\lambda,c}(f)^{[k]} \right)_{\delta}.\]

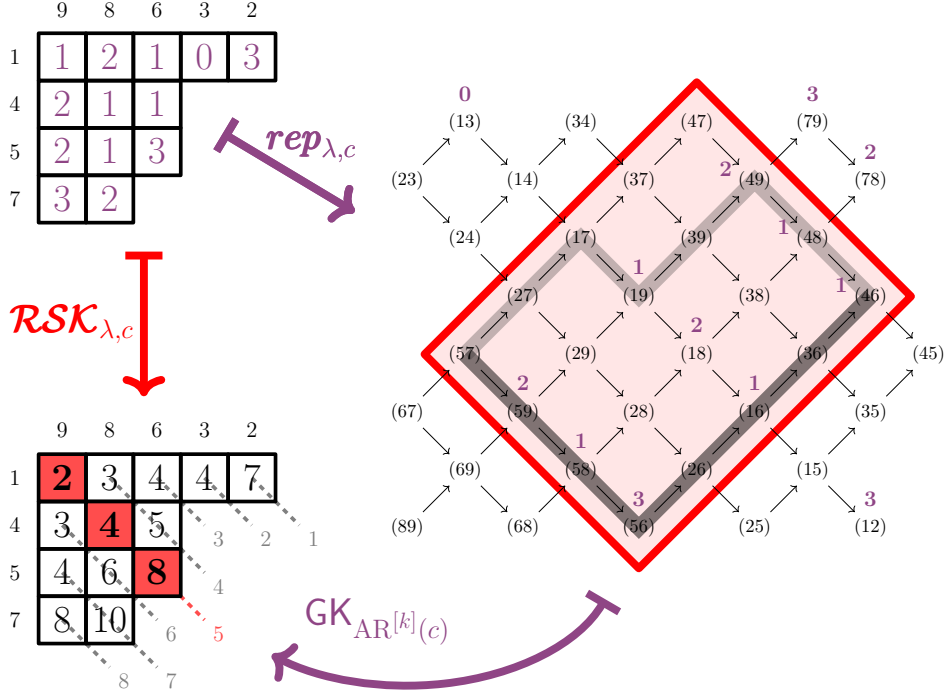
\begin{figure}[!ht]
	\centering
	\[
	\scalebox{0.7}{
		\begin{tikzpicture}
			
			\begin{scope}[xshift=-7cm, yshift = 8.5cm,scale=0.9]
				\tkzDefPoint(0,0){a}
				\tkzDefPoint(0,1){b}
				\tkzDefPoint(1,1){c}
				\tkzDefPoint(1,0){d}
				\tkzDrawPolygon[line width = 0.7mm, color = black](a,b,c,d);
				
				\tkzDefPoint(1,0){a}
				\tkzDefPoint(1,1){b}
				\tkzDefPoint(2,1){c}
				\tkzDefPoint(2,0){d}
				\tkzDrawPolygon[line width = 0.7mm, color = black](a,b,c,d);
				
				\tkzDefPoint(2,0){a}
				\tkzDefPoint(2,1){b}
				\tkzDefPoint(3,1){c}
				\tkzDefPoint(3,0){d}
				\tkzDrawPolygon[line width = 0.7mm, color = black](a,b,c,d);
				
				\tkzDefPoint(3,0){a}
				\tkzDefPoint(3,1){b}
				\tkzDefPoint(4,1){c}
				\tkzDefPoint(4,0){d}
				\tkzDrawPolygon[line width = 0.7mm, color = black](a,b,c,d);
				
				\tkzDefPoint(4,1){a}
				\tkzDefPoint(4,0){b}
				\tkzDefPoint(5,0){c}
				\tkzDefPoint(5,1){d}
				\tkzDrawPolygon[line width = 0.7mm, color = black](a,b,c,d);
				
				\tkzDefPoint(0,0){a}
				\tkzDefPoint(0,-1){b}
				\tkzDefPoint(1,-1){c}
				\tkzDefPoint(1,0){d}
				\tkzDrawPolygon[line width = 0.7mm, color = black](a,b,c,d);
				
				\tkzDefPoint(1,0){a}
				\tkzDefPoint(1,-1){b}
				\tkzDefPoint(2,-1){c}
				\tkzDefPoint(2,0){d}
				\tkzDrawPolygon[line width = 0.7mm, color = black](a,b,c,d);
				
				\tkzDefPoint(2,0){a}
				\tkzDefPoint(2,-1){b}
				\tkzDefPoint(3,-1){c}
				\tkzDefPoint(3,0){d}
				\tkzDrawPolygon[line width = 0.7mm, color = black](a,b,c,d);
				
				\tkzDefPoint(0,-2){a}
				\tkzDefPoint(0,-1){b}
				\tkzDefPoint(1,-1){c}
				\tkzDefPoint(1,-2){d}
				\tkzDrawPolygon[line width = 0.7mm, color = black](a,b,c,d);
				
				\tkzDefPoint(1,-2){a}
				\tkzDefPoint(1,-1){b}
				\tkzDefPoint(2,-1){c}
				\tkzDefPoint(2,-2){d}
				\tkzDrawPolygon[line width = 0.7mm, color = black](a,b,c,d);
				
				\tkzDefPoint(2,-2){a}
				\tkzDefPoint(2,-1){b}
				\tkzDefPoint(3,-1){c}
				\tkzDefPoint(3,-2){d}
				\tkzDrawPolygon[line width = 0.7mm, color = black](a,b,c,d);
				
				\tkzDefPoint(0,-3){a}
				\tkzDefPoint(0,-2){b}
				\tkzDefPoint(1,-2){c}
				\tkzDefPoint(1,-3){d}
				\tkzDrawPolygon[line width = 0.7mm, color = black](a,b,c,d);
				
				\tkzDefPoint(1,-3){a}
				\tkzDefPoint(1,-2){b}
				\tkzDefPoint(2,-2){c}
				\tkzDefPoint(2,-3){d}
				\tkzDrawPolygon[line width = 0.7mm, color = black](a,b,c,d);
				
				\node at (0.5,1.5){{\Large $9$}};
				\node at (1.5,1.5){{\Large $8$}};
				\node at (2.5,1.5){{\Large $6$}};
				\node at (3.5,1.5){{\Large $3$}};
				\node at (4.5,1.5){{\Large $2$}};
				\node at (-0.5,0.5){{\Large $1$}};
				\node at (-0.5,-0.5){{\Large $4$}};
				\node at (-0.5,-1.5){{\Large $5$}};
				\node at (-0.5,-2.5){{\Large $7$}};
				
				\tkzLabelPoint[Plum](0.5,0.95){{\Huge $1$}};
				\tkzLabelPoint[Plum](1.5,0.95){{\Huge $2$}};
				\tkzLabelPoint[Plum](2.5,0.95){{\Huge $1$}};
				\tkzLabelPoint[Plum](3.5,0.95){{\Huge $0$}};
				\tkzLabelPoint[Plum](4.5,0.95){{\Huge $3$}};
				\tkzLabelPoint[Plum](0.5,-0.05){{\Huge $2$}};
				\tkzLabelPoint[Plum](1.5,-0.05){{\Huge $1$}};
				\tkzLabelPoint[Plum](2.5,-0.05){{\Huge $1$}};
				\tkzLabelPoint[Plum](0.5,-1.05){{\Huge $2$}};
				\tkzLabelPoint[Plum](1.5,-1.05){{\Huge $1$}};
				\tkzLabelPoint[Plum](2.5,-1.05){{\Huge $3$}};
				\tkzLabelPoint[Plum](0.5,-2.05){{\Huge $3$}};
				\tkzLabelPoint[Plum](1.5,-2.05){{\Huge $2$}};
			\end{scope}
			
			\begin{scope}[xshift=-7cm, yshift = 0.5cm,scale=0.9]
				
				\draw [line width=0.7mm, red!70, dashed] (0.5,0.5) --  (3.5,-2.5);
				\node[red!70] at (3.8,-2.8){{\Large $5$}};
				
				\draw [line width=0.7mm, gray, dashed] (1.5,0.5) --  (3.5,-1.5);
				\node[gray] at (3.8,-1.8){{\Large $4$}};
				
				\draw [line width=0.7mm, gray, dashed] (2.5,0.5) --  (3.5,-0.5);
				\node[gray] at (3.8,-0.8){{\Large $3$}};
				
				\draw [line width=0.7mm, gray, dashed] (3.5,0.5) --  (4.5,-0.5);
				\node[gray] at (4.8,-0.8){{\Large $2$}};
				
				\draw [line width=0.7mm, gray, dashed] (4.5,0.5) --  (5.5,-0.5);
				\node[gray] at (5.8,-0.8){{\Large $1$}};
				
				\draw [line width=0.7mm, gray, dashed] (0.5,-0.5) --  (2.5,-2.5);
				\node[gray] at (2.8,-2.8){{\Large $6$}};
				
				\draw [line width=0.7mm, gray, dashed] (0.5,-1.5) --  (2.5,-3.5);
				\node[gray] at (2.8,-3.8){{\Large $7$}};
				
				\draw [line width=0.7mm, gray, dashed] (0.5,-2.5) --  (1.5,-3.5);
				\node[gray] at (1.8,-3.8){{\Large $8$}};
				
				\tkzDefPoint(0,0){a}
				\tkzDefPoint(0,1){b}
				\tkzDefPoint(1,1){c}
				\tkzDefPoint(1,0){d}
				\tkzDrawPolygon[line width = 0.7mm, color = black, fill=red!70](a,b,c,d);
				
				\tkzDefPoint(1,0){a}
				\tkzDefPoint(1,1){b}
				\tkzDefPoint(2,1){c}
				\tkzDefPoint(2,0){d}
				\tkzDrawPolygon[line width = 0.7mm, color = black](a,b,c,d);
				
				\tkzDefPoint(2,0){a}
				\tkzDefPoint(2,1){b}
				\tkzDefPoint(3,1){c}
				\tkzDefPoint(3,0){d}
				\tkzDrawPolygon[line width = 0.7mm, color = black](a,b,c,d);
				
				\tkzDefPoint(3,0){a}
				\tkzDefPoint(3,1){b}
				\tkzDefPoint(4,1){c}
				\tkzDefPoint(4,0){d}
				\tkzDrawPolygon[line width = 0.7mm, color = black](a,b,c,d);
				
				\tkzDefPoint(4,1){a}
				\tkzDefPoint(4,0){b}
				\tkzDefPoint(5,0){c}
				\tkzDefPoint(5,1){d}
				\tkzDrawPolygon[line width = 0.7mm, color = black](a,b,c,d);
				
				\tkzDefPoint(0,0){a}
				\tkzDefPoint(0,-1){b}
				\tkzDefPoint(1,-1){c}
				\tkzDefPoint(1,0){d}
				\tkzDrawPolygon[line width = 0.7mm, color = black](a,b,c,d);
				
				\tkzDefPoint(1,0){a}
				\tkzDefPoint(1,-1){b}
				\tkzDefPoint(2,-1){c}
				\tkzDefPoint(2,0){d}
				\tkzDrawPolygon[line width = 0.7mm, color = black, fill = red!70](a,b,c,d);
				
				\tkzDefPoint(2,0){a}
				\tkzDefPoint(2,-1){b}
				\tkzDefPoint(3,-1){c}
				\tkzDefPoint(3,0){d}
				\tkzDrawPolygon[line width = 0.7mm, color = black](a,b,c,d);
				
				\tkzDefPoint(0,-2){a}
				\tkzDefPoint(0,-1){b}
				\tkzDefPoint(1,-1){c}
				\tkzDefPoint(1,-2){d}
				\tkzDrawPolygon[line width = 0.7mm, color = black](a,b,c,d);
				
				\tkzDefPoint(1,-2){a}
				\tkzDefPoint(1,-1){b}
				\tkzDefPoint(2,-1){c}
				\tkzDefPoint(2,-2){d}
				\tkzDrawPolygon[line width = 0.7mm, color = black](a,b,c,d);
				
				\tkzDefPoint(2,-2){a}
				\tkzDefPoint(2,-1){b}
				\tkzDefPoint(3,-1){c}
				\tkzDefPoint(3,-2){d}
				\tkzDrawPolygon[line width = 0.7mm, color = black,fill = red!70](a,b,c,d);
				
				\tkzDefPoint(0,-3){a}
				\tkzDefPoint(0,-2){b}
				\tkzDefPoint(1,-2){c}
				\tkzDefPoint(1,-3){d}
				\tkzDrawPolygon[line width = 0.7mm, color = black](a,b,c,d);
				
				\tkzDefPoint(1,-3){a}
				\tkzDefPoint(1,-2){b}
				\tkzDefPoint(2,-2){c}
				\tkzDefPoint(2,-3){d}
				\tkzDrawPolygon[line width = 0.7mm, color = black](a,b,c,d);
				
				\node at (0.5,1.5){{\Large $9$}};
				\node at (1.5,1.5){{\Large $8$}};
				\node at (2.5,1.5){{\Large $6$}};
				\node at (3.5,1.5){{\Large $3$}};
				\node at (4.5,1.5){{\Large $2$}};
				\node at (-0.5,0.5){{\Large $1$}};
				\node at (-0.5,-0.5){{\Large $4$}};
				\node at (-0.5,-1.5){{\Large $5$}};
				\node at (-0.5,-2.5){{\Large $7$}};
				
				\tkzLabelPoint(0.5,0.95){{\Huge $\mathbf{2}$}};
				\tkzLabelPoint(1.5,0.95){{\Huge $3$}};
				\tkzLabelPoint(2.5,0.95){{\Huge $4$}};
				\tkzLabelPoint(3.5,0.95){{\Huge $4$}};
				\tkzLabelPoint(4.5,0.95){{\Huge $7$}};
				\tkzLabelPoint(0.5,-0.05){{\Huge $3$}};
				\tkzLabelPoint(1.5,-0.05){{\Huge $\mathbf{4}$}};
				\tkzLabelPoint(2.5,-0.05){{\Huge $5$}};
				\tkzLabelPoint(0.5,-1.05){{\Huge $4$}};
				\tkzLabelPoint(1.5,-1.05){{\Huge $6$}};
				\tkzLabelPoint(2.5,-1.05){{\Huge $\mathbf{8}$}};
				\tkzLabelPoint(0.5,-2.05){{\Huge $8$}};
				\tkzLabelPoint(1.5,-2.05){{\Huge $10$}};
			\end{scope}
			
			\begin{scope}[yshift=0cm,scale=1.1]
				\tkzDefPoint(0.3,3){a}
				\tkzDefPoint(4,-0.7){b}
				\tkzDefPoint(8.7,4){c}
				\tkzDefPoint(5,7.7){d}
				\tkzDrawPolygon[line width = 1.5mm, color = red, fill = red!10](a,b,c,d);
				
				\draw[line width = 3mm,  opacity=0.5, black] (1,3) -- (4,0) -- (8,4);
				
				\draw[line width = 3mm, opacity=0.3, black!80] (1,3) -- (3,5) -- (4,4) -- (6,6) -- (8,4);

				\node (67) at (0,2) {(67)};
				\node (57) at (1,3) {(57)};
				\node (27) at (2,4) {(27)};
				\node (17) at (3,5) {(17)};
				\node (37) at (4,6) {(37)};
				\node (47) at (5,7) {(47)};
				
				\node (69) at (1,1) {(69)};
				\node (59) at (2,2) {(59)};
				\node[Plum] at (2,2.5){{\Large $\mathbf 2$}};
				\node (29) at (3,3) {(29)};
				\node (19) at (4,4) {(19)};
				\node[Plum] at (4,4.5){{\Large $\mathbf 1$}};
				
				\node (39) at (5,5) {(39)};
				
				\node (49) at (6,6) {(49)};
				\node[Plum] at (5.5,6.2){{\Large $\mathbf 2$}};
				
				\node (68) at (2,0) {(68)};
				\node (58) at (3,1) {(58)};
				\node[Plum] at (3,1.5){{\Large $\mathbf 1$}};
				\node (28) at (4,2) {(28)};
				\node (18) at (5,3) {(18)};
				\node[Plum] at (5,3.5){{\Large $\mathbf 2$}};
				\node (38) at (6,4) {(38)};
				\node (48) at (7,5) {(48)};
				\node[Plum] at (6.5,5.2){{\Large $\mathbf 1$}};
				\node (78) at (8,6) {(78)};
				\node[Plum] at (8,6.5){{\Large $\mathbf 2$}};
				\node (89) at (0,0) {(89)};
				\node (79) at (7,7) {(79)};
				\node[Plum] at (7,7.5){{\Large $\mathbf 3$}};
				
				\node (56) at (4,0) {(56)};
				\node[Plum] at (4,0.5){{\Large $\mathbf 3$}};
				
				\node (25) at (6,0) {(25)};
				\node (12) at (8,0) {(12)};
				\node[Plum] at (8,0.5){{\Large $\mathbf 3$}};
				
				\node (13) at (1,7) {(13)};
				\node[Plum] at (1,7.5){{\Large $\mathbf 0$}};
				\node (34) at (3,7) {(34)};
				
				\node (26) at (5,1) {(26)};
				\node (15) at (7,1) {(15)};
				\node (23) at (0,6) {(23)};
				\node (14) at (2,6) {(14)};
				
				\node (16) at (6,2) {(16)};
				\node[Plum] at (6,2.5){{\Large $\mathbf 1$}};
				
				\node (35) at (8,2) {(35)};
				\node (24) at (1,5) {(24)};
				
				\node (36) at (7,3) {(36)};
				\node (45) at (9,3) {(45)};
				
				\node (46) at (8,4) {(46)};
				\node[Plum] at (7.5,4.2){{\Large $\mathbf 1$}};
				
				\draw[->] (13)--(14);\draw[->] (14)--(34);\draw[->] (14)--(17);\draw[->] (15)--(35);\draw[->] (15)--(12);\draw[->] (16)--(36);\draw[->] (16)--(15);\draw[->] (17)--(37);\draw[->] (17)--(19);\draw[->] (18)--(38);\draw[->] (18)--(16);\draw[->] (19)--(39);\draw[->] (19)--(18);\draw[->] (23)--(13);\draw[->] (23)--(24);\draw[->] (24)--(14);\draw[->] (24)--(27);\draw[->] (25)--(15);\draw[->] (26)--(16);\draw[->] (26)--(25);\draw[->] (27)--(17);\draw[->] (27)--(29);\draw[->] (28)--(18);\draw[->] (28)--(26);\draw[->] (29)--(19);\draw[->] (29)--(28);\draw[->] (34)--(37);\draw[->] (35)--(45);\draw[->] (36)--(46);\draw[->] (36)--(35);\draw[->] (37)--(47);\draw[->] (37)--(39);\draw[->] (38)--(48);\draw[->] (38)--(36);\draw[->] (39)--(49);\draw[->] (39)--(38);\draw[->] (46)--(45);\draw[->] (47)--(49);\draw[->] (48)--(78);\draw[->] (48)--(46);\draw[->] (49)--(79);\draw[->] (49)--(48);\draw[->] (56)--(26);\draw[->] (57)--(27);\draw[->] (57)--(59);\draw[->] (58)--(28);\draw[->] (58)--(56);\draw[->] (59)--(29);\draw[->] (59)--(58);\draw[->] (67)--(57);\draw[->] (67)--(69);\draw[->] (68)--(58);\draw[->] (69)--(59);\draw[->] (69)--(68);\draw[->] (79)--(78);\draw[->] (89)--(69);
			\end{scope}	
			\draw [|->,line width=1.5mm,red] (-5,5.25) -- node[left]{{\Huge $\GRSK_{\lambda,c} $}} (-5,2.5);	
			
			\draw [|->,line width=1.5mm,Plum] (-3.5,7.5) -- node[above right,xshift=-0.6cm]{{\Huge $\crep_{\lambda,c} $}} (-1,6);	
			
			\draw [|->,line width=1.5mm,Plum] (-3.5,7.5) -- node[above right,xshift=-0.6cm]{{\Huge $\crep_{\lambda,c} $}} (-1,6);	
			
			\tikzset{
				partial ellipse/.style args={#1:#2:#3}{
					insert path={+ (#1:#3) arc (#1:#2:#3)}
				}
			}
			
			\draw [|->,line width=1.5mm,Plum] (0,-0.5)  [partial ellipse=340:230:4cm and 2.5cm] node[above,xshift=2cm]{{\Huge $\GK_{\AR^{[k]}(c)}$}};
			
	\end{tikzpicture} }\]
	\caption{\label{fig:RSKandAR} Explicit calculation of $\GRSK_{\lambda,c}(f)$ for the boxes in $D_k(\lambda)$, with $k = 5$, from a filling of $\lambda = (5,3,3,2)$, with $c = (1,3,4,7,9,8,6,5,2)$.}
\end{figure}

See \cref{fig:RSKandAR} for an explicit example.

\begin{theorem} \label{mainthm}
	Let $n \geqslant 1$ and $\lambda \in \Hk_n$ . Let $c \in \mathfrak{S}_{n+1}$ be a Coxeter element. The map $\GRSK_{\lambda,c}$ realizes a one-to-one correspondence from fillings of shape $\lambda$ to reverse plane partitions of shape $\lambda$.
\end{theorem}

\begin{proof}
	See the Coxeter element $c = c_1 \in \mathfrak{S}_{n+1}$, as the $A_n$ type quiver $Q(c_1)$. See $\lambda$ as the choice of maximal canonically Jordan recoverable subcategory $\Cat_{Q(c_1)}(c_2)$ of $\rep_\mathbb{K}(Q(c_1))$ where $c_2 \in \mathfrak{S}_{n+1}$ such that $\oplamb(c_2) = \lambda$.
	
	Then the map $\crep_{\lambda,c}$ gives a one-to-one correspondence from fillings of $\lambda$ to representations in $\Cat_{Q(c_1)}(c_2)$, up to isomorphism. Set $E = \opE_{Q(c_1)}(\crep_{\lambda,c_1}(f))$ to be the representation of $\Cat_{Q(c_1)}(c_2)$ whose multiplicities are given by $\crep_{\lambda,c_1}(f)$. The calculations of $\GRSK_{\lambda,c_1}(f)$ correspond exactly to the calculations of $\GenJF(E)$ by \cref{thm:GenJFGK} and \ref{thm:bijLRstorRPP}. As stated by \cref{thm:CJRprecise}, we are done.
\end{proof}

\begin{remark} \label{rem:anybipart}
	One can notice that we could apply slightly broadened generalization by allowing to take any bipartition $(\L,\R)$ such that $\min(\L \cup \R) \in \L$ and $\max(\L\cup\R) \in \R$, and an Coxeter element $c \in \mathfrak{S}_{n+1}$ with $n \geqslant \max(R)$. In fact, this slightly deformed correspondence corresponds to the application of $\GRSK_{\lambda,c}$, for $\lambda \in \Hk_n$ such that $\L \subseteq \L_{\opc(\lambda)} \cup \{1\}$ and $\R \subseteq \R_{\opc(\lambda)} \cup \{n+1\}$, with its domain restricted to fillings that are vanishing on boxes $\LR{\ell,r}_{\opc(\lambda)}$ where $\ell \notin \L$ or $r \notin \R$.
\end{remark}

\subsection{Link to previous well-known combinatorial bijections}
\label{ss:link}

Let us start this section with an obvious proposition from the construction of $\GRSK_{\lambda, c}$.

\begin{proposition} \label{prop:conjandinvonGRSK}
	Let $n \geqslant 1$. Consider a Coxeter element $c \in \mathfrak{S}_{n+1}$, and $\lambda \in \Hk_n$. Then $\GRSK_{\lambda, c^{-1}} = \GRSK_{\lambda,c}$ and  $\GRSK_{\lambda',c} = (\GRSK_{\lambda,c})^{\top}$ where $f^{\top}$ is the filling of $\lambda'$ obtained by transposing the filling $f$ of $\lambda$.
\end{proposition}

The following result shows that we extended Gansner's RSK correspondence.

\begin{proposition} \label{prop:GRSKandnewGRSK}
	Let $n \geqslant 1$ and $c \in \mathfrak{S}_{n+1}$ be a Coxeter element. Set $\lambda = \oplamb(c)$. Then $\GRSK_{\lambda,c} = \GRSK_\lambda$.
\end{proposition}

\begin{proof} In this configuration, by setting $\L = \{1\} \cup \L_c = \{\ell_1 < \ldots \ell_p \}$ and $\R = \R_c \cup \{n+1\} = \{r_1 < \ldots < r_q\}$, we can check that: 
	\begin{enumerate}[label = $\bullet$, itemsep=1mm]
		\item for $i \in \{1,\ldots, p-1\}$, $c(\ell_i) = \ell_{i+1}$;
		\item for $j \in \{1, \ldots, q-1\}$, $c(r_{j+1}) = r_j$.
	\end{enumerate}
	Thus, the orientation of the Auslander--Reiten quiver $\AR(c)$ corresponds precisely to the reverse of the one used to calculate the $\GRSK_{\oplamb(c)}$. It induces the desired result.
\end{proof}

\begin{proposition} \label{prop:HGandnewGRSK}
	Let $n \geqslant 1$ and $\lambda \in \Hk_n$. Consider $c = (1,2,\ldots,n+1)$. Then $\GRSK_{\lambda,c}$ corresponds exactly to the Hillman-Grassl correspondence.
\end{proposition}

\begin{proof} Set $\L = \{\ell_1 < \ldots \ell_p\}$ and $\R=\{r_1 < \ldots < r_q\}$ such that $(\L,\R)$ is the elementary bipartition of $\{1,\ldots,n+1\}$ such that $\lambda = \oplamb(\L,\R)$. We can check that, in $\AR(c)$, $(\ell_{i+1},r_j)$ follows directly $(\ell_i,r_j)$, and $(\ell_i, r_{j+1})$ follows directly $(\ell_i, r_j)$. Then, the orientation we must choose to realize the Hillman--Grassl correspondence coincides with the one given by $\AR(c)$. Therefore, we get the desired result.
\end{proof}

Now we highlight the link with the work of \cite{GPT19}

\begin{proposition} \label{prop:GPTandnewGRSK}
	Let $n \geqslant 1$ and $\lambda \in \Hk_n$. Consider a Coxeter element $c \in \mathfrak{S}_{n+1}$. Then, for any filling $f$ of $\lambda$, \[\Phi_{\opc(\lambda)}^{-1} \circ \GRSK_{\lambda,c}(f) = \GenJF(\opE_{Q(c)}(f)).\] 
\end{proposition}

Then, the following corollary holds.

\begin{corollary} \label{cor:ScrambledRSK}
	Let $n \geqslant 1$ and $\lambda \in \Hk_n$. Consider a Coxeter element $c \in \mathfrak{S}_{n+1}$. If $\opc(\lambda) = (1,2,\ldots,m,n+1,n,\ldots,m+1)$ for some $m \in \{1,\ldots, n\}$, then \[\GRSK_{\lambda,c} = \RSK_{m,c}\] where $\RSK_{m,c}$ corresponds to Scrambled RSKs.
\end{corollary}

Note that choosing a filling of such a $\lambda$ corresponds to a representation, up to isomorphism, of $\mathscr{C}_{Q(c),m}$.

\subsection{Enumerative properties}
\label{ss:Combidentity}

\begin{lemma} \label{lem:sumofdiagonals}
	Let $n \geqslant 1$, and $\lambda \in \Hk_n$. Consider a Coxeter elements $c \in \mathfrak{S}_{n+1}$. Then, for any filling $f$ of $\lambda$, and for any $k \in \{1,\ldots,n\}$, we have \[\sum_{\varepsilon=1}^{\delta_k} \GRSK_{\lambda,c}(f)(\ldiag{k,\varepsilon}_\lambda) = \sum_{b \in \square_k(\lambda)} f(b) .\]
\end{lemma}

\begin{proof}
	This is a direct consequence of \cref{GKprop} and the definition of the extended generalization. Precisely, by setting $\L = \L_c \cup \{1\}$ and $\R = \R_c \cup \{n+1\}$, 
	\[\sum_{\varepsilon=1}^{\delta_k} \GRSK_{\lambda,c}(f)(\ldiag{k,\epsilon}_\lambda) = \sum_{(\ell,r) \in \L \times \R, \ell \leqslant k < r} \rep_{\lambda,c}(f)(\ell,r) = \sum_{b \in \square_k(\lambda)} f(b)\]
\end{proof}

\begin{remark} \label{rem:dimsum} From the previous proof, one can notice that, given $\lambda \in \Hk_n$ and a Coxeter element $c \in \mathfrak{S}_{n+1}$, we have 
	\[\dim(\opE_{Q(c)}(f)_k) = \sum_{b \in \square_k(\lambda)} f(b).\]
	This value does not depend on $c$.
\end{remark}

The following corollary extends the product formula of the trace generating function of reverse plane partitions of minuscule posets. Proctor originally established this version \cite{P84}, and Garver, Patrias, and Thomas proved it differently \cite{GPT19}.

\begin{corollary} Let $n \geqslant 1$, and $\lambda \in \Hk_n$. For any $A_n$ type quiver $Q$, we have \[\rho_\lambda(x_1,\ldots,x_n) = \sum_{E\in \Cat_Q(\opc(\lambda))} \prod_{\ell=1}^n x_\ell^{\dim(E_\ell)} = \prod_{M \in \Cat_Q(\opc(\lambda)), \text{ indec.}} \dfrac{1}{1 - \prod_{\ell=1}^n x_\ell^{\dim(M_\ell)}}, \] where the sum is over the isomorphism classes of representations in $\Cat_Q(\opc(\lambda))$, and the product is over isomorphism classes of indecomposable representations in $\Cat_Q(\opc(\lambda))$
\end{corollary}

\begin{proof}
	Let $Q$ be an $A_n$ type quiver. Consider the Coxeter element $c \in \mathfrak{S}_{n+1}$ such that $Q(c) = Q$. By \cref{lem:sumofdiagonals} and \cref{rem:dimsum}, we have \[\rho_\lambda(x_1,\ldots,x_n) = \sum_{f \in \RPP(\lambda)} \prod_{\ell=1}^n x_\ell^{\sum_{\varepsilon =1}^{\delta_l} f(\ldiag{\ell,\varepsilon}_\lambda)} = \sum_{f \in \RPP(\lambda)} \prod_{\ell=1}^n x_\ell^{\dim(E_{Q(c)}(f)_\ell))}.\]
	
	So \[\rho_\lambda(x_1,\ldots,x_n) = \sum_{E \in \Cat_Q(\opc(\lambda))} \prod_{\ell=1}^n x_\ell^{\dim(E_\ell)}.\]
	
	Moreover, each box $b \in \Fer(\lambda)$ corresponds bijectively to an indecomposable representation, using the $\opc(\lambda)-$coordinates and $\opE_Q \circ \crep_{\lambda, c}$. Set, for any $b \in \Fer(\lambda)$, $M_b$ its associated indecomposable representation in $\rep_\mathbb{K}(Q)$, up to isomorphism. We have that:
	\begin{enumerate}[label=$\bullet$, itemsep=1mm]
		\item $b \in \square_\ell(\lambda)$ if and only if $\dim((M_b)_\ell) = 1$, and,
		\item $b \notin \square_\ell(\lambda)$ if and only if $\dim((M_b)_\ell) = 0$.
	\end{enumerate} Therefore, \[\omega_{\lambda,b}(x_1,\ldots,x_n) = \prod_{\ell =1}^n x_\ell^{\dim((M_b)_\ell)}.\]
	
	We get the desired result by \cref{cor:tracegenformula}.
\end{proof}

\begin{remark}
	Note that this identity (the second equality) is a consequence of our extended generalization $\GRSK_{\lambda,c}$ where the Coxeter element $c \in \mathfrak{S}_{n+1}$ corresponds to the choice of an $A_n$ type quiver. The product formula comes from the fact that any representation $E \in \rep_\mathbb{K}(Q)$ can be decomposed, in a unique way, as a sum of indecomposable representations (See \cref{ss:JRandCJR}).
\end{remark}

	\addtocontents{toc}{\protect\setcounter{tocdepth}{1}}
	
	\section{To go further}
	\label{sec:Further}
	
In this section, we suggest some research directions that could follow this work.

\subsection*{A realization of $\GRSK_{\lambda,c}$ via local transformations}

\cite{H14} gives a complete realization of $\GRSK_\lambda$ via local transformations, for any integer partitions. In quiver representation settings, \cite{GPT19} give a realization of $\RSK_{m,c}$ via local transformations (``toggles"), for any $n \in \mathbb{N}^*$, and any Coxeter element $c \in \mathfrak{S}_{n+1}$. In parallel, \cite{DNV21} give a realization of $\RSK_{m,c}$ (mentioned as \emph{Scrambled RSK}) via the same ``toggles" but described using combinatorial terms.

Thinking about the crucial link with quiver representations and the local description given by \cite[Section 4.3]{GPT19}, one could think about considering a sequence of ``toggles", translated from operations preserving canonical Jordan recoverability (see \cref{ss:CJRpreserving}), based on the linear order on transpositions of $\mathfrak{S}_{n+1}$, compatible with (the opposite of) the Auslander--Reiten quiver. 

We hope to elaborate such a realization of $\GRSK_{\lambda,c}$ in the near future.

\subsection*{Dynkin type (or other) variations of $\GRSK_{\lambda,c}$}

Instead of considering Coxeter elements of the symmetric group, one could ask if it is possible to consider Coxeter elements of any Weyl group. We see at least two ways to think about that.

The first way is to adapt the setting to the Weyl group we are considering. For instance, if we work with the signed symmetric group, a $B$ type Weyl group, we could think about type $B$ RSK, seen as the domino correspondence. We refer the reader to the work of Garfinkle \cite{G90,G92} Bonafé, Geck, Iancu, and Lam \cite{BGIL10} for more details. In the type $D$ case, Gern did a few studies in his PhD thesis \cite{G13} in terms of Kazhdan--Lusztig polynomials.

The second way could be to work in the quiver representation setting and to determine all the canonically Jordan recoverable subcategories of any type $D$ and type $E$ quivers. \cite[Theorem 1.3]{GPT19} gives us a type $D$ and  a type $E$ versions of the RSK correspondence.

From the perspective of quiver representation theory, we can go deeper into it and work with well-behaved and well-known bounded quivers. For instance, for \emph{string quivers}, we have a complete classification of the isomorphism classes of indecomposable representations \cite{BR87} and of the morphisms between them \cite{K88,CB89}. Some work was already done, for a subfamily of string quivers, called \emph{gentle quivers}, on canonically Jordan recoverable subcategories in this domain \cite{D22}, and we hope to characterize all of them. One interesting combinatorial outcome could be the construction of an ``RSK correspondence via string algebras".

	\addtocontents{toc}{\protect\setcounter{tocdepth}{8}}
	
	\section*{Acknowledgements}
	
	The author acknowledges the \emph{CHARMS program grant (ANR-19-CE40-0017-02)} and the \emph{Engineering and Physical Sciences Research Council (EP/W007509/1)} for their partial funding support. He thanks the selection committee of the 36th edition of the FPSAC Conference (Bochum, 2024) for its insightful comments on my extended abstract \cite{DFPSAC24}, accepted for the poster session. 
	
	He is grateful to Ben Adenbaum, Emily Gunawan, Florent Hivert, Yann Palu, GaYee Park, and Michael Schoonheere for their interest and discussions on this project. He especially thanks Phillippe Nadeau for his advice, which allowed him to simplify some combinatorial proofs. Finally, he acknowledges Hugh Thomas for his advice and comments on this work.
	\bibliography{References3}

@book{ASS06, place={Cambridge}, series={London Mathematical Society Student Texts}, title={Elements of the Representation Theory of Associative Algebras: Techniques of Representation Theory}, volume={1}, DOI={10.1017/CBO9780511614309}, publisher={Cambridge University Press}, author={Assem, Ibrahim and Skowro{\`{n}}ski, Andrzej and Simson, Daniel}, year={2006}, collection={London Mathematical Society Student Texts}}

@book{F96, place={Cambridge}, series={London Mathematical Society Student Texts}, title={Young Tableaux: With Applications to Representation Theory and Geometry}, DOI={10.1017/CBO9780511626241}, publisher={Cambridge University Press}, author={Fulton, William}, year={1996}, collection={London Mathematical Society Student Texts}}

@book{St99,
	author = {Stanley, Richard P.},
	isbn = {9781139811729 113981172X},
	publisher = {Cambridge University Press},
	refid = {821883332},
	title = {Enumerative Combinatorics},
	url = {http://math.mit.edu/~rstan/ec/},
	volume = 2,
	year = 1999
}

@article{P84,
	title = {Bruhat Lattices, Plane Partition Generating Functions, and Minuscule Representations},
	journal = {European Journal of Combinatorics},
	volume = {5},
	number = {4},
	pages = {331-350},
	year = {1984},
	issn = {0195-6698},
	doi = {https://doi.org/10.1016/S0195-6698(84)80037-2},
	url = {https://www.sciencedirect.com/science/article/pii/S0195669884800372},
	author = {Proctor, Robert A.}
}

@article{D22,
	title = {Jordan recoverability of some subcategories of modules over gentle algebras},
	journal = {Journal of Pure and Applied Algebra},
	volume = {228},
	number = {3},
	pages = {107446},
	year = {2024},
	issn = {0022-4049},
	doi = {https://doi.org/10.1016/j.jpaa.2023.107446},
	url = {https://www.sciencedirect.com/science/article/pii/S0022404923001299},
	author = {Benjamin Dequêne}
}

@article{Dauv20,
	author = {Dauvergne, Duncan},
	title = {{Hidden invariance of last passage percolation and directed polymers}},
	volume = {50},
	journal = {The Annals of Probability},
	number = {1},
	publisher = {Institute of Mathematical Statistics},
	pages = {18 -- 60},
	year = {2022},
	doi = {10.1214/21-AOP1527},
	URL = {https://doi.org/10.1214/21-AOP1527}}

@article{Ga81Ma,
	author = {Gansner, Emden},
	year = {1981},
	month = {02},
	pages = {295-315},
	title = {Matrix correspondences of plane partitions},
	volume = {92},
	journal = {Pacific Journal of Mathematics},
	doi = {10.2140/pjm.1981.92.295}
}

@book{St72,
	author = {Stanley, Richard P.},
	address = {Providence, R.I},
	booktitle = {Ordered structures and partitions},
	publisher = {American Mathematical Society},
	series = {Memoirs of the American Mathematical Society ; no. 119},
	isbn = {9780821818190},
	language = {eng},
	lccn = {52042839},
	title = {Ordered structures and partitions },
	year = {1972},
	
}

@article{Ga81Hi,
	title = "{The Hillman-Grassl correspondence and the enumeration of reverse plane partitions}",
	journal = {Journal of Combinatorial Theory, Series A},
	volume = {30},
	number = {1},
	pages = {71-89},
	year = {1981},
	issn = {0097-3165},
	doi = {https://doi.org/10.1016/0097-3165(81)90041-8},
	url = {https://www.sciencedirect.com/science/article/pii/0097316581900418},
	author = {Emden R Gansner}
}

@article{GPT19,
	author={Garver, Alexander
	and Patrias, Rebecca
	and Thomas, Hugh},
	title={Minuscule reverse plane partitions via quiver representations},
	journal={Selecta Mathematica},
	year={2023},
	month={Apr},
	day={26},
	volume={29},
	number={3},
	pages={37},
	issn={1420-9020},
	doi={10.1007/s00029-023-00831-4},
	url={https://doi.org/10.1007/s00029-023-00831-4}
}

@article{Kr06,
	title = {Growth diagrams, and increasing and decreasing chains in fillings of Ferrers shapes},
	journal = {Advances in Applied Mathematics},
	volume = {37},
	number = {3},
	pages = {404-431},
	year = {2006},
	note = {Special Issue in honor of Amitai Regev on his 65th Birthday},
	issn = {0196-8858},
	doi = {https://doi.org/10.1016/j.aam.2005.12.006},
	url = {https://www.sciencedirect.com/science/article/pii/S0196885806000546},
	author = {Krattenthaler, Christian}
}

@article{HG76,
	author  = {Hillman, Abraham P. and Grassl, Richard M.},
	title   = {Reverse Plane Partitions and Tableau Hook Numbers},
	journal = {Journal of Combinatorial Theory, Series A},
	volume  = {21},
	number  = {2},
	pages   = {216--221},
	year    = {1976},
	doi     = {10.1016/0097-3165(76)90065-0}
}

@article{B72,
	author  = {Burge, William H.},
	title   = {Four Correspondences Between Graphs and Generalized Young Tableaux},
	journal = {Journal of Combinatorial Theory, Series A},
	volume  = {17},
	number  = {1},
	pages   = {12--30},
	year    = {1974},
	doi     = {10.1016/0097-3165(74)90024-7}
}

@article{K70,
	author  = {Knuth, Donald E.},
	title   = {Permutations, Matrices, and Generalized Young Tableaux},
	journal = {Pacific Journal of Mathematics},
	volume  = {34},
	number  = {3},
	pages   = {709--727},
	year    = {1970},
	doi     = {10.2140/pjm.1970.34.709}
}

@article{Deq23,
	author       = {Dequ{\^e}ne, Benjamin},
	title        = {Canonically Jordan recoverable categories for modules over the path algebra of {$A_n$} type quivers},
	journal      = {arXiv preprint},
	year         = {2023},
	eprint       = {2308.16626},
	archivePrefix= {arXiv},
	primaryClass = {math.RT},
	doi          = {10.48550/arXiv.2308.16626}
}

@article{GK76,
	title = {The structure of {S}perner $k$-families},
	journal = {Journal of Combinatorial Theory, Series A},
	volume = {20},
	number = {1},
	pages = {41-68},
	year = {1976},
	issn = {0097-3165},
	doi = {https://doi.org/10.1016/0097-3165(76)90077-7},
	url = {https://www.sciencedirect.com/science/article/pii/0097316576900777},
	author = {Curtis Greene and Daniel J Kleitman}
}

@book{Sch14,
	author    = {Schiffler, Ralf},
	title     = {Quiver {R}epresentations},
	series    = {CMS Books in Mathematics},
	publisher = {Springer},
	address   = {Cham},
	year      = {2014},
	isbn      = {978-3-319-09203-2},
	doi       = {10.1007/978-3-319-09204-9}
}

@Article{R38,
	author={Robinson, Giles de B.},
	title={On the Representations of the Symmetric Group},
	journal={American Journal of Mathematics},
	year={1938},
	month={2024/07/02/},
	publisher={Johns Hopkins University Press},
	volume={60},
	number={3},
	pages={745-760},
	note={Full publication date: Jul., 1938},
	doi={10.2307/2371609},
	url={http://www.jstor.org/stable/2371609}
}

@article{S61, 
	author={Schensted, Craige},
	journal={Canadian Journal of Mathematics},  
	title={Longest Increasing and Decreasing Subsequences}, 
	volume={13}, 
	pages={179–191},
	year={1961},
	DOI={10.4153/CJM-1961-015-3}, 
}

@InCollection{ES87,
	author = {Erd{\H{o}}s, Paul and Szekeres, George},
	editor = {Gessel, Ira and Rota, Gian-Carlo},
	title = {A Combinatorial Problem in Geometry},
	booktitle = {Classic Papers in Combinatorics},
	publisher = {Birkh{\"a}user Boston},
	address = {Boston, MA},
	year = {1987},
	pages = {49--56},
	isbn = {978-0-8176-4842-8},
	doi = {10.1007/978-0-8176-4842-8_3},
	url = {https://doi.org/10.1007/978-0-8176-4842-8_3}
}

@book{S13,
	title={The Symmetric Group: Representations, Combinatorial Algorithms, and Symmetric Functions},
	author={Sagan, Bruce E.},
	isbn={9781475768046},
	lccn={00040042},
	series={Graduate Texts in Mathematics},
	url={https://books.google.fr/books?id=Y6vTBwAAQBAJ},
	year={2013},
	publisher={Springer New York}
}

@Article{S97,
	author = {Sch\"utzenberger, Marcel-Paul},
	title = {Pour le mono{\"\i}de plaxique},
	journal = {Math\'ematiques informatique et sciences humaines},
	pages = {5--10},
	publisher = {Ecole des hautes-\'etudes en sciences sociales},
	volume = {140},
	year = {1997},
	mrnumber = {1627563},
	zbl = {0948.20040},
	language = {fr},
	url = {http://www.numdam.org/item/MSH_1997__140__5_0/}
}

@Article{LS81,
	author = {Lascoux, Alain and Schützenberger, Marcel-P.},
	title = {Le monoïde plaxique},
	year = {1981},
	volume = {109},
	pages={129–-156},
	journal ={Noncommutative structures in algebra and geometric combinatorics (Naples, 1978)},
	publisher = {Ricercia Sci.},
}

@incollection{V77,
	author    = {Viennot, G{\'e}rard},
	editor    = {Foata, Dominique},
	title     = {Une forme g{\'e}om{\'e}trique de la correspondance de {R}obinson--{S}chensted},
	booktitle = {Combinatoire et Repr{\'e}sentation du Groupe Sym{\'e}trique},
	publisher = {Springer},
	address   = {Berlin},
	year      = {1977},
	pages      = {29--58},
	doi        = {10.1007/BFb0069178}
}

@article{AF22,
	author = {Aigner, Florian and Frieden, Gabriel},
	title = "{qRSt: A Probabilistic Robinson–Schensted Correspondence for Macdonald Polynomials}",
	journal = {International Mathematics Research Notices},
	volume = {2022},
	number = {17},
	pages = {13505-13568},
	year = {2021},
	month = {05},
	issn = {1073-7928},
	doi = {10.1093/imrn/rnab083},
	url = {https://doi.org/10.1093/imrn/rnab083},
	eprint = {https://academic.oup.com/imrn/article-pdf/2022/17/13505/45617599/rnab083.pdf},
}

@article{P01,
	author = {Pak, Igor},
	journal = {Séminaire Lotharingien de Combinatoire [electronic only]},
	language = {eng},
	pages = {B46f, 13 p.},
	publisher = {Universität Wien, Fakultät für Mathematik},
	title = {Hook length formula and geometric combinatorics.},
	url = {http://eudml.org/doc/121696},
	volume = {46},
	year = {2001},
}

@article{GRB23,
	author = {Gillespie, Maria and Reimer-Berg, Andrew},
	title = {A {Generalized} {RSK} for {Enumerating} {Linear} {Series} on $n$-pointed {Curves}},
	journal = {Algebraic Combinatorics},
	pages = {1--16},
	publisher = {The Combinatorics Consortium},
	volume = {6},
	number = {1},
	year = {2023},
	doi = {10.5802/alco.250},
	language = {en},
	url = {https://alco.centre-mersenne.org/articles/10.5802/alco.250/}
}

@article{R07,
	author  = {Reading, Nathan},
	title   = {Clusters, {C}oxeter-Sortable Elements and Noncrossing Partitions},
	journal = {Transactions of the American Mathematical Society},
	volume  = {359},
	number  = {12},
	pages   = {5931--5958},
	year    = {2007},
	doi     = {10.1090/S0002-9947-07-04264-0}
}

@article{DNV21,
	author = {Dauvergne, Duncan and Nica, Mihai and  Vir{\'a}g, B{\'a}lint},
	title = {{RSK in last passage percolation: a unified approach}},
	volume = {19},
	journal = {Probability Surveys},
	number = {none},
	publisher = {Institute of Mathematical Statistics and Bernoulli Society},
	pages = {65 -- 112},
	year = {2022},
	doi = {10.1214/22-PS4},
	URL = {https://doi.org/10.1214/22-PS4}
}

@misc{H14,
	author = {Sam Hopkins},
	title = {RSK via local transformations},
	year = {2014},
	note = {\url{https://www.samuelfhopkins.com/docs/rsk.pdf}}
}

@InCollection{BGIL10,
	author = {Bonnaf{\'e}, C{\'e}dric and Geck, Meinolf and Iancu, Lacrimioara and Lam, Thomas},
	editor = {Gyoja, Akihiko and Nakajima, Hiraku and Shinoda, Ken-ichi and Shoji, Toshiaki and Tanisaki, Toshiyuki},
	title = {On Domino Insertion and {K}azhdan--{L}usztig Cells in Type {$B_n$}},
	booktitle = {Representation Theory of Algebraic Groups and Quantum Groups},
	publisher = {Birkh{\"a}user Boston},
	address = {Boston},
	year = {2010},
	pages = {33--54},
	isbn = {978-0-8176-4697-4},
	doi = {10.1007/978-0-8176-4697-4_3},
	url = {https://doi.org/10.1007/978-0-8176-4697-4_3}
}

@article{G90,
	author = {Garfinkle, Devra},
	title = {On the classification of primitive ideals for complex classical {Lie} algebras, {I}},
	journal = {Compositio Mathematica},
	pages = {135--169},
	publisher = {Kluwer Academic Publishers},
	volume = {75},
	number = {2},
	year = {1990},
	mrnumber = {1065203},
	zbl = {0737.17003},
	language = {en},
	url = {http://www.numdam.org/item/CM_1990__75_2_135_0/}
}

@article{G92,
	author = {Garfinkle, Devra},
	title = {On the classification of primitive ideals for complex classical {Lie} algebras, {II}},
	journal = {Compositio Mathematica},
	pages = {307--336},
	publisher = {Kluwer Academic Publishers},
	volume = {81},
	number = {3},
	year = {1992},
	mrnumber = {1149172},
	zbl = {0762.17007},
	language = {en},
	url = {http://www.numdam.org/item/CM_1992__81_3_307_0/}
}

@phdthesis{G13,
	title        = "{Leading Coefficients of Kazhdan--Lusztig Polynomials in Type $D$}",
	author       = {Tyson C. Gern},
	year         = {2013},
	school       = {University of Colorado Boulder},
	type         = {PhD thesis}
}

@inproceedings{DFPSAC24,
	author = {Dequ{\^e}ne, Benjamin},
	title = {An extended generalization of RSK via the combinatorics of type {A} quiver representations},
	booktitle = {S{\'e}minaire Lotharingien de Combinatoire: Proceedings of the 36th Conference on Formal Power Series and Algebraic Combinatorics},
	year = {2024},
	month = jul,
	address = {Bochum, Germany},
	organization = {FPSAC},
	volume = {91B (76)},
	url = {https://drive.google.com/file/d/1ozqqLae3_WqWzDbKEBwI6FF94tmDKsXc/view?usp=sharing}
}

@article{BR87,
	author = {Butler, Michael Charles Richard and Ringel, Claus Michael},
	title = {Auslander-{R}eiten sequences with few middle terms and applications to string algebras},
	journal = {Communications in Algebra},
	volume = {15},
	number = {1-2},
	pages = {145-179},
	year  = {1987},
	publisher = {Taylor & Francis}
}

@article{CB89,
	title = {Maps between representations of zero-relation algebras},
	journal = {Journal of Algebra},
	volume = {126},
	number = {2},
	pages = {259-263},
	year = {1989},
	issn = {0021-8693},
	doi = {https://doi.org/10.1016/0021-8693(89)90304-9},
	author = {Crawley-Boevey, William W.}
}

@article{K88,
	author  = {Krause, Henning},
	title   = {Maps Between Tree and Band Modules},
	journal = {Journal of Algebra},
	volume  = {137},
	number  = {1},
	pages   = {186--194},
	year    = {1991},
	doi     = {10.1016/0021-8693(91)90088-P}
}

@article{BGP73,
	author  = {Bernstein, I. N. and Gel'fand, I. M. and Ponomarev, V. A.},
	title   = {Coxeter Functors and Gabriel's Theorem},
	journal = {Russian Mathematical Surveys},
	volume  = {28},
	number  = {2},
	year    = {1973},
	pages   = {17--32},
	doi     = {10.1070/RM1973v028n02ABEH001526}
}
	\bibliographystyle{alpha}
\end{document}